\documentclass[reqno, a4paper, 10pt]{amsart}

\usepackage[margin=1in]{geometry}
\usepackage[algo2e, linesnumbered, noline, noend, ruled]{algorithm2e}
\usepackage[foot]{amsaddr}
\usepackage{amsfonts}
\usepackage{amsmath}
\usepackage{amssymb}
\usepackage{amsthm}
\usepackage[english]{babel}
\usepackage{booktabs}
\usepackage{braket}
\usepackage{float}
\usepackage[T1]{fontenc}
\usepackage{graphicx}
\usepackage[utf8]{inputenc}
\usepackage{nicefrac}
\usepackage[retain-explicit-plus, retain-zero-exponent, per-mode=symbol, output-exponent-marker=\text{e}, group-digits=false]{siunitx}
\usepackage[table]{xcolor}
\usepackage{enumitem}
\usepackage{subcaption}

\usepackage{cite}
\usepackage{hyperref}

	\newtheorem{Theorem}{Theorem}[section]
	\newtheorem{Proposition}[Theorem]{Proposition}

	\theoremstyle{definition}
	\newtheorem{Definition}[Theorem]{Definition}
	\newtheorem{Remark}[Theorem]{Remark}
	\newtheorem*{Remark*}{Remark}
	
	\setenumerate{label=(\arabic*), ref=(\arabic*)}
	
	\newcommand{\R}{\mathbb{R}}
	
	\newcommand{\norm}[2]{\left\lvert \left\lvert #1 \right\rvert \right\rvert_{#2}}
	
	\newcommand{\dmeas}[1]{\ \mathrm{d}#1}

	\makeatletter
	\DeclareOldFontCommand{\rm}{\normalfont\rmfamily}{\mathrm}
	\DeclareOldFontCommand{\sf}{\normalfont\sffamily}{\mathsf}
	\DeclareOldFontCommand{\tt}{\normalfont\ttfamily}{\mathtt}
	\DeclareOldFontCommand{\bf}{\normalfont\bfseries}{\mathbf}
	\DeclareOldFontCommand{\it}{\normalfont\itshape}{\mathit}
	\DeclareOldFontCommand{\sl}{\normalfont\slshape}{\@nomath\sl}
	\DeclareOldFontCommand{\sc}{\normalfont\scshape}{\@nomath\sc}
	\makeatother
	
	\numberwithin{equation}{section}
	
	\allowdisplaybreaks

	\newcommand{\tablegray}{gray!25}
	
	\newcommand{\subin}{^\mathrm{in}}
	\newcommand{\subwall}{^\mathrm{wall}}
	\newcommand{\subout}{^\mathrm{out}}
	
	\newcommand{\subfix}{^\mathrm{fix}}
	\newcommand{\subdef}{^\mathrm{def}}
	\newcommand{\subobs}{^\mathrm{obs}}

	\newcommand{\reynolds}{\mathrm{Re}}
	\newcommand{\conductivity}{\kappa}
	
	\newcommand{\normal}{n}

	\newcommand{\velocity}{u}
	\newcommand{\pressure}{p}

	\newcommand{\laplace}{\Delta}
	\newcommand{\grad}{\nabla}
	\newcommand{\divergence}[1]{\mathrm{div}\left(#1\right)}
	\newcommand{\integral}[1]{\int_{#1}}

	\newcommand{\costfunction}{\mathcal{J}}
	\newcommand{\reducedcostfunction}{J}
	\newcommand{\state}{u}
	\newcommand{\adjoint}{p}
	\newcommand{\test}{v}
	\newcommand{\holdall}{D}

	\newcommand{\admissiblegeom}{\mathcal{A}}
	\newcommand{\vectorfield}{\mathcal{V}}
	\newcommand{\flow}{\Phi_t}
	\newcommand{\shapedistro}{g}
	\newcommand{\gradientdefo}{\mathcal{G}}
	\newcommand{\searchdefo}{\mathcal{D}}
	\newcommand{\scalarproduct}{g^S}
	\newcommand{\stekpoinc}{S^p}
	\newcommand{\rieshapegrad}{\gamma}
	\newcommand{\geodesic}{\gamma}

	\newcommand{\retraction}{\mathrm{R}}
	\newcommand{\numretraction}{\tilde{\mathrm{R}}}
	\newcommand{\vectortransport}{\mathcal{T}}
	\newcommand{\numtransport}{\tilde{\mathcal{T}}}
	\newcommand{\searchdirection}{d}
	\newcommand{\increment}{\eta}
	\newcommand{\iidx}[1]{_{#1}}
	\newcommand{\shapemanifold}{B_e}
	\newcommand{\tangentspace}{T}
	\newcommand{\graddiff}{\mathcal{Y}}
	
	\newcommand{\lamefirst}{\lambda_{\mathrm{elas}}}
	\newcommand{\lamesecond}{\mu_{\mathrm{elas}}}
	\newcommand{\damping}{\delta_{\mathrm{elas}}}
	\newcommand{\symgrad}{\varepsilon}
	
	\newcommand{\armijoparam}{\sigma}
	\newcommand{\stepparam}{\omega}
	
	\newcommand{\cgiter}{k_{\mathrm{cg}}}
	\newcommand{\cgtol}{\varepsilon_\mathrm{cg}}
	
	\newcommand{\transposed}{^\top}
	
\begin{document}

{\small \begin{center}
		This is a post-peer-review, pre-copyedit version of an article published in SIAM Journal on Optimization. 
		
		The final version is available online at \url{https://doi.org/10.1137/20M1367738}.
	\end{center} }

\title[NCG Methods for PDE constrained Shape Optimization]{Nonlinear Conjugate Gradient Methods for PDE Constrained Shape Optimization Based on Steklov-Poincar\'e-Type Metrics}
\author{Sebastian Blauth$^{1,2}$}
\address{$^1$ Fraunhofer ITWM, Kaiserslautern, Germany}
\address{$^2$ TU Kaiserslautern, Kaiserslautern, Germany}
\email{\href{mailto:sebastian.blauth@itwm.fraunhofer.de}{sebastian.blauth@itwm.fraunhofer.de}}


\begin{abstract}
	Shape optimization based on shape calculus has received a lot of attention in recent years, particularly regarding the development, analysis, and modification of efficient optimization algorithms. In this paper we propose and investigate nonlinear conjugate gradient methods based on Steklov-Poincar\'e-type metrics for the solution of shape optimization problems constrained by partial differential equations. We embed these methods into a general algorithmic framework for gradient-based shape optimization methods and discuss the numerical discretization of the algorithms. We numerically compare the proposed nonlinear conjugate gradient methods to the already established gradient descent and limited memory BFGS methods for shape optimization on several benchmark problems. The results show that the proposed nonlinear conjugate gradient methods perform well in practice and that they are an efficient and attractive addition to already established gradient-based shape optimization algorithms.
	
	
	
	
\end{abstract}

\maketitle

\vspace{-0.5cm}
\section{Introduction}
\label{sec:introduction}

Shape optimization problems constrained by partial differential equations (PDEs) and their solution based on shape calculus have received a lot of attention in recent years. They are used for many industrial applications, e.g., to optimize electric motors \cite{Gangl2015Shape, Gangl2018Sensitivity}, polymer spin packs \cite{Hohmann2019Shape, Leithaeuser2018Shape, Leithaeuser2018Designing}, aircrafts and automobiles \cite{Othmer2014Adjoint, Gauger2012Non}, or microchannel heat exchangers \cite{Blauth2019Model, Blauth2020Shape}. Alternatively, such shape optimization techniques are also used to solve inverse problems, e.g., in image segmentation \cite{Hintermueller2003second, Dogan2008variational} or electrical impedance tomography \cite{Hintermueller2008Electrical, Laurain2016Distributed}. The development, analysis, and modification of algorithms for the efficient solution of PDE constrained shape optimization problems has also attracted a lot of interest in recent literature, e.g., in \cite{Schulz2014Riemannian}, where a Riemannian view on shape optimization and corresponding shape Newton methods are analyzed, in \cite{Schulz2016Efficient, Schulz2016Computational}, where quasi-Newton methods based on Steklov-Poincar\'e-type metrics are proposed and investigated numerically, in \cite{Etling2020First}, where restricted mesh deformations are considered to obtain better meshes and to avoid remeshing, or in \cite{Hiptmair2015Comparison}, where the numerical approximation of shape derivatives is investigated in the finite element context.

In this paper, we propose nonlinear conjugate gradient (NCG) methods for shape optimization problems based on the Steklov-Poincar\'e-type metrics introduced in \cite{Schulz2016Efficient}. These methods are of particular interest for large-scale optimization problems since they have comparatively low memory requirements, and yet are very efficient. 
In particular, NCG methods usually outperform the widely used gradient descent method and often exhibit superlinear convergence behavior. 
To the best of our knowledge, such nonlinear CG methods for PDE constrained shape optimization have not been investigated in the literature so far.

We present the NCG methods in the context of a general algorithmic framework for gradient-based shape optimization methods, into which a gradient descent method, limited memory BFGS (L-BFGS) methods from \cite{Schulz2016Efficient}, and the novel NCG methods are embedded. In particular, in this paper we consider five popular NCG methods, namely the Fletcher-Reeves \cite{Fletcher1964Function}, Polak-Ribi\`ere \cite{Polak1969Note, Polyak1971conjugate}, Hestenes-Stiefel \cite{Hestenes1952Methods}, Dai-Yuan \cite{Dai1999nonlinear}, and Hager-Zhang \cite{Hager2005new} NCG variants. As our algorithmic framework is based on the Steklov-Poincar\'e-type metrics from \cite{Schulz2016Efficient}, it is well-suited for numerical discretization and leads to little computational overhead. 

We investigate the numerical performance of the NCG methods on four benchmark shape optimization problems: A shape optimization problem constrained by a Poisson equation from \cite{Etling2020First}, a shape identification problem in electrical impedance tomography based on the ones considered in \cite{Laurain2016Distributed, Hintermueller2008Electrical, Schulz2016Efficient}, the shape optimization of an obstacle in Stokes flow from \cite{Schulz2016Computational}, and the shape optimization of a pipe with Navier-Stokes flow from \cite{Schmidt2010Efficient}. For each of these problems, we compare the NCG methods with the gradient descent and L-BFGS methods to evaluate their performance. 
The obtained results show that the proposed NCG methods perform very well in practice. In particular, the NCG methods always significantly outperform the gradient descent method, but perform slightly worse compared to the L-BFGS methods, which is to be expected from their finite-dimensional counterparts. However, the NCG methods require considerably less memory than the L-BFGS methods, which makes them an attractive addition to gradient-based shape optimization algorithms.



This paper is structured as follows. In Section~\ref{sec:preliminaries}, we recall NCG methods for finite-dimensional optimization problems as well as shape calculus and the Steklov-Poincar\'e-type metrics for Riemannian shape optimization from \cite{Schulz2016Efficient}. A general algorithmic framework for shape optimization as well as the NCG methods embedded into it are presented in Section~\ref{sec:algorithms}. Finally, in Section~\ref{sec:numerics}, we numerically compare the NCG methods with the gradient descent and L-BFGS methods on the four shape optimization problems described previously and investigate their performance.

\section{Preliminaries}
\label{sec:preliminaries}

First, we recapitulate NCG methods for finite-dimensional nonlinear optimization problems. Afterwards, we recall shape calculus as well as the Riemannian view on it from \cite{Schulz2014Riemannian} and the Steklov-Poincar\'e-type metrics from \cite{Schulz2016Efficient}.

%

\subsection{NCG Methods for Finite Dimensional Optimization Problems}
\label{ssec:ncg_finite}

Let us briefly recall NCG methods for finite-dimensional optimization problems. A classical nonlinear optimization problem in $\R^n$ is given by
\begin{equation*}
	\min_x\ f(x) \qquad \text{ s.t. } x\in \R^n,
\end{equation*}
where $f\in C^1(\R^n; \R)$. Starting from an initial guess $x\iidx{0}$, NCG methods attempt to solve this problem through the iteration
\begin{equation*}
	x\iidx{k+1} = x\iidx{k} + \alpha\iidx{k} \searchdirection\iidx{k}.
\end{equation*}
Here, $\alpha\iidx{k} > 0$ is a step size that is computed, e.g., by means of a backtracking line search, and $\searchdirection\iidx{k}$ is the search direction defined by
\begin{equation*}
	\searchdirection\iidx{k} = - g\iidx{k} + \beta\iidx{k} \searchdirection\iidx{k-1}, \qquad \text{ with } \qquad d\iidx{0} = - g\iidx{0},
\end{equation*}
where $g\iidx{k} = \grad f(x\iidx{k})$ and $\beta\iidx{k}$ is the update parameter for the NCG methods. In the literature, there are several update formulas available, each leading to a slightly different NCG method. For the description of $\beta\iidx{k}$, we denote by $\norm{\cdot}{}$ the Euclidean norm in $\R^n$ and write $y\iidx{k} = g\iidx{k+1} - g\iidx{k}$. In this paper, we consider the following five popular variants given by
\begin{alignat*}{2}
	\beta^{\text{FR}}\iidx{k} &= \frac{\norm{g\iidx{k}}{}^2}{\norm{g\iidx{k-1}}{}^2} \qquad &&\text{ (Fletcher and Reeves \cite{Fletcher1964Function})},\\
	\beta^{\text{PR}}\iidx{k} &= \frac{g\iidx{k}\transposed y\iidx{k-1}}{\norm{g\iidx{k-1}}{}^2} \qquad &&\text{ (Polak and Ribi\`ere \cite{Polak1969Note} and Polyak \cite{Polyak1971conjugate})}, \\
	\beta^{\text{HS}}\iidx{k} &= \frac{g\iidx{k}\transposed y\iidx{k-1}}{\searchdirection\iidx{k-1}\transposed y\iidx{k-1}} \qquad &&\text{ (Hestenes and Stiefel \cite{Hestenes1952Methods})}, \\
	\beta^{\text{DY}}\iidx{k} &= \frac{\norm{g\iidx{k}}{}^2}{\searchdirection\iidx{k-1}\transposed y\iidx{k-1}} \qquad &&\text{ (Dai and Yuan \cite{Dai1999nonlinear})}, \\
	\beta^{\text{HZ}}\iidx{k} &= \left( y\iidx{k-1} - 2 \searchdirection\iidx{k-1} \frac{\norm{y\iidx{k-1}}{}^2}{d\iidx{k-1}\transposed y\iidx{k-1}} \right)\transposed \frac{g\iidx{k}}{\searchdirection\iidx{k-1}\transposed y\iidx{k-1}} \qquad &&\text{ (Hager and Zhang \cite{Hager2005new})}.
\end{alignat*}
Compared to the gradient descent method, NCG methods only need to store one or two additional vectors, depending on which method is used, while usually being significantly more efficient than the former. In contrast, L-BFGS methods need $2 m$ additional vectors of storage, where $m$ is the size of the memory, but usually perform slightly better than NCG methods. However, for very large-scale problems arising, e.g., from industrial applications, only the gradient descent method, L-BFGS method with $m=1$, and NCG methods may be feasible (cf. \cite{Kelley1999Iterative}). Finally, we remark that a detailed description of NCG methods for finite-dimensional optimization problems can be found in, e.g., \cite{Hager2006survey, Kelley1999Iterative, Nocedal2006Numerical}.

\subsection{Shape Calculus}
\label{ssec:shape_calculus}

A general PDE constrained shape optimization problem can be written in the form
\begin{equation*}
	\min_{\Omega \in \admissiblegeom} \costfunction(\Omega, \state) \qquad \text{ s.t. } \qquad e(\Omega, \state, \test) = 0 \quad \text{ for all } \test,
\end{equation*}
where $\costfunction$ is a cost functional which is to be optimized over a set of admissible geometries $\admissiblegeom$. Moreover, $e(\Omega, \state, \test)$ is a state equation, i.e., a PDE constraint, given on the domain $\Omega\subset \R^d$ with state $\state$ and test function $\test$ which is interpreted in the weak form
\begin{equation*}
	\text{Find } \state \text{ such that } \qquad e(\Omega, \state, \test) = 0 \qquad \text{ for all } \test.
\end{equation*}
As usual, we assume that the PDE constraint admits a unique solution $\state = \state(\Omega)$ so that $e(\Omega, \state(\Omega), \test) = 0$ for all $\test$. Hence, we can introduce the so-called reduced cost functional $\reducedcostfunction(\Omega) = \costfunction(\Omega, \state(\Omega))$ and consider the equivalent reduced problem
\begin{equation}
	\label{eq:abstract_sop}
	\min_{\Omega \in \admissiblegeom}\ \reducedcostfunction(\Omega).
\end{equation}

As a prototypical example for such a problem, we consider the following one from \cite{Etling2020First},
\begin{equation}
	\label{eq:poisson_sop}
	\min_{\Omega \in \admissiblegeom}\ \costfunction(\Omega, \state) = \integral{\Omega} \state \dmeas{x} \qquad \text{ s.t. } \qquad -\laplace \state = f \text{ in } \Omega, \quad \state = 0 \text{ on } \Gamma,
\end{equation}
where $\Gamma = \partial\Omega$ denotes the boundary of $\Omega$ and the set of admissible geometries is given by
\begin{equation*}
	\admissiblegeom = \Set{\Omega \subset \holdall | \Omega \text{ is a bounded Lipschitz domain}},
\end{equation*}
for some bounded hold-all domain $\holdall \subset \R^d$. The state equation for this problem is given by a Poisson problem, whose weak form reads
\begin{equation}
	\label{eq:weak_poisson}
	\text{Find } \state \in H^1_0(\Omega) \text{ such that } \qquad \integral{\Omega} \grad \state \cdot \grad \test \dmeas{x} = \integral{\Omega} f \test \dmeas{x} \qquad \text{ for all } \test\in H^1_0(\Omega),
\end{equation}
where $f\in H^1(\holdall)$, which is needed to ensure that the shape derivative exists (cf. \cite{Delfour2011Shapes}). Due to the Lax-Milgram Lemma (see, e.g., \cite{Evans2010Partial}), we know that \eqref{eq:weak_poisson} has a unique weak solution $\state=\state(\Omega)$ for any bounded Lipschitz domain $\Omega$, in particular, for any $\Omega \in \admissiblegeom$. Hence, problem \eqref{eq:poisson_sop} can be recast into the form \eqref{eq:abstract_sop} with $\reducedcostfunction(\Omega) = \costfunction(\Omega, \state(\Omega))$. 

Shape calculus is used to compute sensitivities of shape functionals, i.e., functionals defined on a subset of the power set of $\R^d$, such as the one in \eqref{eq:poisson_sop}, w.r.t.\ infinitesimal deformations of the domain. For a detailed introduction to this topic we refer the reader, e.g., to the textbooks \cite{Sokolowski1992Introduction, Delfour2011Shapes}. We utilize the so-called speed method which transforms a domain $\Omega \subset \holdall$ by the flow of a vector field $\vectorfield$. In particular, we consider a vector field $\vectorfield \in C^k_0(\holdall; \R^d)$ for $k\geq1$, i.e., the space of all $k$-times continuously differentiable functions from $\holdall$ to $\R^d$ with compact support in $\holdall$. A point $x_0 \in \Omega$ is transported along the flow of $\vectorfield$ as described by the following initial value problem (IVP)
\begin{equation}
	\label{eq:ivp_trafo}
	\dot{x}(t) = \vectorfield(x(t)), \qquad x(0) = x_0.
\end{equation}
From classical ODE theory (cf. \cite{Deuflhard2002Scientific}) we know that the above IVP has a unique solution $x(t)$ for all $t \in [0, \tau]$ if $\tau > 0$ is sufficiently small. This allows us to define the flow of $\vectorfield$ as the mapping
\begin{equation*}
	\flow^\vectorfield \colon \R^d \to \R^d; \quad x_0 \mapsto \flow^\vectorfield x = x(t),
\end{equation*}
where $x(t)$ is the solution of \eqref{eq:ivp_trafo}. The flow $\flow^\vectorfield$ is a diffeomorphism (cf. \cite{Delfour2011Shapes}), which we now use to define the shape derivative.

\begin{Definition}
	\label{def:shape_derivative}
	Let $\mathcal{S} \subset \set{\Omega | \Omega \subset \holdall}$, $\reducedcostfunction\colon S \to \R$, and $\Omega \in \mathcal{S}$. Moreover, let $\vectorfield \in C^k_0(\holdall;\R^d)$ for $k\geq 1$, let $\flow^\vectorfield$ be the flow associated to $\vectorfield$, and let $\flow^\vectorfield(\Omega) \in \mathcal{S}$ for all $t \in [0, \tau]$ for a sufficiently small $\tau > 0$. We say that the shape functional $\reducedcostfunction$ has a Eulerian semi-derivative at $\Omega$ in direction $\vectorfield$ if the following limit exists
	\begin{equation*}
		d\reducedcostfunction(\Omega)[\vectorfield] := \lim_{t\searrow 0} \frac{\reducedcostfunction(\flow^\vectorfield(\Omega)) - \reducedcostfunction(\Omega)}{t} = \left. \frac{d}{dt} \reducedcostfunction(\flow^\vectorfield(\Omega)) \right\rvert_{t=0^+}.
	\end{equation*}
	Furthermore, let $\Xi$ be a topological vector subspace of $C^\infty_0(\holdall; \R^d)$. We say that $\reducedcostfunction$ is shape differentiable at $\Omega$ w.r.t.\ $\Xi$ if it has an Eulerian semi-derivative at $\Omega$ in all directions $\vectorfield \in \Xi$ and, additionally, the mapping
	\begin{equation*}
		\Xi \to \R; \quad \vectorfield \mapsto d\reducedcostfunction(\Omega)[\vectorfield]
	\end{equation*}
	is linear and continuous. In this case, we call $d\reducedcostfunction(\Omega)[\vectorfield]$ the shape derivative of $\reducedcostfunction$ at $\Omega$ w.r.t.\ $\Xi$ in direction $\vectorfield$. 
\end{Definition}

For our model shape optimization problem, we have the following shape derivative (cf. \cite{Etling2020First}).

\begin{Proposition}
	The reduced cost functional $\reducedcostfunction$ corresponding to problem \eqref{eq:poisson_sop} is shape differentiable and has the following shape derivative
	\begin{equation}
		\label{eq:vol_sd_poisson}
		d\reducedcostfunction(\Omega)[\vectorfield] = \integral{\Omega} \state\ \divergence{\vectorfield} \dmeas{x} + \integral{\Omega} \left( \left( \divergence{\vectorfield}I - (D\vectorfield + D\vectorfield\transposed) \right) \grad \state \right) \cdot \grad \adjoint - \divergence{f\vectorfield} \adjoint \dmeas{x},
	\end{equation}
	where $I$ is the identity matrix in $\R^d$, $\state$ solves the state equation \eqref{eq:weak_poisson}, and $\adjoint$ solves the adjoint equation
	\begin{equation}
		\label{eq:weak_adjoint_poisson}
		\text{Find } \adjoint \in H^1_0(\Omega) \text{ such that } \qquad \integral{\Omega} \grad \adjoint \cdot \grad \varphi \dmeas{x} = - \integral{\Omega} \varphi \dmeas{x} \qquad \text{ for all } \varphi \in H^1_0(\Omega).
	\end{equation}
\end{Proposition}

The above formula for the shape derivative of problem \eqref{eq:poisson_sop} is also known as volume formulation since it involves integrals over the domain $\Omega$. However, under certain smoothness assumptions on the domain $\Omega$ and the data for the PDE constraint, there exists an equivalent boundary formulation of the shape derivative due to the Structure Theorem, which we recall in the following. 

\begin{Theorem}[Structure Theorem]
	\label{thm:hadamard}
	Let $\reducedcostfunction$ be a shape functional which is shape differentiable at some $\Omega\subset \R^d$ and let $\Gamma = \partial \Omega$ be compact. Further, let $k\geq 0$ be an integer for which 
	\begin{equation*}
		d\reducedcostfunction(\Omega) \colon C^\infty_0(\holdall; \R^d) \to \R; \quad \vectorfield \mapsto d\reducedcostfunction(\Omega)[\vectorfield]
	\end{equation*}
	is continuous w.r.t.\ the $C^k_0(\holdall;\R^d)$ topology, and assume that $\Gamma$ is of class $C^{k+1}$. Then, there exists a continuous functional $\shapedistro\colon C^k(\Gamma) \to \R$ so that
	\begin{equation*}
		d\reducedcostfunction(\Omega)[\vectorfield] = \shapedistro[\vectorfield \cdot \normal],
	\end{equation*}
	where $n$ is the outer unit normal vector on $\Gamma$.
	In particular, if $\shapedistro \in L^1(\Gamma)$, we have
	\begin{equation*}
		d\reducedcostfunction(\Omega)[\vectorfield] = \integral{\Gamma} \shapedistro\ \vectorfield \cdot \normal \dmeas{s}.
	\end{equation*}
\end{Theorem}
Note, that the proof of the structure theorem can be found in \cite[Theorem 3.6 and Corollary 1]{Delfour2011Shapes}.

This form of the shape derivative given in the structure theorem is known as boundary formulation. For the model problem \eqref{eq:poisson_sop}, we have the following result (cf. \cite{Etling2020First}).
\begin{Proposition}
	Let $\Omega$ be a bounded domain with $C^{1,1}$ boundary $\Gamma$. Then, the shape derivative \eqref{eq:vol_sd_poisson} has the following equivalent boundary representation
	\begin{equation}
		\label{eq:surf_sd_poisson}
		d\reducedcostfunction(\Omega)[\vectorfield] = \integral{\Gamma} -\frac{\partial \state}{\partial \normal} \frac{\partial \adjoint}{\partial \normal}\ \vectorfield \cdot \normal \dmeas{s},
	\end{equation}
	where $\state$ solves \eqref{eq:weak_poisson} and $\adjoint$ solves \eqref{eq:weak_adjoint_poisson}. In the framework of Theorem~\ref{thm:hadamard}, we have $\shapedistro = - \frac{\partial \state}{\partial \normal} \frac{\partial \adjoint}{\partial \normal} \in L^2(\Gamma)$ as discussed in \cite{Etling2020First}.
\end{Proposition}

\begin{Remark}
	For our model problem \eqref{eq:poisson_sop}, there exist two equivalent formulations of the shape derivative, the volume formulation \eqref{eq:vol_sd_poisson} and the boundary formulation \eqref{eq:surf_sd_poisson}. Thanks to the structure theorem, we know that this is the case for any sufficiently smooth shape functional of the form
	\begin{equation*}
		\reducedcostfunction(\Omega) = \integral{\Omega} j_\Omega \dmeas{x},
	\end{equation*}
	with a sufficiently smooth function $j_\Omega \colon \Omega \to \R$, whose shape derivative can be written as
	\begin{equation*}
		d\reducedcostfunction_\Omega(\Omega)[\vectorfield] = \integral{\Omega} G [\vectorfield] \dmeas{x} \qquad \text{ and } \qquad d\reducedcostfunction_\Gamma(\Omega)[\vectorfield] = \integral{\Gamma} \shapedistro\ \vectorfield\cdot \normal \dmeas{s}.
	\end{equation*}
	Here, $G$ represents a linear (differential) operator acting on $\vectorfield$, and $\shapedistro$ is the functional from Theorem~\ref{thm:hadamard}. If $\Omega$ has a sufficiently regular boundary (cf. Theorem~\ref{thm:hadamard}) and the functional $\reducedcostfunction$ is sufficiently smooth, both formulations are equivalent and we have $d\reducedcostfunction(\Omega)[\vectorfield] = d\reducedcostfunction_\Omega(\Omega)[\vectorfield] = d\reducedcostfunction_\Gamma(\Omega)[\vectorfield]$. However, for the boundary formulation $d\reducedcostfunction_\Gamma$ we usually need more smoothness assumptions than for the volume formulation $d\reducedcostfunction_\Omega$. In particular, for a given setting there may only be a volume, but no boundary formulation of the shape derivative available. Finally, we remark that the volume formulation shows better approximation properties for the numerical solution of shape optimization problems with the finite element method (cf. \cite{Hiptmair2015Comparison}), which is why we consider only volume formulations for our numerical experiments in Section~\ref{sec:numerics}.
\end{Remark}

\subsection{Riemannian Shape Optimization and Steklov-Poinar\'e-Type Metrics}
\label{ssec:riemannian}

We start by recalling the Riemannian view on shape optimization introduced in \cite{Schulz2014Riemannian}, following the notations used, e.g., in \cite{Schulz2016Computational, Schulz2016Efficient, Schulz2014Riemannian, Geiersbach2020Stochastic}. We first consider connected and compact subsets $\Omega \subset \holdall \subset \R^2$ with $C^\infty$ boundary, where $\holdall$ is a bounded hold-all domain. As in \cite{Michor2006Riemannian}, we define the space of all smooth two-dimensional shapes as
\begin{equation*}
	\shapemanifold(S^1;\R^2) := \text{Emb}(S^1;\R^2) / \text{Diff}(S^1),
\end{equation*}
i.e., the set of all equivalence classes of $C^\infty$ embeddings of the unit circle $S^1 \subset \R^2$ into $\R^2$, given by $\text{Emb}(S^1;\R^2)$, where the equivalence relation is defined via the set of all $C^\infty$ diffeomorphisms of $S^1$ into itself, given by $\text{Diff}(S^1)$. Note, that this equivalence relation factors out reparametrizations as these do not change the underlying shape. In \cite{Kriegl1997convenient} it is shown that $\shapemanifold$ is in fact a smooth manifold. An element of $\shapemanifold(S^1;\R^2)$ is represented by a smooth curve $\Gamma\colon S^1 \to \R^2; \theta \mapsto \Gamma(\theta)$. Due to the equivalence relation, the tangent space at $\Gamma \in \shapemanifold$ is isomorphic to the set of all $C^\infty$ normal vector fields along $\Gamma$, i.e.,
\begin{equation*}
	\tangentspace_\Gamma \shapemanifold \cong \Set{h | h=\alpha \normal, \alpha\in C^\infty(\Gamma;\R)} \cong \Set{\alpha | \alpha \in C^\infty(\Gamma;\R)},
\end{equation*}
where $\normal$ is the unit outer normal vector on $\Gamma$, i.e., $\normal(\theta) \perp \Gamma'(\theta)$ for all $\theta \in S^1$, with $\Gamma'$ being the circumferential derivative as in \cite{Michor2006Riemannian}. Note, that for this manifold, several metrics are discussed in \cite{Michor2006Riemannian}, and that this viewpoint can also be extended to higher dimensions (cf. \cite{Michor2005Vanishing}). Furthermore, we remark that a similar idea is given by the Courant metric and the quotient space $\mathcal{F} / \mathcal{G}(\Omega_0)$, with the latter being similar to the shape space $B_e$ (see, e.g., \cite[Chapter~3]{Delfour2011Shapes}).


As in \cite{Schulz2016Efficient}, we consider the following Steklov-Poincar\'e-type metric $\scalarproduct_\Gamma$ at some $\Gamma \in \shapemanifold$ which is defined as
\begin{equation}
	\label{eq:def_riemannian}
	\scalarproduct_\Gamma \colon H^{\nicefrac{1}{2}}(\Gamma) \times H^{\nicefrac{1}{2}}(\Gamma) \to \R; \quad (\alpha, \beta) \mapsto \integral{\Gamma} \alpha (\stekpoinc_\Gamma)^{-1} \beta \dmeas{s}.
\end{equation}
Here $\stekpoinc_\Gamma$ is a symmetric and coercive operator defined by
\begin{equation*}
	\stekpoinc_\Gamma \colon H^{-\nicefrac{1}{2}}(\Gamma) \to H^{\nicefrac{1}{2}}(\Gamma); \quad \alpha \mapsto U \cdot \normal,
\end{equation*}
where $U \in H^1(\Omega; \R^d)$ solves the problem
\begin{equation}
	\label{eq:sp_bilinear}
	\text{Find } U \in H^1(\Omega;\R^d) \text{ such that } \qquad a_\Omega(U, V) = \integral{\Gamma} \alpha (V\cdot \normal) \dmeas{s} \qquad \text{ for all } V \in H^1(\Omega; \R^d),
\end{equation}
for a symmetric, continuous, and coercive bilinear form $a_\Omega\colon H^1(\Omega; \R^d) \times H^1(\Omega;\R^d) \to \R$. Note, that $\stekpoinc_\Gamma$ corresponds to a projected Steklov-Poincar\'e operator, as discussed in \cite{Schulz2016Efficient} and that \eqref{eq:sp_bilinear} is well-posed due to the Lax-Milgram Lemma (see, e.g., \cite{Evans2010Partial}). The norm induced by this metric is denoted by
\begin{equation*}
	\norm{\alpha}{\Gamma} = \sqrt{\scalarproduct_{\Gamma}(\alpha, \alpha)}, \qquad \text{ for } \alpha \in H^{\nicefrac{1}{2}}(\Gamma).
\end{equation*}
To define a Riemannian metric on $\shapemanifold$, we restrict the Steklov-Poincar\'e metric to the tangent space $\tangentspace_\Gamma \shapemanifold$, i.e., we consider $\scalarproduct_\Gamma$ as a mapping $\scalarproduct_\Gamma \colon \tangentspace_\Gamma \shapemanifold \times \tangentspace_\Gamma \shapemanifold \to \R$.

\begin{Remark}
	\label{rem:boundaries}
	The scalar product $\scalarproduct_\Gamma$ in \eqref{eq:def_riemannian} is given for the case that the entire boundary of $\Omega$ is variable, as it is the case, e.g., for our model problem \eqref{eq:poisson_sop}. If $\Gamma = \Gamma\subfix \cup \Gamma\subdef$, where $\Gamma\subfix$ is fixed and $\Gamma\subdef$ is variable, we can instead use the space
	\begin{equation*}
		H^1_{\Gamma\subfix}(\Omega; \R^d) = \Set{V \in H^1(\Omega;\R^d) | V = 0 \text{ on } \Gamma\subfix}
	\end{equation*}
	in \eqref{eq:sp_bilinear} and modify the bilinear form $a_\Omega$, the scalar product $\scalarproduct_\Gamma$, and the operator $\stekpoinc_\Gamma$ accordingly.
\end{Remark}

Let us now discuss the relation between the metric $\scalarproduct_\Gamma$ and shape calculus. To do so, we assume that the shape functional $\reducedcostfunction$ is shape differentiable and has a shape derivative of the form
\begin{equation*}
	d\reducedcostfunction(\Omega)[\vectorfield] = \integral{\Gamma} \shapedistro\ \vectorfield\cdot \normal \dmeas{s},
\end{equation*}
with $\shapedistro \in L^2(\Gamma)$. Then, the Riemannian shape gradient w.r.t.\ $\scalarproduct_\Gamma$ is given by $\rieshapegrad \in \tangentspace_\Gamma \shapemanifold$, which is the solution of the following Riesz problem
\begin{equation}
	\label{eq:riesz_riemannian}
	\text{Find } \rieshapegrad \in \tangentspace_\Gamma \shapemanifold \text{ such that } \qquad \scalarproduct_\Gamma(\rieshapegrad, \phi) = \integral{\Gamma} \shapedistro \phi \dmeas{s} \qquad \text{ for all } \phi \in \tangentspace_\Gamma \shapemanifold.
\end{equation}
Due to the definition of $\scalarproduct_\Gamma$, the solution of \eqref{eq:riesz_riemannian} is given by $\rieshapegrad = \stekpoinc_\Gamma \shapedistro$, in particular, we have that $\rieshapegrad = \gradientdefo \cdot \normal$, where $\gradientdefo$ solves
\begin{equation}
	\label{eq:def_gradient_defo}
	\text{Find } \gradientdefo \in H^1(\Omega; \R^d) \text{ such that } \qquad a_\Omega(\gradientdefo, \vectorfield) = d\reducedcostfunction(\Omega)[\vectorfield] \qquad \text{ for all } \vectorfield \in H^1(\Omega;\R^d).
\end{equation}
Due to the Lax-Milgram Lemma, this problem has a unique solution $\gradientdefo$ which we call the gradient deformation of $\reducedcostfunction$ at $\Omega$. The gradient deformation $\gradientdefo$ can be interpreted as an extension of the Riemannian shape gradient $\rieshapegrad$ to the entire domain $\Omega$. Note, that due to the coercivity of $a_\Omega$ there exists a constant $C>0$ so that
\begin{equation}
	\label{eq:gradient_defo_descent}
	d\reducedcostfunction(\Omega)[-\gradientdefo] = a_\Omega(-\gradientdefo, \gradientdefo) \leq -C \norm{\gradientdefo}{H^1(\Omega)}^2 \leq 0,
\end{equation}
i.e., a infinitesimal transformation of the domain $\Omega$ along the flow associated to the negative gradient deformation yields a descent in the shape functional $\reducedcostfunction$. This fact is often used for the numerical solution of shape optimization problems, e.g., in \cite{Etling2020First, Schulz2016Computational, Blauth2019Model, Hohmann2019Shape}, where the domain is successively deformed according to the negative gradient deformation in the sense of a gradient descent method. We discuss this more detailedly in Section~\ref{sec:algorithms}. 

Let us end this section with the following remarks.
\begin{Remark}
	For the numerical solution of shape optimization problems, where the domain $\Omega$ is discretized by a mesh, it is desirable to obtain gradient deformations that lead to smooth mesh deformations. For this reason, the bilinear form $a_\Omega$ is often chosen according to the linear elasticity equations (see, e.g., \cite{Schulz2016Efficient, Schulz2016Computational, Gangl2015Shape, Etling2020First, Blauth2019Model}), i.e., $a_\Omega$ is given by
	\begin{equation}
		\label{eq:linear_elasticity}
		\begin{aligned}
			&a_\Omega \colon H^1(\Omega;\R^d) \times H^1(\Omega;\R^d) \to \R; \\
			&(V, W) \mapsto a_\Omega(V, W) = \integral{\Omega} 2 \lamesecond\ \symgrad(V) : \symgrad(W) + \lamefirst\ \divergence{V} \divergence{W} + \damping\ V \cdot W \dmeas{x},
		\end{aligned}
	\end{equation}
	where $A:B$ denotes the Frobenius inner product between matrices $A$ and $B$, $\symgrad(V) = \nicefrac{1}{2}\left( DV + DV\transposed \right)$ is the symmetric part of the Jacobian $DV$, $\lamefirst$ and $\lamesecond$ are the so-called Lam\'e parameters, for which we assume $\lamesecond > 0$ and $2\lamesecond + d \lamefirst > 0$, and $\damping \geq 0$ is a damping parameter. The latter is required to be positive in case $a_\Omega$ is defined on $H^1(\Omega;\R^d)$ as this case corresponds to a pure Neumann problem, where $a_\Omega$ is not coercive in case $\damping = 0$. If we consider a setting as in Remark~\ref{rem:boundaries} and $a_\Omega$ is given on $H^1_{\Gamma\subfix}(\Omega;\R^d)$, the damping parameter is allowed to vanish as $a_\Omega$ is also coercive for $\damping = 0$ in this case due to the Poincar\'e inequality. We remark that for the numerical experiments in this paper (cf. Section~\ref{sec:numerics}) we only consider bilinear forms $a_\Omega$ of the form \eqref{eq:linear_elasticity}. Finally, note that each different choice for the bilinear form $a_\Omega$ leads to a different Steklov-Poincar\'e operator $\stekpoinc_\Gamma$ and, hence, to a different Riemannian metric $\scalarproduct_{\Gamma}$. 
\end{Remark}

\begin{Remark}
	As discussed in \cite{Schulz2016Efficient}, the Riemannian shape gradient $\rieshapegrad = \gradientdefo \cdot \normal$ is not necessarily an element of the tangent space $\tangentspace_\Gamma \shapemanifold$ as it is not guaranteed that $\gradientdefo$ is in $C^\infty$. However, if the bilinear form $a_\Omega$ arises from a second order PDE with sufficiently smooth coefficients and the distribution of the shape derivative $\shapedistro$ as well as the domain $\Omega$ are sufficiently smooth, then the gradient deformation $\gradientdefo \in H^1(\Omega; \R^d)$ is indeed an element of $C^\infty(\overline{\Omega}; \R^d)$ by the theorem of infinite differentiability (cf. \cite[Theorem~6, Section~6.3]{Evans2010Partial}), in which case $\rieshapegrad \in \tangentspace_\Gamma \shapemanifold$.
\end{Remark}

\section{Algorithmic Solution of Shape Optimization Problems}
\label{sec:algorithms}

Having recalled shape calculus as well the Riemannian view on it from \cite{Schulz2014Riemannian} and the Steklov-Poincar\'e-type metrics from \cite{Schulz2016Efficient}, we now formulate a general algorithmic framework for gradient-based shape optimization methods. Throughout this section we consider the general form of a PDE constrained shape optimization problem given by the reduced problem
\begin{equation*}
	\min_{\Omega \in \admissiblegeom} \reducedcostfunction(\Omega).
\end{equation*}
We follow the algorithmic ideas from \cite{Ring2012Optimization} and formulate the general descent algorithm on the shape space $\shapemanifold$ based on the Steklov-Poincar\'e metric $\scalarproduct_\Gamma$. Afterwards, we propose NCG methods for shape optimization in this setting. Finally, we give an alternative formulation of the descent algorithm based on domain deformations which is well-suited for numerical discretization.

\subsection{General Descent Algorithm for Shape Optimization}
\label{ssec:general_descent_algo}

Before we describe our algorithmic framework, we briefly recall the concepts of a retraction and vector transport. For a detailed discussion of these topics we refer the reader, e.g., to \cite{Absil2008Optimization, Ring2012Optimization}. For some $\Gamma \in \shapemanifold$, a retraction $\retraction_\Gamma$ is a smooth mapping $\retraction_\Gamma \colon \tangentspace_\Gamma \shapemanifold \to \shapemanifold$ which satisfies $\retraction_{\Gamma} (0_\Gamma) = \Gamma$ and $D\retraction_\Gamma (0_\Gamma) = \text{id}_{\tangentspace_\Gamma \shapemanifold}$, where $D\retraction_\Gamma$ denotes the derivative of $\retraction_\Gamma$, and $0_\Gamma$ and $\text{id}_{\tangentspace_\Gamma \shapemanifold}$ are the zero element and identity mapping in $\tangentspace_\Gamma \shapemanifold$, respectively. A special case of a retraction is given by the so-called exponential map $\exp_\Gamma \colon \tangentspace_\Gamma \shapemanifold \to \shapemanifold; \alpha \mapsto \exp_\Gamma \alpha = \geodesic_\Gamma^\alpha(1)$, where $\geodesic_\Gamma^\alpha \colon [0,1] \to \shapemanifold$ is the unique geodesic starting at $\geodesic_\Gamma^\alpha(0) = \Gamma$ with $\dot{\geodesic}_\Gamma^\alpha(0) = \alpha$. Moreover, as in \cite{Schulz2015Structured} we denote by $\vectortransport \colon \tangentspace \shapemanifold \oplus \tangentspace \shapemanifold \to \tangentspace \shapemanifold; (\alpha, \beta) \mapsto \vectortransport_{\alpha} \beta$ a vector transport which satisfies the following properties. For $\alpha, \beta \in \tangentspace_\Gamma \shapemanifold$ it holds that $\vectortransport_\alpha \beta$ is an element of $\tangentspace_{\retraction_\Gamma (\alpha)} \shapemanifold$. Further, we have the relations $\vectortransport_{0_\Gamma} \alpha = \alpha$ and $\vectortransport_{\gamma} (a \alpha + b \beta) = a \vectortransport_\gamma \alpha + b \vectortransport_{\gamma} \beta$ for $a, b \in \R$. Note, that the parallel transport or parallel translation, as defined in, e.g., \cite{Absil2008Optimization, Ring2012Optimization}, is a vector transport associated to the exponential map.

\begin{algorithm2e}[b]
	\KwIn{Initial geometry, represented by $\Gamma\iidx{0}$, initial step size $t\iidx{0}$, tolerance $\texttt{tol} \in (0,1)$, maximum number of iterations $k_{\text{max}}$, parameters for the Armijo rule $\armijoparam \in (0, 1)$ and $\stepparam \in (0,1)$ \label{line:input}}
	\For{k=0,1,2,\dots, $k_\textrm{max}$}{
		Compute the solution of the state and adjoint systems \label{line:pde_solves} \\
		Compute the Riemannian shape gradient $\rieshapegrad\iidx{k} \in \tangentspace_{\Gamma\iidx{k}} \shapemanifold$ from \eqref{eq:riesz_riemannian} \label{line:shape_gradient} \\
		\If{$\norm{\rieshapegrad\iidx{k}}{\Gamma\iidx{k}} \leq \textup{\texttt{tol}} \norm{\rieshapegrad\iidx{0}}{\Gamma\iidx{0}}$ \label{line:stopping}}{
			Stop with approximate solution $\Gamma\iidx{k}$
		}
		Compute a search direction $\searchdirection\iidx{k} \in \tangentspace_{\Gamma\iidx{k}} \shapemanifold$ \label{line:search_direction} \\
		\If{$\scalarproduct_{\Gamma\iidx{k}}(\searchdirection\iidx{k}, \rieshapegrad\iidx{k}) > 0$ \label{line:descent_test}}{
			Set $\searchdirection\iidx{k} = - \rieshapegrad\iidx{k}$ \label{line:descent_failed} \\
		}
		\While{$\reducedcostfunction(\retraction_{\Gamma\iidx{k}}(t \searchdirection\iidx{k})) > \reducedcostfunction(\Gamma\iidx{k}) + \armijoparam t\ \scalarproduct_{\Gamma\iidx{k}}(\rieshapegrad\iidx{k}, \searchdirection\iidx{k})$ \label{line:armijo}}{
			Decrease the step size: $t = \stepparam t$ \label{line:armijo_fail} \\
		}
		Set $t\iidx{k} = t, \increment\iidx{k} = t\iidx{k} \searchdirection\iidx{k}$, and update the geometry via $\Gamma\iidx{k+1} = \retraction_{\Gamma\iidx{k}}(\increment\iidx{k})$ \label{line:armijo_success} \\
		Increase the step size for the next iteration: $t = \nicefrac{t\iidx{k}}{\stepparam}$ \label{line:update_stepsize} \\
	}
	\caption{General descent algorithm for shape optimization problems.}
	\label{algo:descent}
\end{algorithm2e}

Let us now investigate the general descent algorithm for shape optimization, given in Algorithm~\ref{algo:descent}. After specifying the initial domain $\Omega\iidx{0}$ via its boundary $\Gamma\iidx{0}$ as well as several parameters for the input, the algorithm proceeds as follows: in iteration $k$ we consider the iterate $\Gamma\iidx{k}$ with corresponding domain $\Omega\iidx{k}$. First, the state and adjoint systems corresponding to the shape optimization problem are solved in line~\ref{line:pde_solves} so that we can calculate the shape derivative $d\reducedcostfunction(\Omega\iidx{k})[\cdot]$. In line~\ref{line:shape_gradient} we solve the Riesz problem \eqref{eq:riesz_riemannian} to obtain the Riemannian shape gradient $\rieshapegrad\iidx{k}$, which by definition of $\scalarproduct_\Gamma$ involves the computation of the gradient deformation as an intermediate step. Afterwards, in line~\ref{line:stopping}, the convergence of the algorithm is tested by a relative stopping criterion involving the norm of the Riemannian shape gradient $\norm{\rieshapegrad\iidx{k}}{\Gamma\iidx{k}}$. In line~\ref{line:search_direction} the search direction $\searchdirection\iidx{k}$ is computed, based on the information from $\rieshapegrad\iidx{k}$. In this paper, we consider different choices for the search direction which lead to either the gradient descent, L-BFGS, or NCG methods, as detailed below. Since, in general, the search direction cannot be guaranteed to be a descent direction we need to make sure that we exclude the case where it leads to an ascent in $\reducedcostfunction$, which is done in line~\ref{line:descent_test}. As is explained in Section~\ref{ssec:algo_numerical}, if we have $\searchdirection\iidx{k} = \searchdefo\iidx{k} \cdot \normal$ for some vector field $\searchdefo\iidx{k}$, then it holds that $\scalarproduct_{\Gamma\iidx{k}}(\searchdirection\iidx{k}, \rieshapegrad\iidx{k}) = d\reducedcostfunction(\Omega\iidx{k})[\searchdefo\iidx{k}]$ (cf. \eqref{eq:relation_scalar_products}). Hence, an infinitesimal deformation of $\Omega\iidx{k}$ along the flow of $\searchdefo\iidx{k}$ leads to an ascent in the cost functional $\reducedcostfunction$ if $\scalarproduct_{\Gamma\iidx{k}}(\searchdirection\iidx{k}, \rieshapegrad\iidx{k}) > 0$. In this case we reinitialize the search direction to the negative Riemannian shape gradient in line~\ref{line:descent_failed}, which is guaranteed to be a descent direction thanks to $d\reducedcostfunction(\Omega\iidx{k})[-\gradientdefo\iidx{k}] = \scalarproduct_{\Gamma\iidx{k}}(\rieshapegrad\iidx{k}, -\rieshapegrad\iidx{k}) = - \norm{\rieshapegrad\iidx{k}}{\Gamma\iidx{k}}^2 \leq 0$. After having computed the search direction, we employ an Armijo line search in lines~\ref{line:armijo} and~\ref{line:armijo_fail} (cf. \cite{Ring2012Optimization}). For the sake of better readability, in line~\ref{line:armijo} of the algorithm we write $\reducedcostfunction(\Gamma)$ instead of $\reducedcostfunction(\Omega)$ to be compatible with the Riemannian framework. If $\searchdirection\iidx{k} = \searchdefo\iidx{k} \cdot \normal$ for some vector field $\searchdefo\iidx{k}$, then using the same arguments as above reveals that the inequality from line~\ref{line:armijo} reduces to the classical form
\begin{equation*}
	\reducedcostfunction(\retraction_{\Gamma\iidx{k}}(t \searchdirection\iidx{k})) > \reducedcostfunction(\Gamma\iidx{k}) + \armijoparam t\ d\reducedcostfunction(\Gamma\iidx{k})[\searchdefo\iidx{k}].
\end{equation*}
Note, that we terminate the algorithm if the trial step size $t$ becomes too small in order to prevent being stuck in the line search. If the Armijo line search is terminated successfully, we update the geometry in line~\ref{line:armijo_success} via $\Gamma\iidx{k+1} = \retraction_{\Gamma\iidx{k}}(\increment\iidx{k})$, where $\increment\iidx{k} = t\iidx{k} \searchdirection\iidx{k}$ is the corresponding increment and $\retraction$ is a retraction as defined previously. Finally, we increase the step size in line~\ref{line:update_stepsize} so that the algorithm uses a larger initial trial step size in the next iteration. 

For the choice of the search direction we have the following remarks. For the gradient descent method, we use $\searchdirection\iidx{k} = - \rieshapegrad\iidx{k}$, which is guaranteed to be a descent direction as discussed above. The L-BFGS methods are detailed, e.g., in \cite{Schulz2016Efficient, Schulz2016Computational, Schulz2015Structured}, where a double loop for the computation of the search direction is given. We have two remarks concerning our implementation details of the L-BFGS methods. First, we note that the curvature condition for the L-BFGS methods reads 
\begin{equation}
	\label{eq:curvature_condition}
	\scalarproduct_{\Gamma\iidx{k}}(s\iidx{k-1}, y\iidx{k-1}) > 0,
\end{equation}
where we use the notation
\begin{equation*}
	s\iidx{k-1} = \vectortransport_{\increment\iidx{k-1}} \increment\iidx{k-1} \in \tangentspace_{\Gamma\iidx{k}} \shapemanifold \qquad \text{ and } \qquad y\iidx{k-1} = \rieshapegrad\iidx{k} - \vectortransport_{\increment\iidx{k-1}} \rieshapegrad\iidx{k-1} \in \tangentspace_{\Gamma\iidx{k}} \shapemanifold,
\end{equation*}
for a vector transport $\vectortransport$ with associated retraction $\retraction$. If condition \eqref{eq:curvature_condition} is not satisfied for some $k$, we restart the L-BFGS methods with a gradient step, as discussed in \cite{Kelley1999Iterative}, i.e., we set $s\iidx{i} = 0$ and $y\iidx{i} = 0$ for all $i \leq k$. Second, as discussed in, e.g., \cite{Nocedal2006Numerical}, the L-BFGS methods have, similarly to a Newton method, a built-in scaling of the search direction so that we always consider a step size of $t=1$ as initial guess for the Armijo line search in case the methods have a non-empty memory.

\subsection{NCG Methods for Shape Optimization}
\label{ssec:ncg_shape}

We now formulate nonlinear conjugate gradient methods for the solution of shape optimization problems, which we embed into the algorithmic framework discussed previously. To do so, we only have to specify the computation of the search direction in line~\ref{line:search_direction} of Algorithm~\ref{algo:descent}, the rest of the algorithm remains unchanged. Analogously to \cite{Ring2012Optimization, Absil2008Optimization}, the search direction for the NCG methods is computed by the following formula
\begin{equation*}
	\searchdirection\iidx{k} = - \rieshapegrad\iidx{k} + \beta\iidx{k} \vectortransport_{\increment\iidx{k-1}} \searchdirection\iidx{k-1} \in \tangentspace_{\Gamma\iidx{k}} \shapemanifold,
\end{equation*}
where, as before, $\increment\iidx{k} = t\iidx{k}\searchdirection\iidx{k}$ is the increment of the iterate from line~\ref{line:armijo_success}. Analogously to the finite-dimensional case, we initialize the search direction with the negative shape gradient, i.e., $\beta\iidx{0} = 0$ or, equivalently, $\searchdirection\iidx{0} = -\rieshapegrad\iidx{0}$. For the computation of the the update parameter $\beta\iidx{k}$ we define
\begin{equation*}
	y\iidx{k-1} = \rieshapegrad\iidx{k} - \vectortransport_{\increment\iidx{k-1}} \rieshapegrad\iidx{k-1} \in \tangentspace_{\Gamma\iidx{k}} \shapemanifold.
\end{equation*}
The parameter $\beta\iidx{k}$ is given in analogy to the finite-dimensional ones from Section~\ref{ssec:ncg_finite} and \cite{Ring2012Optimization, Absil2008Optimization} as follows
\begin{equation}
	\label{eq:id_ncg}
	\begin{aligned}
		\beta^\text{FR}\iidx{k} &= \frac{\norm{\rieshapegrad\iidx{k}}{\Gamma\iidx{k}}^2}{\norm{\vectortransport_{\increment\iidx{k-1}} \rieshapegrad\iidx{k-1}}{\Gamma\iidx{k}}^2}, \\
		\beta^\text{PR}\iidx{k} &= \frac{\scalarproduct_{\Gamma\iidx{k}}\left( \rieshapegrad\iidx{k}, y\iidx{k-1} \right)}{\norm{\vectortransport_{\increment\iidx{k-1}} \rieshapegrad\iidx{k-1}}{\Gamma\iidx{k}}^2}, \\
		\beta^\text{HS}\iidx{k} &= \frac{\scalarproduct_{\Gamma\iidx{k}}\left( \rieshapegrad\iidx{k}, y\iidx{k-1} \right)}{\scalarproduct_{\Gamma\iidx{k}}\left(\vectortransport_{\increment\iidx{k-1}} \searchdirection\iidx{k-1}, y\iidx{k-1} \right)}, \\
		\beta^\text{DY}\iidx{k} &= \frac{\norm{\rieshapegrad\iidx{k}}{\Gamma\iidx{k}}^2}{\scalarproduct_{\Gamma\iidx{k}}\left( \vectortransport_{\increment\iidx{k-1}} \searchdirection\iidx{k-1}, y\iidx{k-1} \right)}, \\
		\beta^\text{HZ}\iidx{k} &= \scalarproduct_{\Gamma\iidx{k}} \left( y\iidx{k-1} - 2\vectortransport_{\increment\iidx{k-1}} \searchdirection\iidx{k-1} \frac{\norm{y\iidx{k-1}}{\Gamma\iidx{k}}^2}{\scalarproduct_{\Gamma\iidx{k}}\left( \vectortransport_{\increment\iidx{k-1}} \searchdirection\iidx{k-1}, y\iidx{k-1} \right)} , \frac{\rieshapegrad\iidx{k}}{\scalarproduct_{\Gamma\iidx{k}}\left( \vectortransport_{\increment\iidx{k-1}} \searchdirection\iidx{k-1}, y\iidx{k-1} \right)} \right).
	\end{aligned}
\end{equation}

As described in \cite{Nocedal2006Numerical}, it can be beneficial to restart the NCG methods with a gradient step, for which there are two popular methods. First, one can reinitialize the search direction to the negative gradient direction every $\cgiter$ iterations, which can possibly enhance the convergence (cf. \cite{Nocedal2006Numerical}). Second, NCG methods try to generate directions for which the corresponding gradients are orthogonal, which indeed holds for the classical linear CG method, i.e., if the cost functional is a strictly convex quadratic function, but is only satisfied approximately for general nonlinear cost functionals. Hence, we reinitialize the search direction with the negative shape gradient if the following criterion (cf. \cite{Nocedal2006Numerical}) is satisfied
\begin{equation}
	\label{eq:orthogonality_test}
	\frac{\scalarproduct_{\Gamma\iidx{k}}\left( \rieshapegrad\iidx{k}, \vectortransport_{\increment\iidx{k-1}} \rieshapegrad\iidx{k-1} \right)}{\norm{\rieshapegrad\iidx{k}}{\Gamma\iidx{k}}^2} \geq \cgtol,
\end{equation}
where $\cgtol$ is a parameter usually chosen in $(0,1)$. We use the notation $\cgiter = \infty$ and $\cgtol = \infty$ in case we do not restart the methods after a fixed amount of iterations or via the condition \eqref{eq:orthogonality_test}, respectively. Note, that the convergence properties of the Fletcher-Reeves and Polak-Ribi\`ere methods on infinite-dimensional manifolds are investigated in \cite{Ring2012Optimization}.

\begin{Remark}
	For the Fletcher-Reeves and Polak-Ribi\`ere methods, one usually considers the term $\norm{\rieshapegrad\iidx{k-1}}{\Gamma\iidx{k-1}}^2$ as the denominator for $\beta\iidx{k}$. However, as all other quantities for the algorithm are computed at the current iterate $\Gamma\iidx{k}$, we instead use the term $\norm{\vectortransport_{\increment\iidx{k-1}} \rieshapegrad\iidx{k-1}}{\Gamma\iidx{k}}^2$ as this is also evaluated at $\Gamma\iidx{k}$. Note, that both expressions are equivalent if the vector transport $\vectortransport$ is an isometry, which is, e.g., satisfied by the parallel transport (cf. \cite{Absil2008Optimization}). Moreover, we note that for the numerical experiments in Section~\ref{sec:numerics} the formulation given in \eqref{eq:id_ncg} yields slightly better results.
\end{Remark}

\subsection{Volume-Based Description of the Descent Algorithm}
\label{ssec:algo_numerical}

The description of the descent algorithm given in Algorithm~\ref{algo:descent} is not yet well-suited for numerical discretization as it focused on an abstract infinite-dimensional setting involving the shape space $\shapemanifold$. In particular, the deformation of the geometry via retractions, which represents only a deformation of the boundary, is not efficient for the numerical solution of shape optimization problems due to the following. From the numerical point of view, the domain $\Omega$ is discretized, e.g., by some kind of mesh. Only deforming the boundary of this mesh is not suitable as, in general, this causes the mesh to degenerate even for comparatively small deformations, which necessitates costly remeshing. Hence, it is more efficient to deform the entire mesh directly, and not only its boundary. As mentioned in Section~\ref{ssec:riemannian}, this is usually done by computing the gradient deformation $\gradientdefo$ via \eqref{eq:def_gradient_defo}, where the bilinear form $a_\Omega$ is chosen according to the linear elasticity equations \eqref{eq:linear_elasticity}, and deforming the the geometry accordingly, which results in smooth mesh deformations. For these reasons, we now make use of the fact that the gradient deformation can be viewed as extension of the Riemannian shape gradient to the entire domain (cf. Section~\ref{ssec:riemannian}) and formulate the descent algorithm in this volume-based setting, which is well-suited for numerical discretization. 

Due to the definition of the Steklov-Poincar\'e metric $\scalarproduct_\Gamma$, we have to compute the gradient deformation $\gradientdefo$ as an intermediate result for the computation of the Riemannian shape gradient $\rieshapegrad = \gradientdefo \cdot \normal$. Moreover, we will see that for the search directions $\searchdirection$ used in Algorithm~\ref{algo:descent} we also have $\searchdirection = \searchdefo \cdot \normal$ for some vector field $\searchdefo$. We then have the following relations between the scalar product $\scalarproduct_\Gamma$ and the bilinear form $a_\Omega$. 
Let the Riemannian shape gradient be given by $\rieshapegrad \in \tangentspace_\Gamma \shapemanifold$, where $\rieshapegrad = \gradientdefo \cdot \normal$ as in Section~\ref{ssec:riemannian}, and consider an element $\beta \in \tangentspace_\Gamma \shapemanifold$ given by $\beta = \vectorfield\cdot \normal$ for some vector field $\vectorfield$. Then, we have that
\begin{equation}
	\label{eq:relation_scalar_products}
	\scalarproduct_\Gamma\left( \rieshapegrad, \beta \right) = \integral{\Gamma} \shapedistro \beta \dmeas{s} = \integral{\Gamma} \shapedistro\ \vectorfield \cdot \normal \dmeas{s} = d\reducedcostfunction(\Omega)[\vectorfield] = a_\Omega(\gradientdefo, \vectorfield).
\end{equation}
In particular, we can compute the norm of the Riemannian shape gradient via
\begin{equation}
	\label{eq:relation_norm}
	\norm{\rieshapegrad}{\Gamma}^2 = a_\Omega(\gradientdefo, \gradientdefo) = \norm{\gradientdefo}{a_\Omega}^2,
\end{equation}
where $\norm{\cdot}{a_\Omega}$ is a norm on $H^1(\Omega;\R^d)$ defined as
\begin{equation*}
	\norm{\vectorfield}{a_\Omega} = \sqrt{a_\Omega\left( \vectorfield, \vectorfield \right)}, \qquad \vectorfield \in H^1(\Omega;\R^d).
\end{equation*}
From these considerations we see that the metric $\scalarproduct_\Gamma$ embeds the gradient deformation $\gradientdefo$ into a Riemannian framework for shape optimization. 

The only thing left to do before we can give the volume-based description of Algorithm~\ref{algo:descent} is to extend the retraction and vector transport so that they act on vector fields instead of elements of the tangent space $\tangentspace_\Gamma \shapemanifold$. We denote these extensions by
\begin{equation*}
	\numretraction_\Omega \colon H^1(\Omega;\R^d) \to \admissiblegeom \qquad \text{ and } \qquad \numtransport \colon H^1(\Omega;\R^d) \times H^1(\Omega;\R^d) \to H^1(\numretraction_\Omega(\vectorfield); \R^d),
\end{equation*}
where the retraction $\numretraction_\Omega$ maps a vector field $\vectorfield \in H^1(\Omega;\R^d)$ to a domain $\numretraction_\Omega(\vectorfield) \subset \admissiblegeom$ deformed by it, and $\numtransport$ maps a vector field $\mathcal{W} \in H^1(\Omega,\R^d)$ to a vector field $\numtransport_{\vectorfield} \mathcal{W} \in H^1(\numretraction_\Omega(\vectorfield))$. Moreover, to be compatible with the previous notions of a retraction and vector transport, we assume that for $\alpha = \vectorfield \cdot \normal \in \tangentspace_\Gamma \shapemanifold$ and $\beta = \mathcal{W} \cdot \normal \in \tangentspace_\Gamma \shapemanifold$ with $\vectorfield, \mathcal{W} \in H^1(\Omega;\R^d)$ it holds that $\retraction_\Gamma (\alpha)$ is the boundary of the domain $\numretraction_\Omega (\vectorfield)$ and that 
\begin{equation*}
	\vectortransport_\alpha \beta = \left(\numtransport_\vectorfield \mathcal{W}\right) \cdot \normal \in \tangentspace_{\retraction_\Gamma (\alpha)} \shapemanifold.
\end{equation*}
This extends the previous notions of retraction and vector transport to our volume-based setting.

\begin{algorithm2e}[b]
	\KwIn{Initial geometry $\Omega\iidx{0}$, initial step size $t\iidx{0}$, tolerance $\texttt{tol} \in (0,1)$, maximum number of iterations $k_{\text{max}}$, parameters for the Armijo rule $\armijoparam \in (0, 1)$ and $\stepparam \in (0,1)$ }
	\For{k=0,1,2,\dots, $k_\textrm{max}$}{
		Compute the solution of the state and adjoint systems \label{lin:pde_solves} \\
		Compute the gradient deformation $\gradientdefo\iidx{k}$ by solving \eqref{eq:def_gradient_defo} \label{lin:gradient_defo} \\
		\If{$\norm{\gradientdefo\iidx{k}}{a_{\Omega\iidx{k}}} \leq \textup{\texttt{tol}} \norm{\gradientdefo\iidx{0}}{a_{\Omega\iidx{0}}}$ \label{lin:stop_test}}{
			Stop with approximate solution $\Omega\iidx{k}$
		}
		Compute a search direction $\searchdefo\iidx{k}$ \label{lin:search_direction} \\
		\If{$a_{\Omega\iidx{k}}\left( \searchdefo\iidx{k} , \gradientdefo\iidx{k} \right) > 0$ \label{lin:descent_test}}{
			Set $\searchdefo\iidx{k} = - \gradientdefo\iidx{k}$ \label{lin:descent_fail} \\
		}
		\While{$\reducedcostfunction(\numretraction_{\Omega\iidx{k}} (t\searchdefo\iidx{k})) > \reducedcostfunction(\Omega\iidx{k}) + \armijoparam t\ a_{\Omega\iidx{k}}\left( \gradientdefo\iidx{k}, \searchdefo\iidx{k} \right)$ \label{lin:armijo}}{
			Decrease the step size: $t = \stepparam t$ \label{lin:armijo_fail} \\
		}
		Set $t\iidx{k} = t, \xi\iidx{k} = t\iidx{k} \searchdefo\iidx{k}$, and update the geometry via $\Omega\iidx{k+1} = \numretraction_{\Omega\iidx{k}} (\xi\iidx{k})$ \label{lin:update_geometry} \\
		Increase the step size for the next iteration: $t = \nicefrac{t\iidx{k}}{\stepparam}$ \label{lin:increase_size} \\
	}
	\caption{Volume-based formulation of the descent algorithm for shape optimization problems.}
	\label{algo:transformed}
\end{algorithm2e}

Following the above discussion, it is now straightforward to modify Algorithm~\ref{algo:descent} to a formulation that uses vector fields as search directions and deforms the domain according to the retraction $\numretraction$, which is given in Algorithm~\ref{algo:transformed}. In particular, the scalar products and norms involving tangent vectors and the metric $\scalarproduct_\Gamma$ can now be rephrased as scalar products and norms involving vector fields and the bilinear form $a_\Omega$ (cf. \eqref{eq:relation_scalar_products} and \eqref{eq:relation_norm}). This formulation is stated in Algorithm~\ref{algo:transformed} and is straightforward to discretize consistently, e.g., by the finite element method, which we briefly discuss in Section~\ref{ssec:implementation}.

We have the following remarks regarding the computation of the search direction in the setting of Algorithm~\ref{algo:transformed}. For the gradient descent method, we use the search direction $\searchdefo\iidx{k} = -\gradientdefo\iidx{k}$, i.e., we use the negative gradient deformation to deform our domain. Due to \eqref{eq:gradient_defo_descent}, we know that this always yields a descent direction. For the computation of the search direction with the L-BFGS methods we refer the reader to \cite{Schulz2016Efficient}, where the corresponding double loop for this volume-based setting is given. In particular, the search direction for the L-BFGS methods is given as a linear combination of previous gradient deformations $\gradientdefo\iidx{i}$ and increments $\xi\iidx{i} = t\iidx{i}\searchdefo\iidx{k}$ (cf. line~\ref{lin:update_geometry} of Algorithm~\ref{algo:transformed}) of the geometry which are transported to the current domain $\Omega\iidx{k}$ via the vector transport $\numtransport$, making it a vector field as well. 

Let us now detail the computation of the search direction for the NCG methods in the setting of Algorithm~\ref{algo:transformed}. Analogously to before, we use the iteration
\begin{equation*}
	\searchdefo\iidx{k} = - \gradientdefo\iidx{k} + \beta\iidx{k} \numtransport_{\xi\iidx{k-1}} \searchdefo\iidx{k-1} \in H^1(\Omega\iidx{k};\R^d),
\end{equation*}
where $\xi\iidx{k} = t\iidx{k} \searchdefo\iidx{k}$ is the increment of the geometry as in line~\ref{lin:update_geometry} of Algorithm~\ref{algo:transformed}. Moreover, we have $\beta\iidx{0} = 0$ or, equivalently, $\searchdefo\iidx{0} = - \gradientdefo\iidx{0}$, as before. Obviously, as for the gradient descent and L-BFGS methods, we observe that this leads to search directions that are represented by vector fields on the domain $\Omega\iidx{k}$. In analogy to Section~\ref{ssec:ncg_shape}, we define $\graddiff\iidx{k-1} = \gradientdefo\iidx{k} - \numtransport_{\xi\iidx{k-1}} \gradientdefo\iidx{k-1} \in H^1(\Omega\iidx{k};\R^d)$ and state the modified update formulas for $\beta\iidx{k}$ in the following
\begin{align*}
		\beta^\text{FR}\iidx{k} &= \frac{\norm{\gradientdefo\iidx{k}}{a_{\Omega\iidx{k}}}^2}{\norm{\numtransport_{\xi\iidx{k-1}} \gradientdefo\iidx{k-1}}{a_{\Omega\iidx{k}}}^2}, \\
		\beta^\text{PR}\iidx{k} &= \frac{a_{\Omega}\left( \gradientdefo\iidx{k}, \graddiff\iidx{k-1} \right)}{\norm{\numtransport_{\xi\iidx{k-1}} \gradientdefo\iidx{k-1}}{a_{\Omega\iidx{k}}}^2}, \\
		\beta^\text{HS}\iidx{k} &= \frac{a_{\Omega}\left( \gradientdefo\iidx{k}, \graddiff\iidx{k-1} \right)}{a_{\Omega\iidx{k}}\left( \numtransport_{\xi\iidx{k-1}} \searchdefo\iidx{k-1}, \graddiff\iidx{k-1} \right)}, \\
		\beta^\text{DY}\iidx{k} &= \frac{\norm{\gradientdefo\iidx{k}}{a_{\Omega\iidx{k}}}^2}{a_{\Omega\iidx{k}}\left( \numtransport_{\xi\iidx{k-1}} \searchdefo\iidx{k-1}, \graddiff\iidx{k-1} \right)}, \\
		\beta^\text{HZ}\iidx{k} &= a_{\Omega\iidx{k}}\left( \graddiff\iidx{k-1} - 2\numtransport_{\xi\iidx{k-1}} \searchdefo\iidx{k-1} \frac{\norm{\graddiff\iidx{k-1}}{a_{\Omega\iidx{k}}}^2}{a_{\Omega\iidx{k}}\left( \numtransport_{\xi\iidx{k-1}} \searchdefo\iidx{k-1} , \graddiff\iidx{k-1} \right)} , \frac{\gradientdefo\iidx{k}}{a_{\Omega\iidx{k}}\left( \numtransport_{\xi\iidx{k-1}} \searchdefo\iidx{k-1}, \graddiff\iidx{k-1} \right)} \right).
\end{align*}
Regarding the restart of the algorithm, condition \eqref{eq:orthogonality_test} can be rewritten in this setting as
\begin{equation}
	\label{eq:restart_ncg}
	\frac{a_{\Omega\iidx{k}}\left( \gradientdefo\iidx{k} , \numtransport_{\xi\iidx{k-1}} \gradientdefo\iidx{k-1} \right)}{\norm{\gradientdefo\iidx{k}}{a_{\Omega\iidx{k}}}^2} \geq \cgtol,
\end{equation}
i.e., we reinitialize $\searchdefo\iidx{k} = - \gradientdefo\iidx{k}$ if \eqref{eq:restart_ncg} holds. 

The above discussions show that the Steklov-Poincar\'e metric $\scalarproduct_\Gamma$ from \cite{Schulz2016Efficient} enables us to define the NCG methods for shape optimization in terms of vector fields directly, which results in efficient methods for PDE constrained shape optimization as demonstrated in the following section.

\section{Numerical Comparison of Gradient-Based Algorithms for Shape Optimization}
\label{sec:numerics}

In this section, we numerically compare the NCG methods to the gradient descent and L-BFGS methods using four benchmark shape optimization problems. Throughout this section, we abbreviate the gradient descent method by GD, the L-BFGS method with memory $m$ by L-BFGS~$m$, and the NCG methods are abbreviated as NCG FR, PR, HS, DY, and HZ, corresponding to the Fletcher-Reeves, Polak-Ribi\`ere, Hestenes-Stiefel, Dai-Yuan, and Hager-Zhang variants, respectively. Finally, note that our implementation of the numerical experiments considered in this section is available freely on GitHub \cite{Blauth2021NCGCode}.

\subsection{Discretization and Setup}
\label{ssec:implementation}

We have implemented the algorithmic framework from Section~\ref{sec:algorithms} in our open-source software package cashocs \cite{Blauth2020CASHOCS}, version 1.2.1, which enables the automated solution of general shape optimization problems with the help of the finite element software FEniCS \cite{Alnes2015FEniCS, Logg2012Automated}. 
Note, that our software cashocs derives the corresponding adjoint system and shape derivatives of the respective problems automatically (cf.\ \cite{Blauth2020CASHOCS}). Since we use only conforming Galerkin finite element methods, the adjoint systems and shape derivatives computed by cashocs are consistent discretizations of the corresponding infinite-dimensional objects. Moreover, we remark that cashocs only uses the volume formulation of the shape derivative as this yields better results for the numerical solution of shape optimization problems with the finite element method (cf.\ \cite{Hiptmair2015Comparison}). For the state and adjoint systems, we detail the finite element discretization for each problem individually later on. For the discretization of the domain $\Omega$ we use the finite element meshes corresponding to the respective state and adjoint systems. Moreover, problem \eqref{eq:def_gradient_defo}, used to determine the gradient deformation, is discretized with linear Lagrange elements. Hence, the gradient deformation $\gradientdefo\iidx{k}$ and the search direction $\searchdefo\iidx{k}$ in Algorithm~\ref{algo:transformed} are discretized as piecewise linear functions. 

For the extended retraction $\numretraction$ we proceed analogously to \cite{Schulz2015Structured} and use
\begin{equation}
	\label{eq:num_retraction}
	\numretraction_\Omega (\vectorfield) = (I + \vectorfield)\Omega = \Set{x + \vectorfield(x) | x\in \Omega}.
\end{equation}
This is a similar concept to the perturbation of identity, which generates a family of transformed domains $\Omega_t$ for $t\geq 0$ through
\begin{equation*}
	\Omega_t = (I + t \vectorfield)\Omega = \Set{x + t \vectorfield(x) | x\in \Omega},
\end{equation*}
and presents an alternative to the speed method for computing first order shape derivatives (cf. \cite{Delfour2011Shapes}). Equation \eqref{eq:num_retraction} implies that we consider a Eulerian setting, where the underlying finite element mesh is deformed and moved in each iteration of the respective optimization algorithm. In particular, the state and adjoint systems as well as the linear elasticity equations for determining the gradient deformation are solved on the deformed domain. From the numerical point of view, we can easily realize this retraction by simply adding the (discretized) vector field $\vectorfield$ to the nodes of the finite element mesh. Note, that the deformation of the geometry only occurs during the Armijo line search in Algorithm~\ref{algo:transformed} and that we reject deformations that would result in inverted or intersecting mesh elements, so that we obtain a conforming finite element mesh for all iterations.

For the vector transport $\numtransport$ we use
\begin{equation}
	\label{eq:num_transport}
	\numtransport_{\vectorfield} \mathcal{W}(y) = \mathcal{W}(x) \qquad \text{ for } y = x + \vectorfield(x) \in \numretraction_\Omega (\vectorfield) \text{ with } x \in \Omega,
\end{equation}
where $\numtransport_{\vectorfield} \mathcal{W} \colon \numretraction_\Omega(\vectorfield) \to \R^d$ is a vector field on $\numretraction_\Omega(\vectorfield)$, and $\vectorfield \colon \Omega \to \R^d$ and $\mathcal{W} \colon \Omega \to \R^d$ are vector fields on $\Omega$. Due to the definition of the retraction and the fact that we exclude deformations that would lead to inverted or intersecting mesh elements as discussed above, the retraction $\numretraction_\Omega$ from \eqref{eq:num_retraction} is in fact invertible so that \eqref{eq:num_transport} is well-defined. From the numerical point of view, a vector field $\mathcal{W}$ defined on $\Omega$ is represented by a vector of nodal values since it is discretized by piecewise linear Lagrange elements. Equation~\eqref{eq:num_transport} then states that the transported vector field $\numtransport_\vectorfield \mathcal{W}$ on the deformed domain $\numretraction_\Omega(\vectorfield)$ is represented by the same vector of nodal values as the original vector field $\mathcal{W}$, where only the position of the corresponding mesh nodes is changed according to \eqref{eq:num_retraction}. Note, that the retraction given by \eqref{eq:num_retraction} and the vector transport given by \eqref{eq:num_transport} are used in \cite{Schulz2016Efficient, Schulz2015Structured, Schulz2016Computational, Geiersbach2020Stochastic} for the numerical realization of the Riemannian framework from Section~\ref{ssec:riemannian}, and we refer the reader to these publications for further details.

The previously described discretizations lead to finite-dimensional (non-)linear systems, for the state and adjoint systems as well as the gradient deformation problem \eqref{eq:def_gradient_defo}, whose solution we briefly discuss in the following. All nonlinear systems are solved by a damped Newton method based on the natural monotonicity criterion from \cite[Chapter 3.3]{Deuflhard2011Newton} with a backtracking line search, giving rise to a sequence of linear systems. For the solution of these and all other linear systems, we use the direct solver MUMPS from the library PETSc \cite{Balay2020PETSc}.

We solve each of the shape optimization problems considered in the subsequent sections numerically using the five NCG methods from Section~\ref{ssec:algo_numerical}. Additionally, we solve each problem with the gradient descent and the L-BFGS methods, where we consider a memory size of one, three, and five for the latter. This enables a detailed comparison of the NCG methods to already established gradient-based shape optimization methods. Note, that from the perspective of memory requirements, the gradient descent and NCG methods are comparable as the latter only need to store one or two additional vectors compared to the former. As remarked in Section~\ref{ssec:ncg_finite}, the L-BFGS~$m$ method needs $2m$ additional vectors of storage over the gradient descent method, which can be prohibitive for very large-scale problems as discussed in \cite{Kelley1999Iterative}. In particular, only the L-BFGS~1 method, which needs to store two additional vectors, is comparable to the NCG and gradient descent methods regarding their memory requirements.

Let us briefly discuss the parameters used for solving the optimization problems. As stopping tolerance we use $\texttt{tol} = \num{5e-4}$ in all cases, which is rather restrictive and ensures (numerical) convergence to a stationary point or local minimizer. For the Armijo line search, we always use $\armijoparam = \num{1e-4}$ and $\stepparam = \nicefrac{1}{2}$, as suggested in \cite{Nocedal2006Numerical, Kelley1999Iterative}. The remaining parameters for the optimization algorithms differ slightly between the problems so that we specify them at the relevant positions below. Note, that the methods only differ in the fact that the initial trial step size for the L-BFGS methods is chosen to be \num{1} in case of a non-empty memory, the remaining parameters are identical for all methods so that a comparison is feasible. Finally, we remark that our implementation of all numerical experiments considered in this paper is available as open-source code on GitHub \cite{Blauth2021NCGCode}.

We present our numerical results regarding the comparison of the methods in the same way for all problems, which we briefly describe in the following. We visualize the history of the methods, i.e., the evolution of the cost functional $\reducedcostfunction(\Omega\iidx{k})$ and relative shape gradient norm $\nicefrac{\norm{\gradientdefo\iidx{k}}{a_{\Omega\iidx{k}}}}{\norm{\gradientdefo\iidx{0}}{a_{\Omega\iidx{0}}}}$, over the optimization. For the sake of better readability, we only show the gradient descent, L-BFGS~5, and NCG methods in these figures and exclude the L-BFGS~3 and~1 methods. We highlight the gradient descent and L-BFGS~5 methods together with the NCG method that performed best by plotting their history in opaque colors, whereas we use transparent colors for the remaining NCG methods. An example for such a figure is given by Figure~\ref{fig:history_poisson}. 

Moreover, we tabulate the amount of iterations the methods require until they first reach a tolerance $\tau$ of
\begin{equation*}
	\tau \in \Set{\num{1e-1},\ \num{5e-2},\ \num{1e-2},\ \num{5e-3},\ \num{1e-3},\ \num{5e-4}}.
\end{equation*}
This allows us to compare how efficient the algorithms are for different tolerances, e.g., if one would want to employ the methods using a less restrictive tolerance. Note, that the main computational cost of Algorithm~\ref{algo:transformed} comes from the solution of the PDEs corresponding to the state and adjoint system as well as the linear elasticity equation, which is used to determine the gradient deformation: Each of these PDEs must be solved once per iteration of the algorithm in order to determine the current gradient deformation. Additionally, the state system possibly has to be solved several times to compute a feasible step size in the Armijo line search. The cost of computing the search direction for the NCG and L-BFGS methods, however, is negligible compared to the cost of the PDE solves. Hence, the number of PDE solves an algorithm performs is a good indicator for its computational cost. For this reason, we also state the number of solves for the state and adjoint systems required by the methods. Note, that the number of solves depicted throughout this section corresponds to the desired tolerance $\texttt{tol} = \num{5e-4}$ in case the respective method converged successfully, or to the number of solves after $k_\textrm{max}$ iterations if the respective method failed to converge. An example for such a table is given by Table~\ref{tab:poisson}.

In the end, we also show plots of the optimized geometries where we use the optimized geometry obtained with the L-BFGS~5 method as reference. Note, that for the sake of brevity, we only show the geometries obtained by the gradient descent method and one or two NCG methods.

\subsection{A Shape Optimization Problem Constrained by a Poisson Equation}
\label{ssec:poisson}

\begin{figure}[!b]
	\centering
	\includegraphics[width=0.75\textwidth]{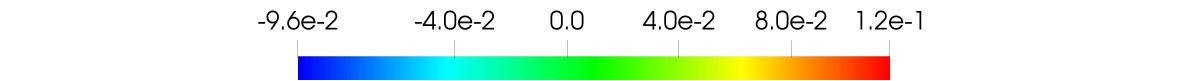}
	\begin{subfigure}[b]{0.49\textwidth}
		\centering
		\includegraphics[width=0.6\textwidth]{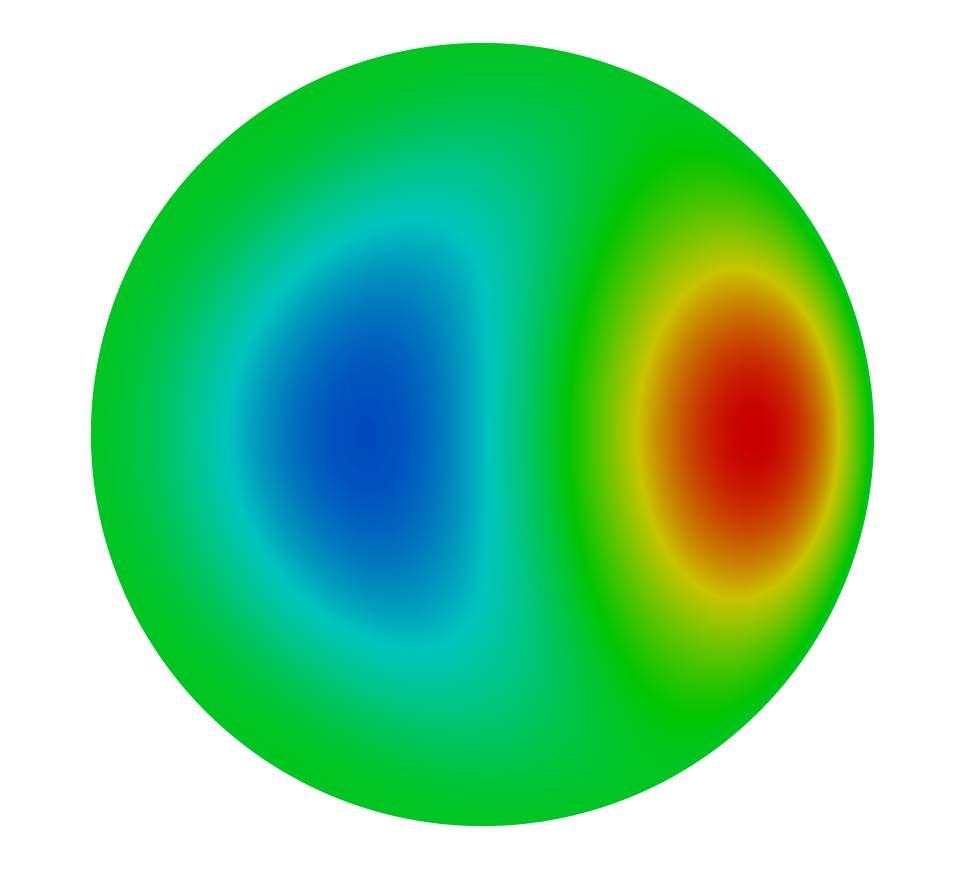}
		\caption{Initial geometry.}
	\end{subfigure}
	\hfil
	\begin{subfigure}[b]{0.49\textwidth}
		\centering
		\includegraphics[width=0.6\textwidth]{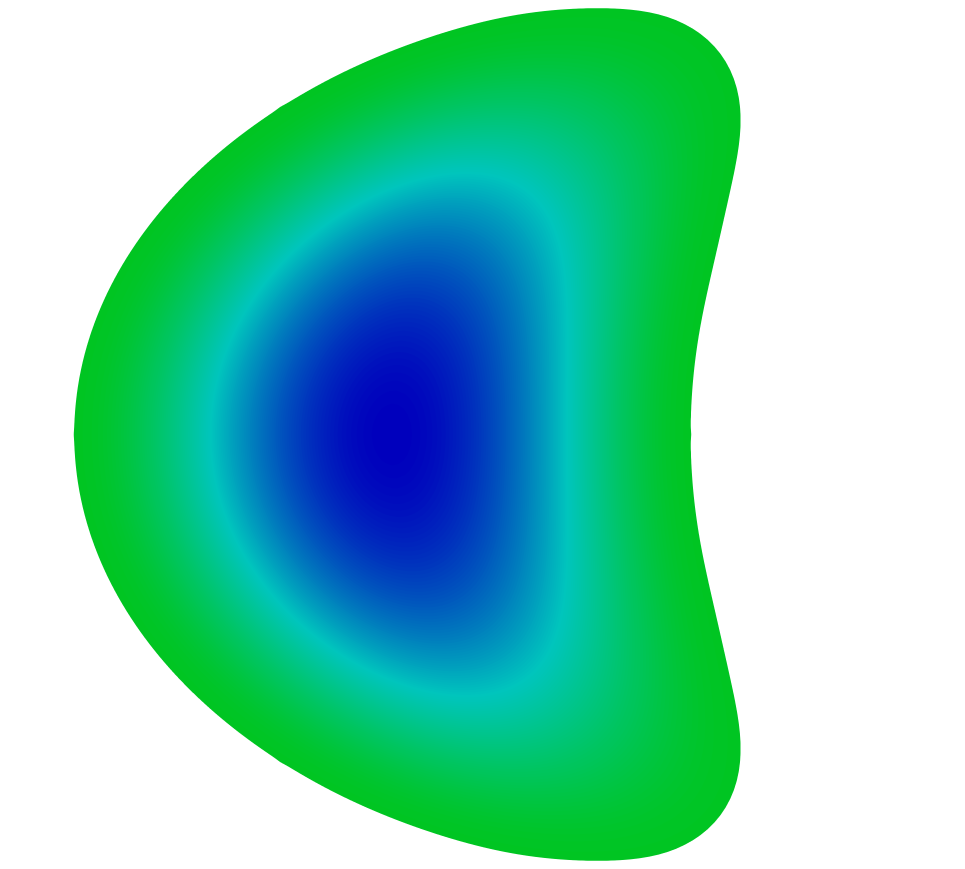}
		\caption{Optimized geometry.}
	\end{subfigure}
	\caption{State variable $u$ for the Poisson problem \eqref{eq:problem_poisson} on the initial and optimized geometries, obtained by the L-BFGS~5 method.}
	\label{fig:state_poisson}
\end{figure}

\begin{table}[t]
	\centering
	{\footnotesize
		\rowcolors{2}{\tablegray}{white}
		\setlength{\tabcolsep}{1em}
		\begin{tabular}{l r}
			\toprule
			parameter & value \\
			\midrule
			initial step size $t\iidx{0}$ & 1.0 \\
			%
			%
			maximum number of iterations $k_{\text{max}}$ & 50 \\
			iterations for NCG restart $\cgiter$ & $\infty$ \\
			tolerance for NCG restart $\cgtol$ & $\infty$ \\
			\midrule
			first Lam\'e parameter $\lamefirst$ & 1.429 \\
			second Lam\'e parameter  $\lamesecond$ & 0.357 \\
			damping parameter $\damping$ & 0.2\\
			\bottomrule
		\end{tabular}
		\caption{Parameters for Algorithm~\ref{algo:transformed} for the Poisson problem \eqref{eq:problem_poisson}.}
		\label{tab:parameters_poisson}
	}
\end{table}

The first problem we consider is the model problem taken from \cite{Etling2020First}, which we used in Section~\ref{ssec:shape_calculus}. It is given by
\begin{equation}
	\label{eq:problem_poisson}
	\min_{\Omega \in \admissiblegeom}\ \costfunction(\Omega, \state) = \integral{\Omega} \state \dmeas{x} \qquad \text{ s.t. } \qquad -\laplace \state = f \text{ in } \Omega, \quad \state = 0 \text{ on } \Gamma,
\end{equation}
where we have
\begin{equation*}
	\admissiblegeom = \Set{\Omega \subset \R^d | \Omega \subset \holdall},
\end{equation*}
for some bounded hold-all domain $\holdall$. The corresponding adjoint equation for problem \eqref{eq:problem_poisson} is given in \eqref{eq:weak_adjoint_poisson} and the volume formulation of the shape derivative is stated in \eqref{eq:vol_sd_poisson}. 

We consider this problem in two dimensions, following the numerical optimization carried out in \cite{Etling2020First}. As initial geometry we use the unit circle $S^1 \subset \R^2$ which we discretize using a uniform mesh consisting of 7651 nodes and 15000 triangles, and for the discretization of the Poisson equation we use piecewise linear Lagrange elements. For the right-hand side $f$ we use 
\begin{equation*}
	f(x) = 2.5 \left( x_1 + 0.4 - x_2^2 \right)^2 + x_1^2 + x_2^2 - 1.
\end{equation*}
Furthermore, we do not consider additional geometric constraints and, thus, use $\holdall = \R^2$. A plot of the state variable $u$ on the initial and optimized domains, obtained by the L-BFGS~5 method, can be found in Figure~\ref{fig:state_poisson}. The parameters used for the optimization algorithm and for the bilinear form $a_\Omega$ are based on the ones used in \cite{Etling2020First} and are summarized in Table~\ref{tab:parameters_poisson}.

\begin{figure}[b]
	\centering
	\begin{subfigure}{0.49\textwidth}
		\includegraphics[width=\textwidth]{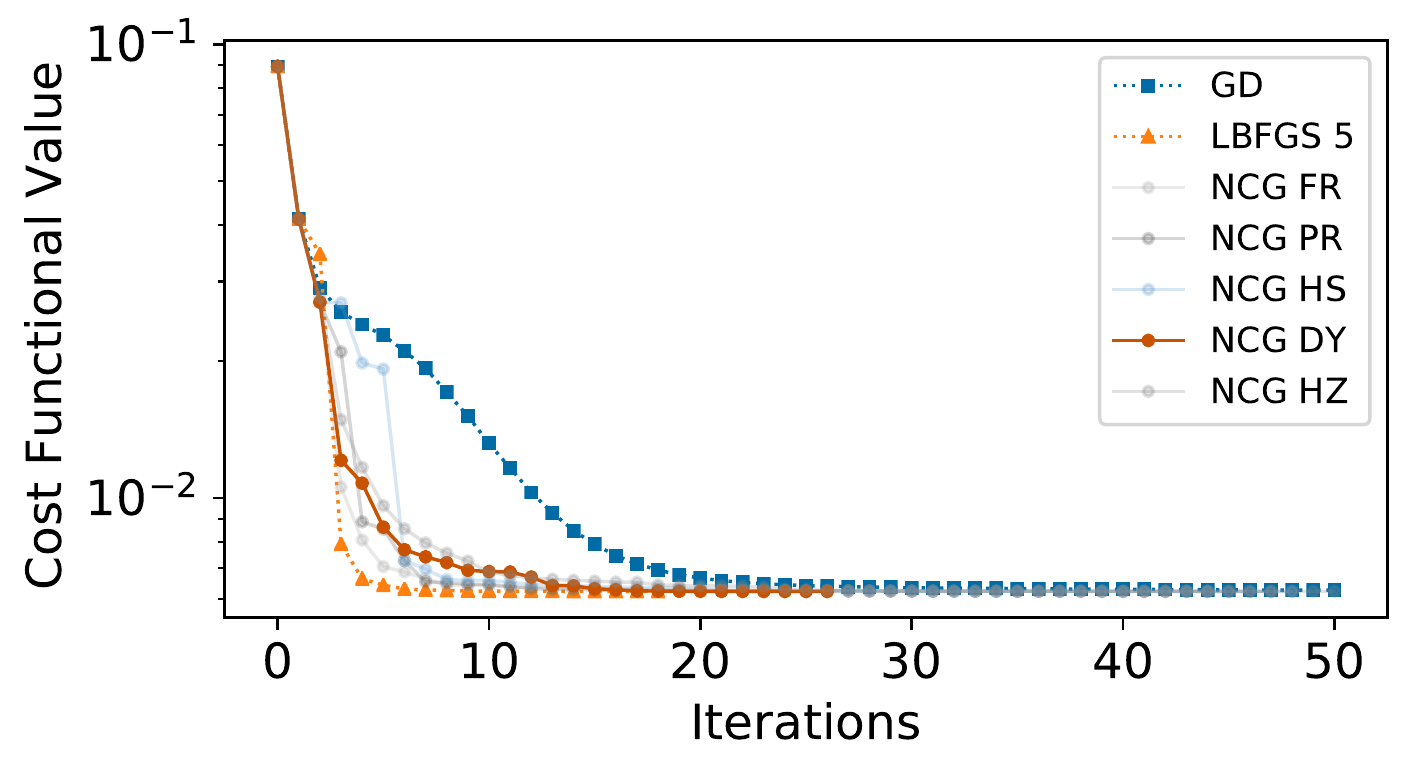}
		\caption{History of the cost functional (shifted by $+0.1$).}
	\end{subfigure}
	\hfil
	\begin{subfigure}{0.49\textwidth}
		\includegraphics[width=\textwidth]{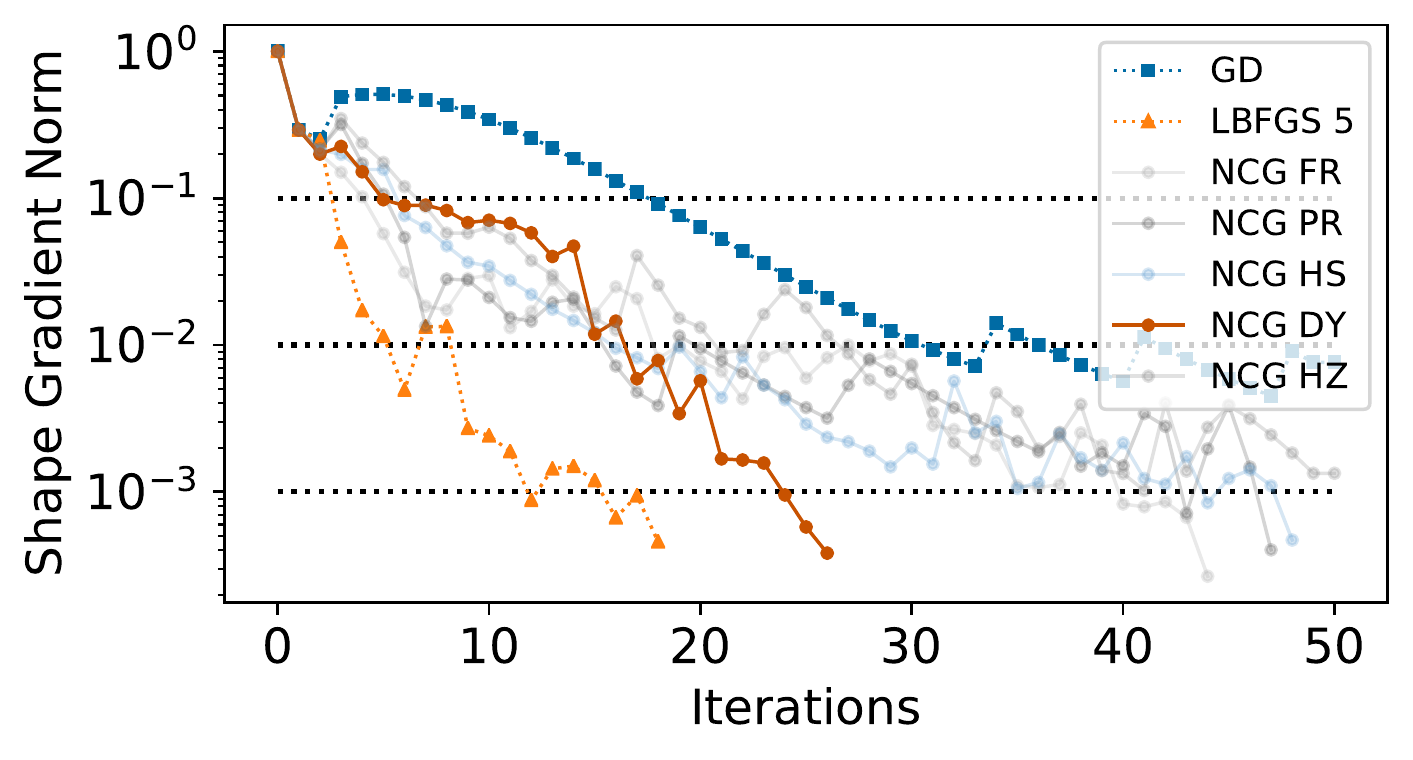}
		\caption{History of the relative shape gradient norm.}
	\end{subfigure}
	\caption{History of the optimization algorithms for the Poisson problem \eqref{eq:problem_poisson}.}
	\label{fig:history_poisson}
\end{figure}

\begin{table}[b]
	\centering
	{\footnotesize
		\rowcolors{2}{\tablegray}{white}
		\setlength{\tabcolsep}{1em}
		\begin{tabular}{c r r r r r r c c}
			\toprule
			tol & \num{1e-01} & \num{5e-02} & \num{1e-02} & \num{5e-03} & \num{1e-03} & \num{5e-04}  &  & state / adjoint solves \\
			\midrule
			GD & 18 & 22 & 31 & 47 & - & - &  & 101 / 50 \\
			\midrule
			L-BFGS 1 & 4 & 5 & 13 & 19 & 28 & 36 &  & 47 / 37 \\
			L-BFGS 3 & 3 & 4 & 6 & 11 & 16 & 22 &  & 29 / 23 \\
			L-BFGS 5 & 3 & 4 & 6 & 6 & 12 & 18 &  & 22 / 19 \\
			\midrule
			CG FR & 5 & 6 & 18 & 22 & 40 & 44 &  & 88 / 45 \\
			CG PR & 6 & 7 & 16 & 17 & 43 & 47 &  & 95 / 48 \\
			CG HS & 6 & 8 & 16 & 21 & 44 & 48 &  & 97 / 49 \\
			CG DY & 5 & 13 & 17 & 19 & 24 & 26 &  & 52 / 27 \\
			CG HZ & 7 & 12 & 21 & 29 & - & - &  & 101 / 50 \\
			\bottomrule
		\end{tabular}
		\caption{Amount of iterations required to reach a prescribed tolerance for the Poisson problem \eqref{eq:problem_poisson}.}
		\label{tab:poisson}
	}
\end{table}

As discussed in Section~\ref{ssec:implementation}, the history and performance of the methods are depicted in Figure~\ref{fig:history_poisson} and Table~\ref{tab:poisson}. 
From the results shown there, we observe that the NCG methods work very well. Each of the NCG methods performs significantly better than the gradient descent method, reaching the investigated tolerances in about half the iterations, and, except for the Hager-Zhang NCG method, all of them reach the desired tolerance of \num{5e-4}. This is also reflected in the evolution of the cost functional which decreases considerably slower for the gradient descent method than for the NCG methods. In particular, the Dai-Yuan NCG method works very well and even outperforms the L-BFGS~1 method. Overall, the performance of the NCG methods is slightly worse than that of the L-BFGS ones, but far better than that of the gradient descent method. Finally, we remark that the NCG methods need to solve the state system more often than the L-BFGS methods. This is due to the fact that the L-BFGS methods have a built-in scaling of the search direction, as remarked in Section~\ref{sec:algorithms}, so that the initial guess of \num{1} for the step size is almost always accepted, whereas the gradient descent and NCG methods need to perform more iterations for the computation of the step size via the Armijo line search.

\begin{figure}[t]
	\centering
	\begin{subfigure}[t]{0.325\textwidth}
		\includegraphics[width=0.8\textwidth]{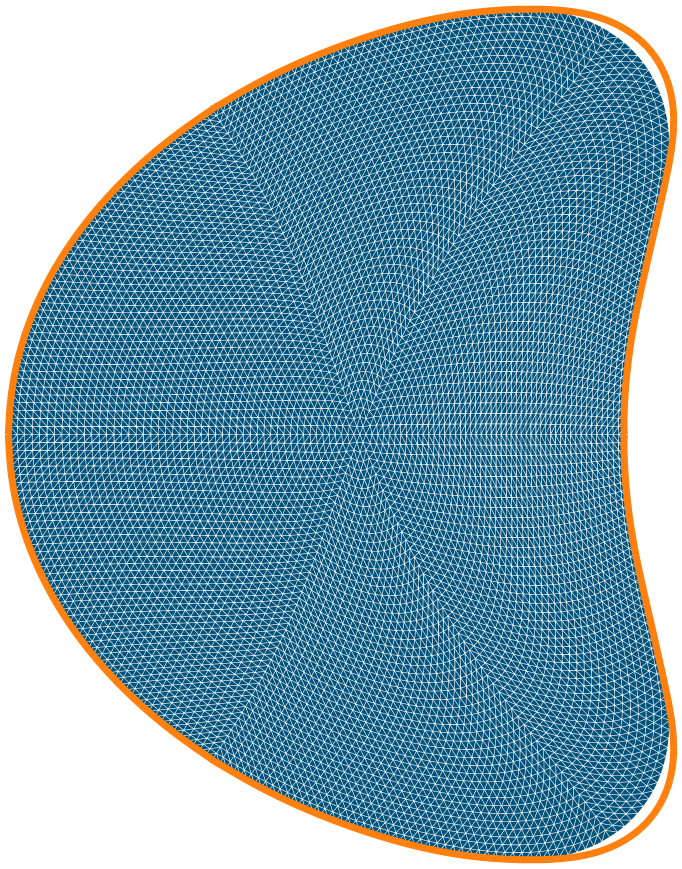}
		\caption{Gradient descent method.}
	\end{subfigure}
	\hfil
	\begin{subfigure}[t]{0.325\textwidth}
		\includegraphics[width=0.8\textwidth]{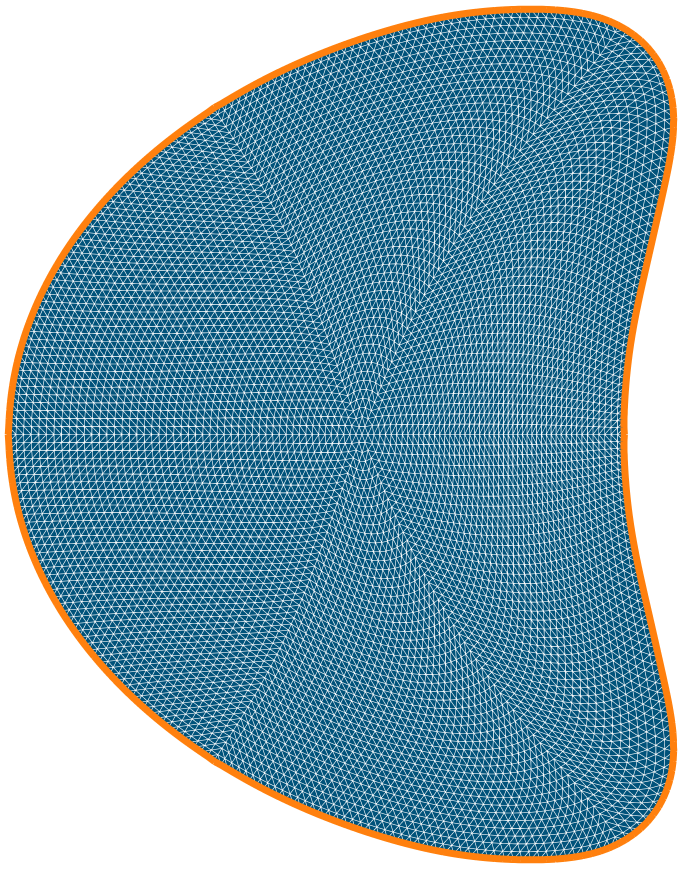}
		\caption{Fletcher-Reeves NCG method.}
	\end{subfigure}
	\hfil
	\begin{subfigure}[t]{0.325\textwidth}
		\includegraphics[width=0.8\textwidth]{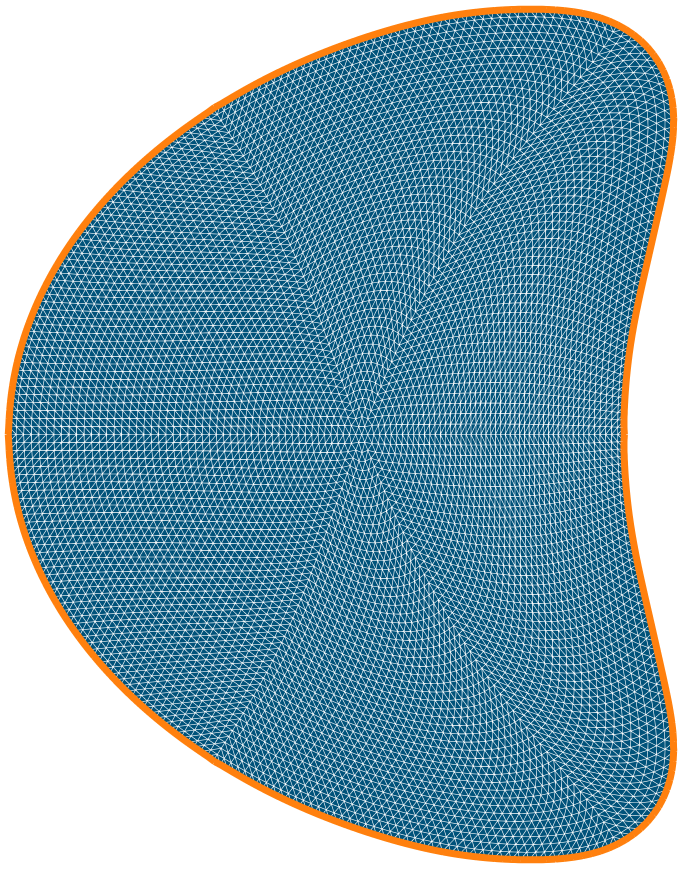}
		\caption{Polak-Ribi\`ere NCG method.}
	\end{subfigure}
	\caption{Optimized Shapes (blue) compared to the solution of the L-BFGS~5 method (orange) for the Poisson problem \eqref{eq:problem_poisson}.}
	\label{fig:optimized_poisson}
\end{figure}

The optimized geometries are shown in Figure~\ref{fig:optimized_poisson}, where we exemplarily compare the gradient descent, Fletcher-Reeves NCG, and Polak-Ribi\`ere NCG methods with the L-BFGS~5 method as reference. The optimized geometries are all rather similar, however, we observe that there is still a visible difference between the geometry obtained by the gradient descent and the reference geometry from the L-BFGS~5 method. This is not the case anymore for both NCG methods depicted there, where the geometries are in perfect agreement with the reference one.

\subsection{Shape Identification in Electrical Impedance Tomography}
\label{ssec:eit}

For our second problem, we consider an inverse problem in electrical impedance tomography, which is based on the ones investigated, e.g., in \cite{Laurain2016Distributed, Schulz2016Efficient, Hintermueller2008Electrical}. For this problem, we consider the hold-all domain $\holdall$ as the domain of the PDE, which is divided into the subdomains $\Omega\subin$ and $\Omega\subout$. Each of these subdomains represents a different material, which we assume to have different, and constant, electric conductivities $\conductivity\subin \in \R^+$ and $\conductivity\subout \in \R^+$, respectively. The goal of this problem is to identify the shape of the interior object $\Omega\subin$ from measurements of the electric potential at the boundary. To do so, usually several experiments are carried out. The electric potential $\state_i$, where the index $i=1,\dots, M$ denotes the number of the experiment, can then be modeled via the following system
\begin{equation}
	\label{eq:pde_eit}
	\begin{alignedat}{2}
		-\conductivity\subin \laplace \state\subin_i &= 0 \quad &&\text{ in } \Omega\subin, \\
		-\conductivity\subout \laplace \state\subout_i &= 0 \quad &&\text{ in } \Omega\subout, \\
		\conductivity\subout \partial_\normal \state\subout_i &= f_i \quad &&\text{ on } \partial\holdall, \\
		\state\subout_i &= \state\subin_i \quad &&\text{ on } \Gamma,\\
		\conductivity\subout \partial_{\normal\subin} \state\subout_i &= \conductivity\subin \partial_{\normal\subin} \state\subin_i \quad &&\text{ on } \Gamma, \\
		\integral{\partial\holdall} \state\subout_i \dmeas{s} &= 0,
	\end{alignedat}
\end{equation}
where $f_i$ denotes the electric current applied on the outer boundary $\partial\holdall$ and $\normal\subin$ is the unit outer normal vector for $\Omega\subin$.
Note, that the fourth and fifth equations of \eqref{eq:pde_eit} are transmission conditions that model the fact that the electric potential and the electric current are continuous over the interface $\Gamma$, and the final equation ensures the unique solvability for this Neumann problem. 

\begin{table}[t]
	\centering
	{\footnotesize
		\rowcolors{2}{\tablegray}{white}
		\setlength{\tabcolsep}{1em}
		\begin{tabular}{l r}
			\toprule
			parameter & value \\
			\midrule
			initial step size $t\iidx{0}$ & 1.0 \\
			%
			%
			maximum number of iterations $k_{\text{max}}$ & 50 \\
			iterations for NCG restart $\cgiter$ & $\infty$ \\
			tolerance for NCG restart $\cgtol$ & $\infty$ \\
			\midrule
			first Lam\'e parameter $\lamefirst$ & 0.0 \\
			second Lam\'e parameter  $\lamesecond$ & 1.0 \\
			damping parameter $\damping$ & 0.0\\
			\bottomrule
		\end{tabular}
		\caption{Parameters for Algorithm~\ref{algo:transformed} for the EIT problem \eqref{eq:problem_eit}.}
		\label{tab:parameters_eit}
	}
\end{table}

For the shape optimization problem, we assume that we are given measurements $m_i$ of the electric potential on $\partial\holdall$ corresponding to the currents $f_i$ and want to identify the shape of the interior object, i.e., the shape of $\Omega\subin$. To do so, we use the following shape optimization problem
\begin{equation}
	\label{eq:problem_eit}
	\min_{\Omega\subin \in \admissiblegeom}\ \costfunction(\Omega\subin, \state) = \sum_{i=1}^{M} \frac{\nu_i}{2} \integral{\partial\holdall} \left( \state_i - m_i \right)^2 \dmeas{s} \qquad \text{ s.t. } \eqref{eq:pde_eit},
\end{equation}
which aims at minimizing the $L^2(\partial\holdall)$ distance between the simulated electric potential $\state_i$ and the corresponding measurement $m_i$. Here, $\nu_i \geq 0$ are weights and the set of admissible geometries is given by
\begin{equation*}
	\admissiblegeom = \Set{\Omega\subin \subset \R^d | \Omega \subset \holdall}.
\end{equation*}
Note, that the corresponding adjoint equation and shape derivatives for \eqref{eq:problem_eit} can be found in \cite{Laurain2016Distributed, Hintermueller2008Electrical}.

\begin{figure}[b]
	\centering
	\begin{subfigure}{0.49\textwidth}
		\includegraphics[width=\textwidth]{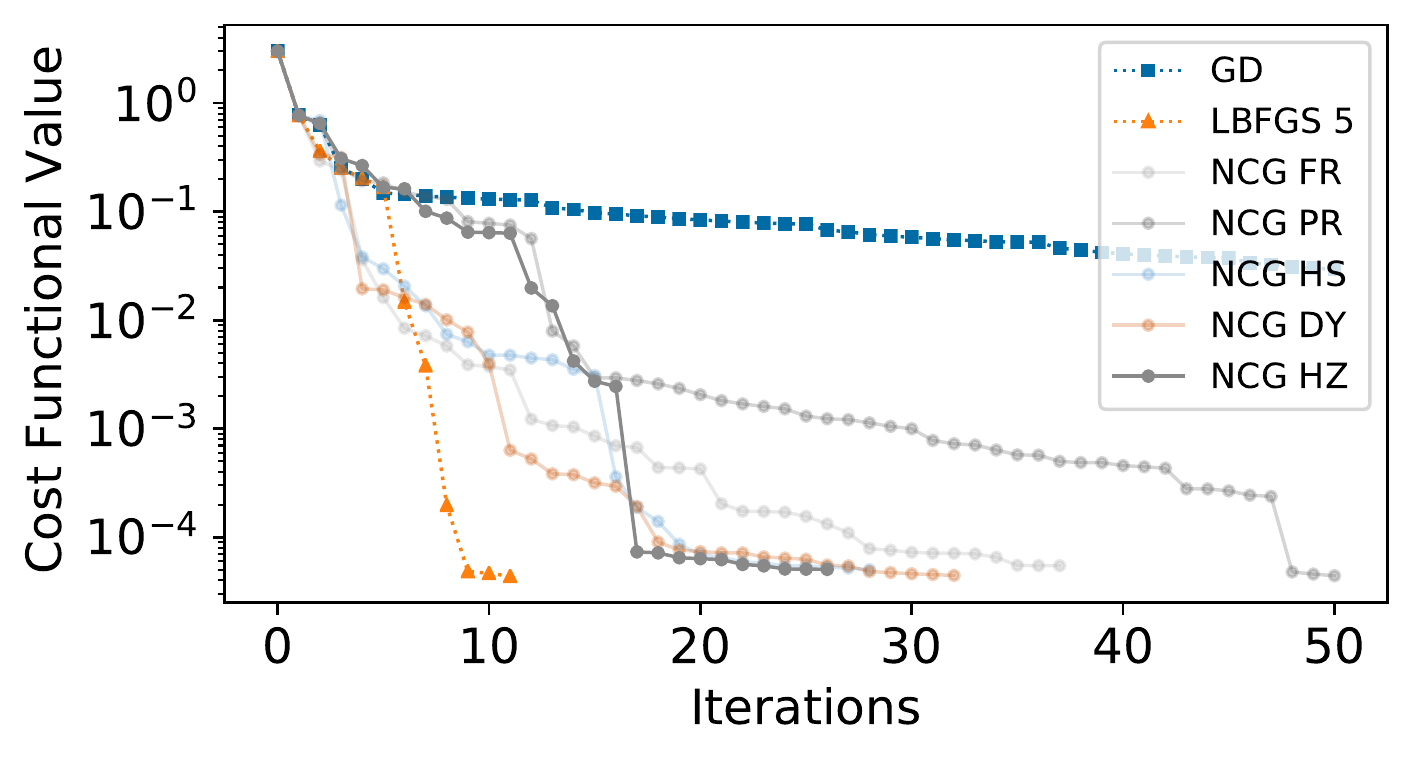}
		\caption{History of the cost functional.}
	\end{subfigure}
	\hfil
	\begin{subfigure}{0.49\textwidth}
		\includegraphics[width=\textwidth]{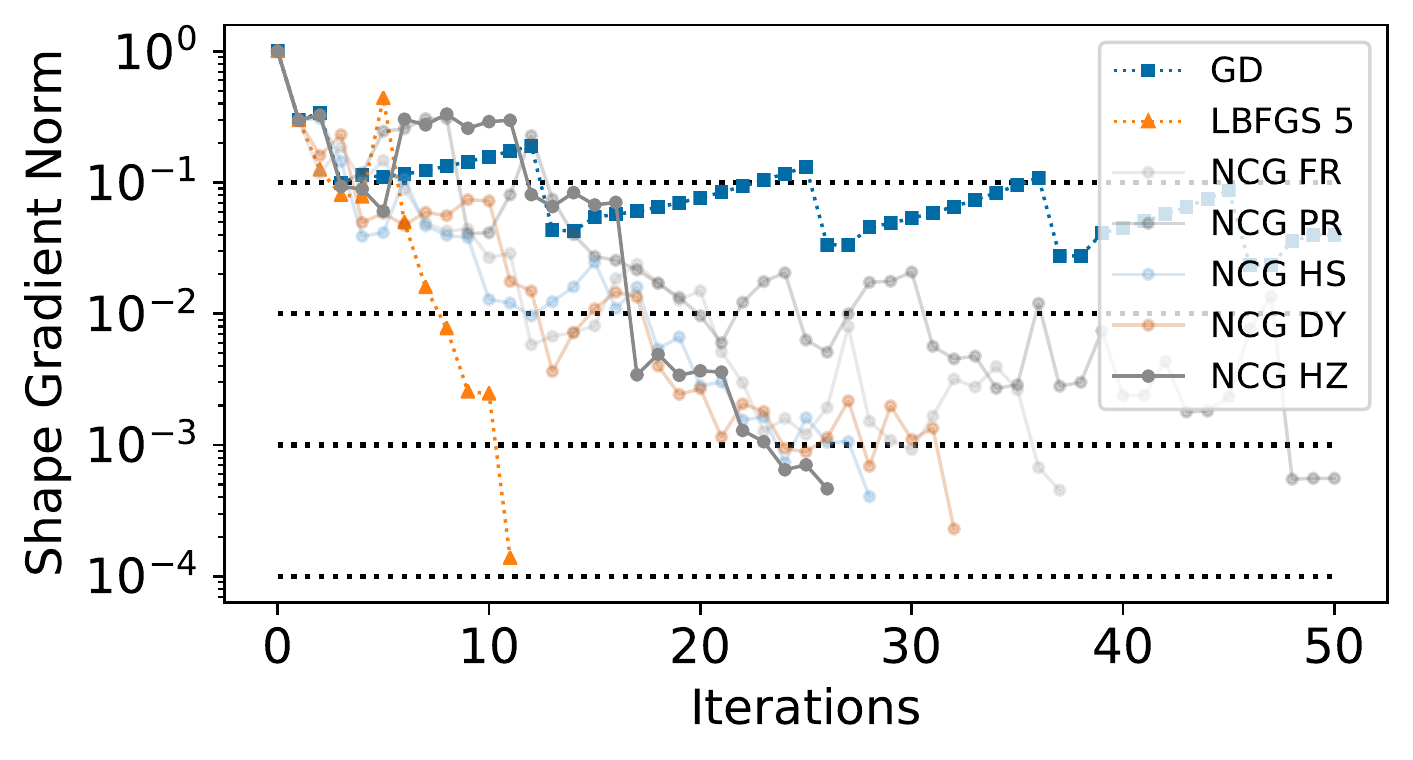}
		\caption{History of the relative shape gradient norm.}
	\end{subfigure}
	\caption{History of the optimization algorithms for the EIT problem \eqref{eq:problem_eit}.}
	\label{fig:history_eit}
\end{figure}

\begin{table}[b]
	\centering
	{\footnotesize
		\rowcolors{2}{\tablegray}{white}
		\setlength{\tabcolsep}{1em}
		\begin{tabular}{c r r r r r r c c}
			\toprule
			tol & \num{1e-01} & \num{5e-02} & \num{1e-02} & \num{5e-03} & \num{1e-03} & \num{5e-04}  &  & state / adjoint solves \\
			\midrule
			GD & 3 & 13 & - & - & - & - &  & 104 / 50 \\
			\midrule
			L-BFGS 1 & 3 & 10 & 25 & 26 & 29 & 30 &  & 39 / 31 \\
			L-BFGS 3 & 3 & 7 & 9 & 10 & 11 & 11 &  & 18 / 12 \\
			L-BFGS 5 & 3 & 6 & 8 & 9 & 11 & 11 &  & 15 / 12 \\
			\midrule
			CG FR & 6 & 7 & 12 & 22 & 30 & 37 &  & 76 / 38 \\
			CG PR & 3 & 9 & 20 & 32 & 48 & - &  & 102 / 50 \\
			CG HS & 4 & 4 & 12 & 20 & 24 & 28 &  & 56 / 29 \\
			CG DY & 4 & 4 & 13 & 13 & 24 & 32 &  & 67 / 33 \\
			CG HZ & 3 & 17 & 17 & 17 & 24 & 26 &  & 53 / 27 \\
			\bottomrule
		\end{tabular}
	}
	\caption{Amount of iterations required to reach a prescribed tolerance for the EIT problem \eqref{eq:problem_eit}.}
	\label{tab:eit}
\end{table}

For the numerical realization of this problem, we choose the following setting closely based on the one of \cite{Laurain2016Distributed}. The hold-all domain is chosen as the unit square $\holdall = (0, 1)^2$, which we discretize using a uniform mesh with 6070 vertices and 11870 triangles, and for the numerical solution of the state and adjoint systems we use piecewise linear Lagrange elements. We divide the boundary of the holdall into $\partial\holdall = \Gamma^t \cup \Gamma^b \cup \Gamma^l \cup \Gamma^r$, corresponding to the top, bottom, left, and right sides of the square. We consider the case of $M=3$ measurements, where the currents $f_i$ are chosen as 
\begin{equation*}
	\begin{aligned}
		&f_1 = 1 \text{ on } \Gamma^l \cup \Gamma^r \qquad \text{ and } \qquad f_1 = -1 \text{ on } \Gamma^t \cup \Gamma^b, \\
		&f_2 = 1 \text{ on } \Gamma^l \cup \Gamma^t \qquad \text{ and } \qquad f_2 = -1 \text{ on } \Gamma^r \cup \Gamma^b, \\
		&f_3 = 1 \text{ on } \Gamma^l \cup \Gamma^b \qquad \text{ and } \qquad f_3 = -1 \text{ on } \Gamma^r \cup \Gamma^t,
	\end{aligned}
\end{equation*}
which satisfies the compatibility condition $\integral{\partial \holdall} f_i \dmeas{s} = 0$, and the electric conductivities are given by $\conductivity\subout = \num{1}$ and $\conductivity\subin = \num{10}$. The measurements are obtained by numerically solving the state system \eqref{eq:pde_eit} with a reference inner domain given by a circle with center $(0.5, 0.5)$ and radius \num{0.2}. For the initial geometry we choose a square with center $(0.5, 0.5)$ and edge length \num{0.4}. The weights $\nu_i$ for the cost functional are chosen so that each summand in $\costfunction$ has value \num{1} after solving $\eqref{eq:pde_eit}$ on the initial geometry. Finally, the parameters for Algorithm~\ref{algo:transformed} and for the choice of the bilinear form $a_\Omega$ are given in Table~\ref{tab:parameters_eit}. 


Analogously to before, we depict the history of the optimization and performance of the methods in Figure~\ref{fig:history_eit} and in Table~\ref{tab:eit}, as discussed in Section~\ref{ssec:implementation}. The results again highlight the capabilities of the NCG methods proposed in Section~\ref{sec:algorithms}. Comparing their performance with the gradient descent method, we observe that each NCG method performs significantly better. In particular, the gradient descent method is only able to reach a tolerance of \num{5e-2} for the relative gradient norm within the maximum number of iterations, whereas all NCG methods, except the Polak-Ribi\`ere one, reach the desired tolerance of \num{5e-4}. A similar trend can be seen for the cost functional, which only decreases by about 1.5 orders of magnitude for the gradient descent method, but by over 4 orders of magnitude for all other methods. Comparing the NCG methods to the L-BFGS ones, we observe that they perform very similarly to the L-BFGS~1 method, particularly the Hestenes-Stiefel and Hager-Zhang method, which need slightly fewer iterations to reach the investigated tolerances. Note, that the L-BFGS~3 and~5 methods perform best of all methods, reaching the desired tolerance after only about a third of the iterations required by the L-BFGS~1 and NCG methods. We also observe that the L-BFGS methods need less solves of the state equation, again due to the same reason as before. However, as remarked earlier, the L-BFGS methods also need significantly more memory than the NCG methods to achieve these results.

\begin{figure}[t]
	\centering
	\begin{subfigure}[t]{0.325\textwidth}
		\includegraphics[width=\textwidth]{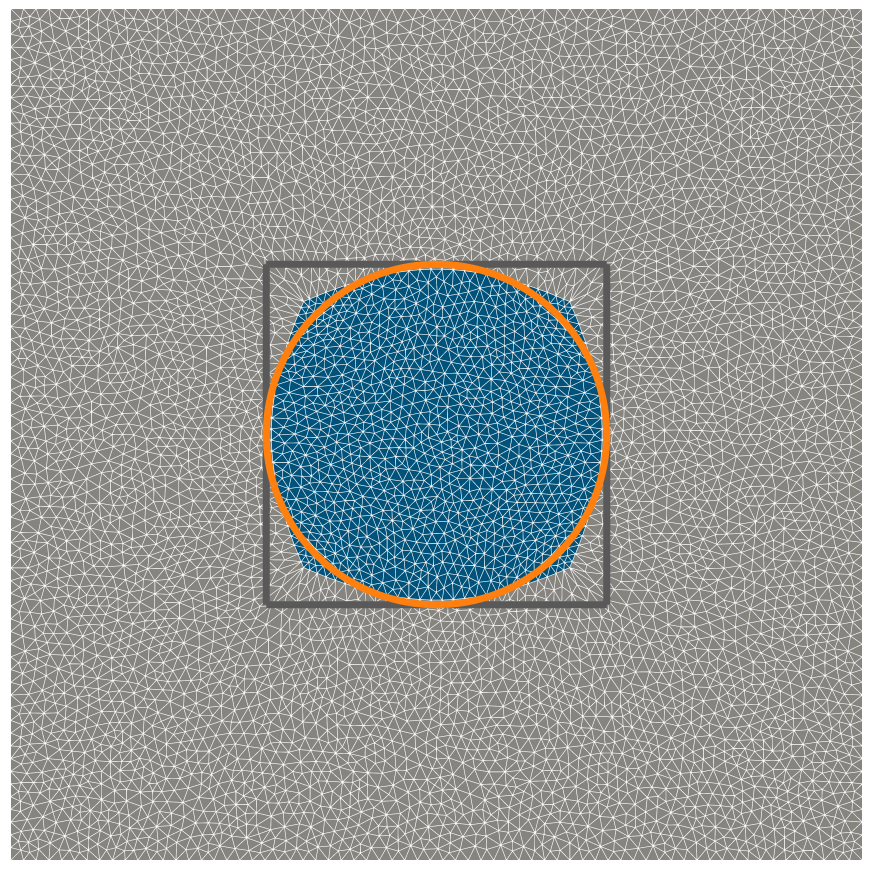}
		\caption{Gradient Descent method.}
	\end{subfigure}
	\hfil
	\begin{subfigure}[t]{0.325\textwidth}
		\includegraphics[width=\textwidth]{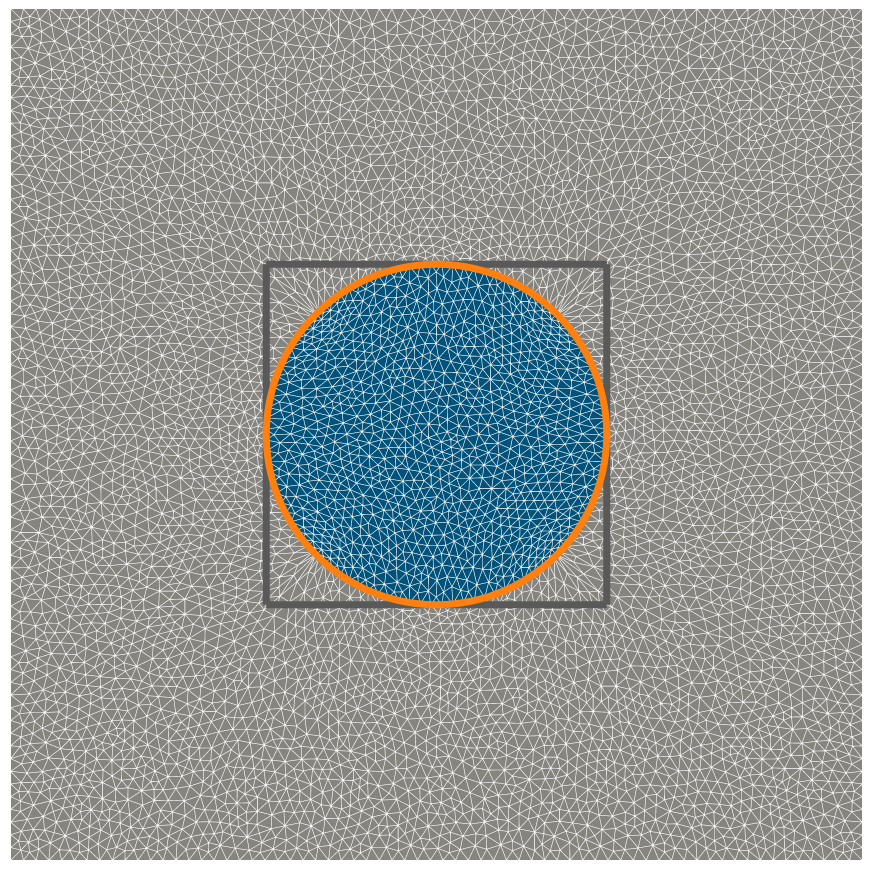}
		\caption{Hestenes-Stiefel NCG method.}
	\end{subfigure}
	\hfil
	\begin{subfigure}[t]{0.325\textwidth}
		\includegraphics[width=\textwidth]{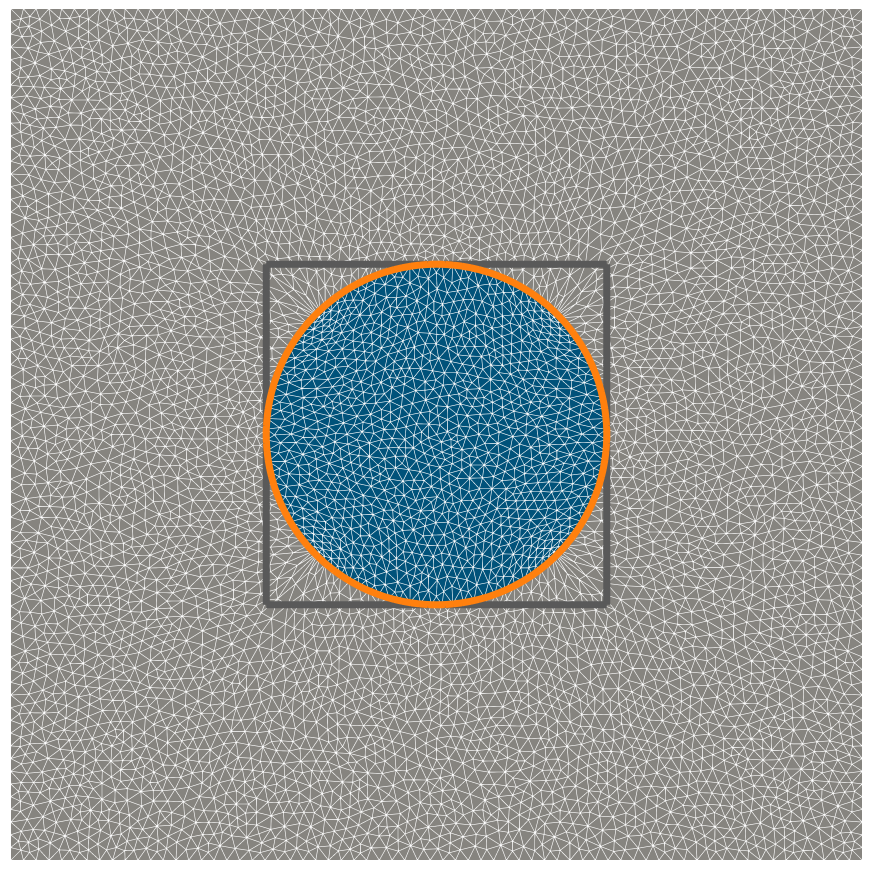}
		\caption{Hager-Zhang NCG method.}
	\end{subfigure}
	\caption{Optimized Shapes for the EIT problem \eqref{eq:problem_eit}, $\Omega\subin$ (blue), $\Omega\subout$ (light gray), together with initial (dark gray) and reference (orange) shape of $\Omega^\text{int}$.}
	\label{fig:optimized_eit}
\end{figure}

The optimized geometries obtained by the methods, shown in Figure~\ref{fig:optimized_eit} for the gradient descent method as well as the Hestenes-Stiefel and Hager-Zhang NCG methods, show a similar picture. For the gradient descent method there is still some difference between the optimized shape and the reference one, in particular, the optimized geometry still exhibits kinks corresponding to those in the initial geometry. On the other hand, both NCG methods yield accurately resolved circles, approximating the reference shape very well.

\subsection{Optimization of an Obstacle in Stokes Flow}
\label{ssec:stokes}

For the next problem, we consider the shape optimization of an obstacle enclosed in Stokes flow, where we aim at minimizing the energy dissipated by the obstacle. For this optimization to be meaningful, we need additional constraints, namely we have to fix the volume and the barycenter of the obstacle, otherwise it would either shrink arbitrarily or move out of the computational domain. This problem is investigated analytically in \cite{Pironneau1974optimum}, where the optimal shape of the obstacle is found to be the well-known ogive. Additionally, this problem is also used in, e.g., \cite{Iglesias2018Two, Dokken2020Automatic, Schulz2016Computational} for validating the performance of shape optimization algorithms. Our formulation of the problem is closely adapted from \cite{Schulz2016Computational}.

\begin{figure}[b]
	\centering
	\includegraphics[width=0.75\textwidth]{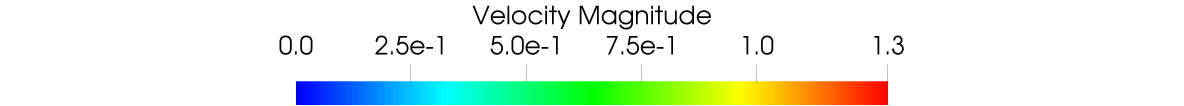}
	\begin{subfigure}{0.49\textwidth}
		\centering
		\includegraphics[width=0.95\textwidth]{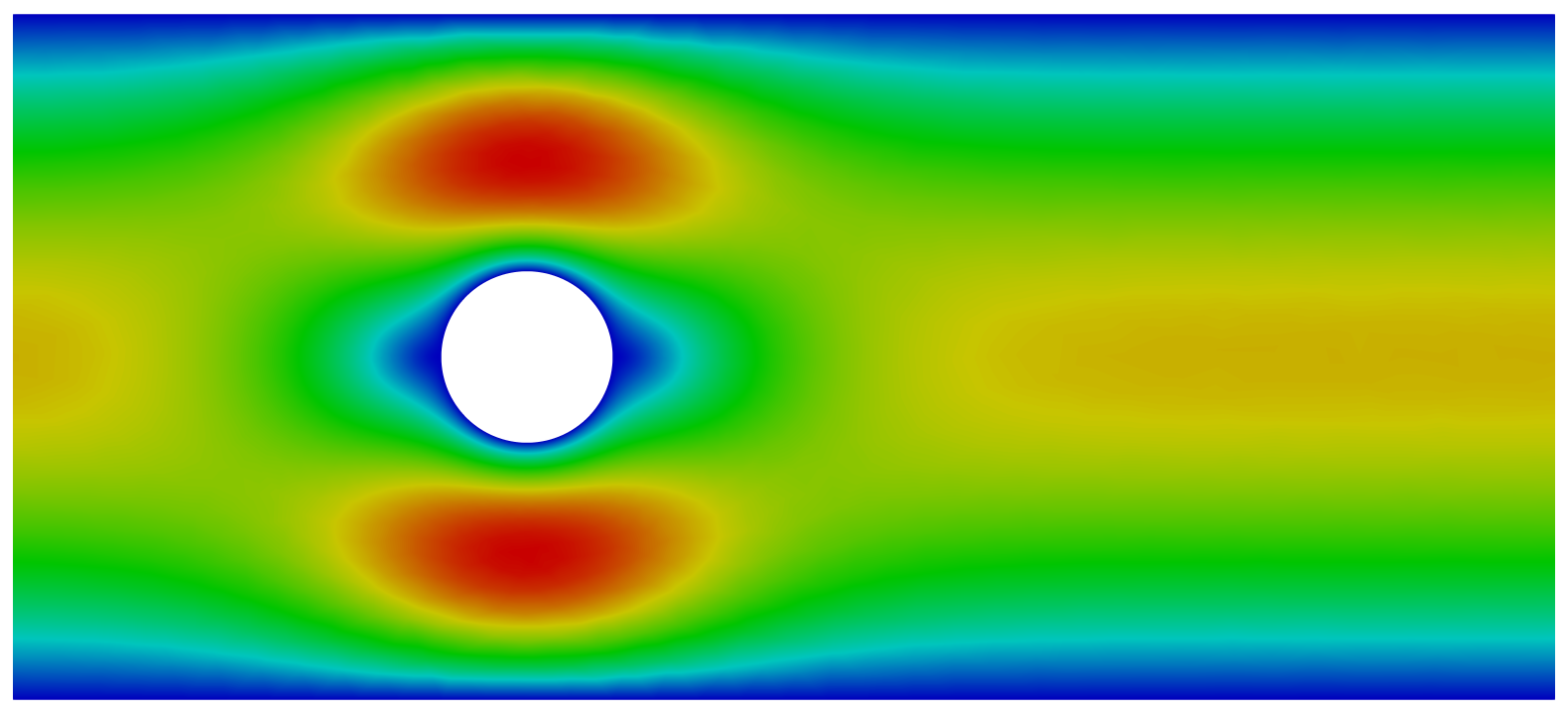}
		\caption{Initial geometry.}
	\end{subfigure}
	\hfil
	\begin{subfigure}{0.49\textwidth}
		\centering
		\includegraphics[width=0.95\textwidth]{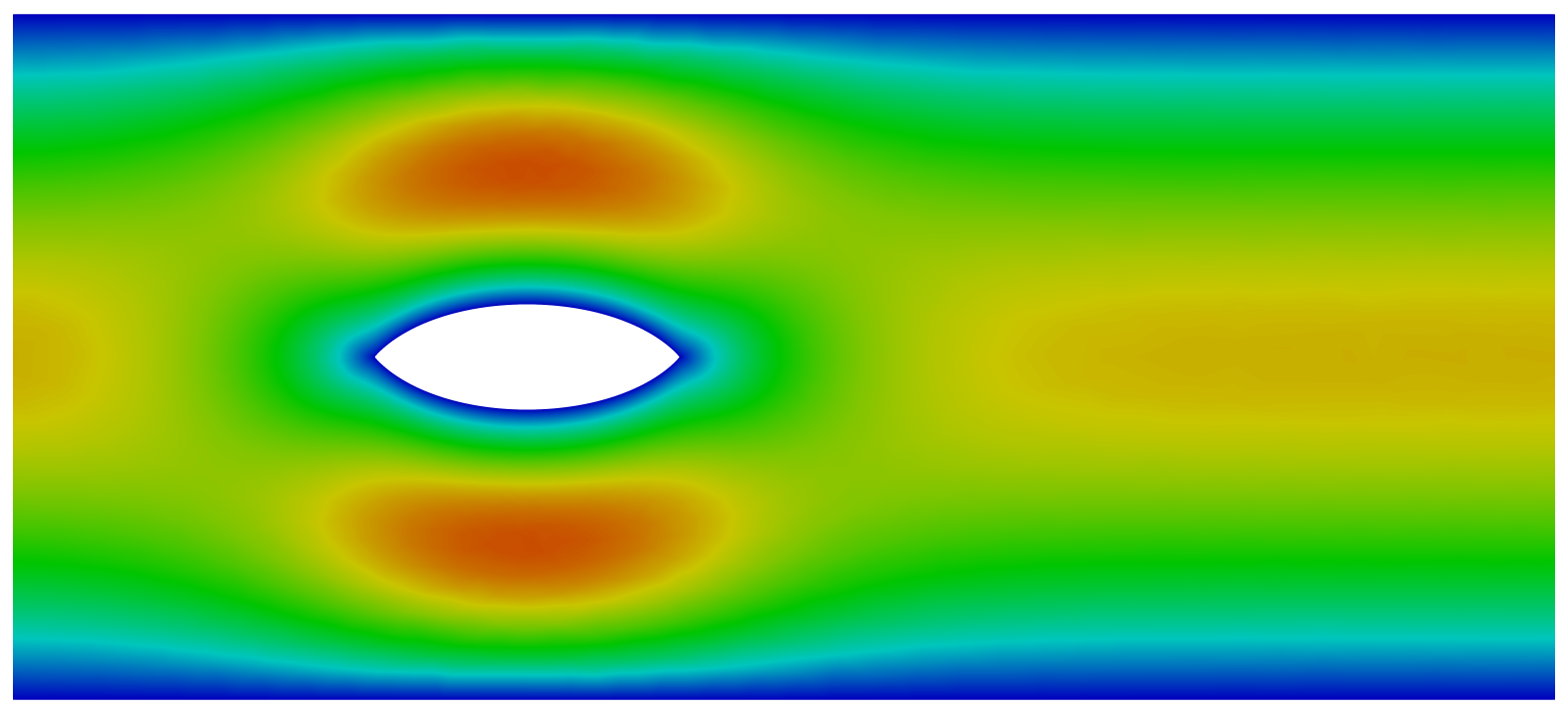}
		\caption{Optimized Geometry.}
	\end{subfigure}
	\caption{Magnitude of the velocity $\norm{\velocity}{}$ for the Stokes problem \eqref{eq:pde_stokes} on the initial and optimized geometries, obtained by the L-BFGS~5 method.}
	\label{fig:state_stokes_velocity}
\end{figure}

\begin{figure}[b]
	\centering
	\includegraphics[width=0.75\textwidth]{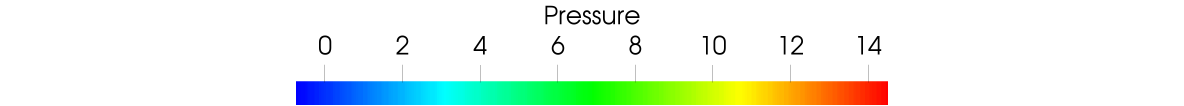}
	\begin{subfigure}{0.49\textwidth}
		\centering
		\includegraphics[width=0.95\textwidth]{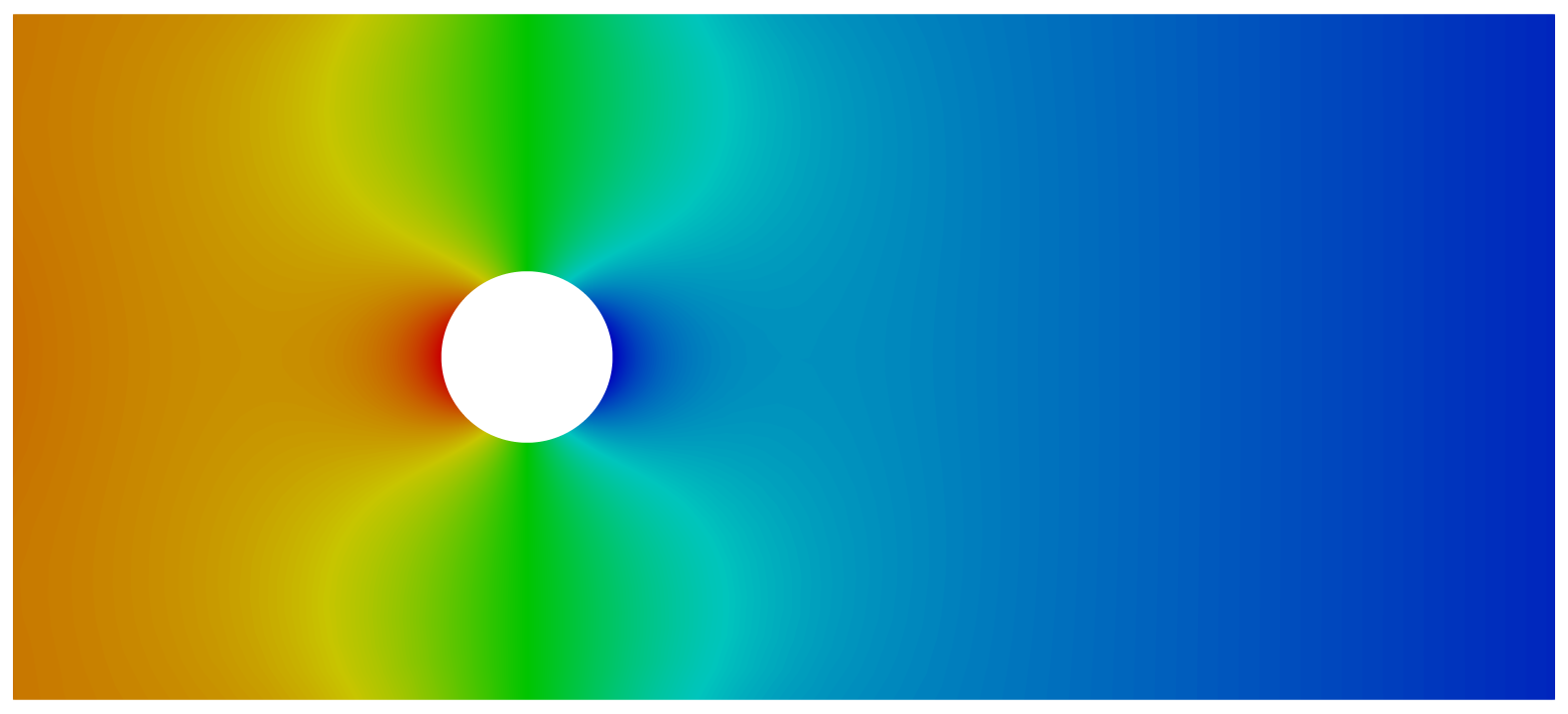}
		\caption{Initial geometry.}
	\end{subfigure}
	\hfil
	\begin{subfigure}{0.49\textwidth}
		\centering
		\includegraphics[width=0.95\textwidth]{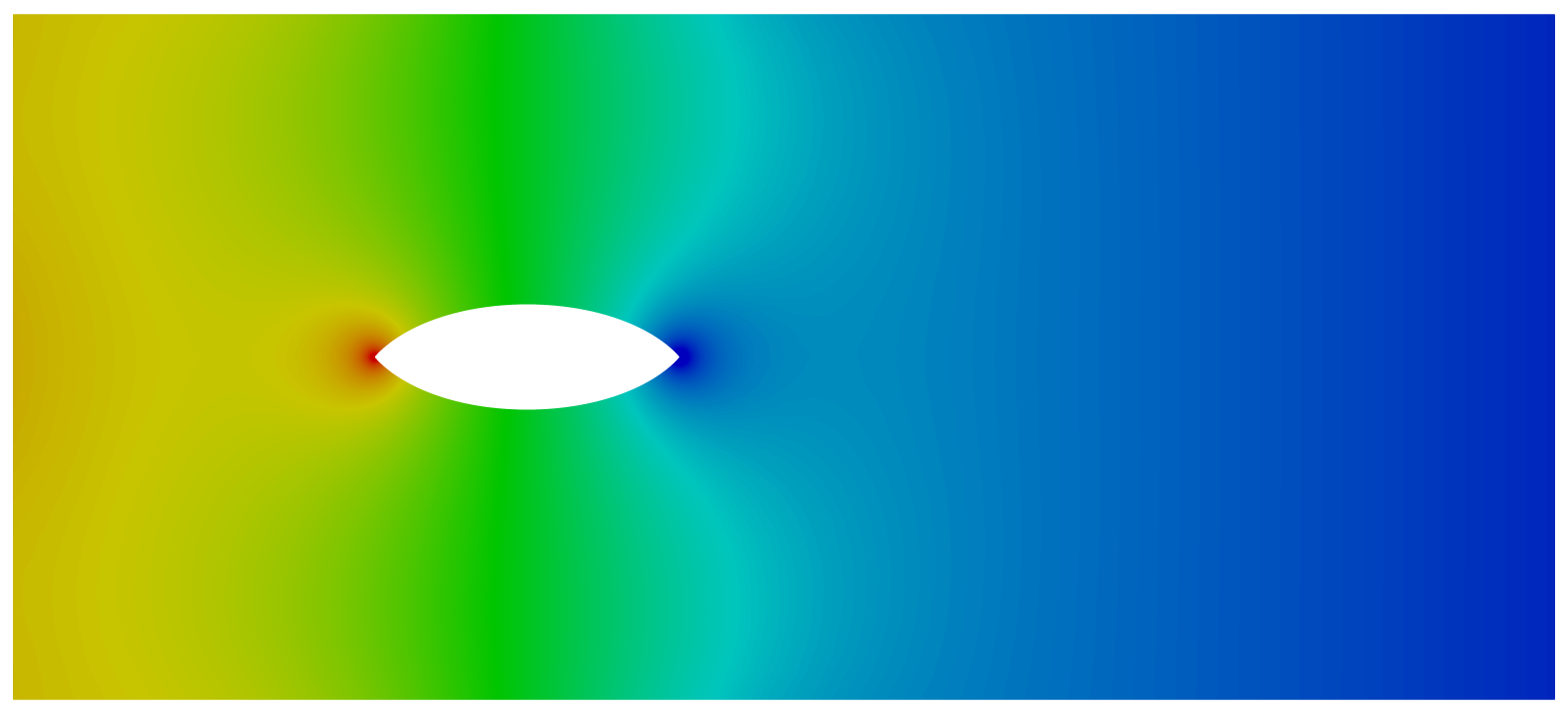}
		\caption{Optimized Geometry.}
	\end{subfigure}
	\caption{Pressure $\pressure$ for the Stokes problem \eqref{eq:pde_stokes} on the initial and optimized geometries, obtained by the L-BFGS~5 method.}
	\label{fig:state_stokes_pressure}
\end{figure}

Regarding the geometrical setup, we have the following. The domain of the flow is denoted by $\Omega$. Its exterior boundary is divided into the inlet $\Gamma\subin$, the wall boundary $\Gamma\subwall$, and the outlet $\Gamma\subout$. Moreover, we have the boundary of the obstacle $\Gamma\subobs$, which lies completely in the interior of $\Omega$ and determines the shape of the obstacle $\Omega\subobs$. The non-dimensionalized Stokes system for this setting is given by
\begin{equation}
	\label{eq:pde_stokes}
	\begin{alignedat}{2}
		- \laplace \velocity + \grad \pressure &= 0 \quad &&\text{ in } \Omega, \\
		\divergence{\velocity} &= 0 \quad &&\text{ in } \Omega, \\
		\velocity &= \velocity\subin \quad &&\text{ on } \Gamma\subin, \\
		\velocity &= 0 \quad &&\text{ on } \Gamma\subwall \cup \Gamma\subobs, \\
		\partial_\normal \velocity - \pressure \normal &= 0 \quad &&\text{ on } \Gamma\subout.
	\end{alignedat}
\end{equation}
The volume and barycenter of the obstacle are defined as 
\begin{equation*}
	\text{vol}(\Omega\subobs) = \integral{\Omega\subobs} 1 \dmeas{x} \qquad \text{ and } \qquad \text{bc}(\Omega\subobs) = \frac{1}{\text{vol}(\Omega\subobs)} \integral{\Omega\subobs} x \dmeas{x}. 
\end{equation*}
Hence, the shape optimization problem can be formulated as
\begin{equation}
	\label{eq:problem_stokes}
	\min_{\Omega \in \admissiblegeom}\ \costfunction(\Omega, \velocity) = \integral{\Omega} \norm{\grad \velocity}{F}^2 \dmeas{x} + \frac{\nu_1}{2} \left( \text{vol}(\Omega\subobs) - \text{vol}(\Omega\subobs_0) \right)^2 + \frac{\nu_2}{2} \left( \text{bc}(\Omega\subobs) - \text{bc}(\Omega\subobs_0) \right)^2 \qquad \text{ s.t. } \eqref{eq:pde_stokes},
\end{equation}
where $\norm{\cdot}{F}$ denotes the Frobenius norm and $\Omega\subobs\iidx{0}$ is the domain of the initial obstacle. Note, that we have regularized the geometrical constraints in the cost function, as in \cite{Dokken2020Automatic}. The set of admissible geometries for this problem is given by
\begin{equation*}
	\admissiblegeom = \Set{\Omega \subset \R^d | \Omega \subset \holdall,\ \Gamma\subin = \Gamma\subin\iidx{0},\ \Gamma\subwall = \Gamma\subwall\iidx{0}, \text{ and } \Gamma\subout = \Gamma\subout\iidx{0}},
\end{equation*}
for some initial flow domain $\Omega\iidx{0}$ with boundary $\Gamma\iidx{0}$, i.e., only the boundary $\Gamma\subobs$ is deformable, the remaining ones are fixed. Note, that the hold-all $\holdall$ for this example is given by $\holdall = \Omega\iidx{0} \cup \Omega\subobs\iidx{0}$. For the shape derivative of this problem we refer the reader to \cite{Schulz2016Computational}.

\begin{table}[t]
	\centering
	{\footnotesize
		\rowcolors{2}{\tablegray}{white}
		\setlength{\tabcolsep}{1em}
		\begin{tabular}{l r}
			\toprule
			parameter & value \\
			\midrule
			initial step size $t\iidx{0}$ & 1.0 \\
			%
			%
			maximum number of iterations $k_{\text{max}}$ & 250 \\
			iterations for NCG restart $\cgiter$ & $\infty$ \\
			tolerance for NCG restart $\cgtol$ & $\infty$ \\
			\midrule
			first Lam\'e parameter $\lamefirst$ & 0.0 \\
			second Lam\'e parameter  $\lamesecond$ & computed from \eqref{eq:lame_siebenborn} \\
			damping parameter $\damping$ & 0.0\\
			\bottomrule
		\end{tabular}
		\caption{Parameters for Algorithm~\ref{algo:transformed} for the Stokes problem \eqref{eq:problem_stokes}.}
		\label{tab:parameters_stokes}
	}
\end{table}

\begin{figure}[b]
	\centering
	\begin{subfigure}[t]{0.49\textwidth}
		\includegraphics[width=\textwidth]{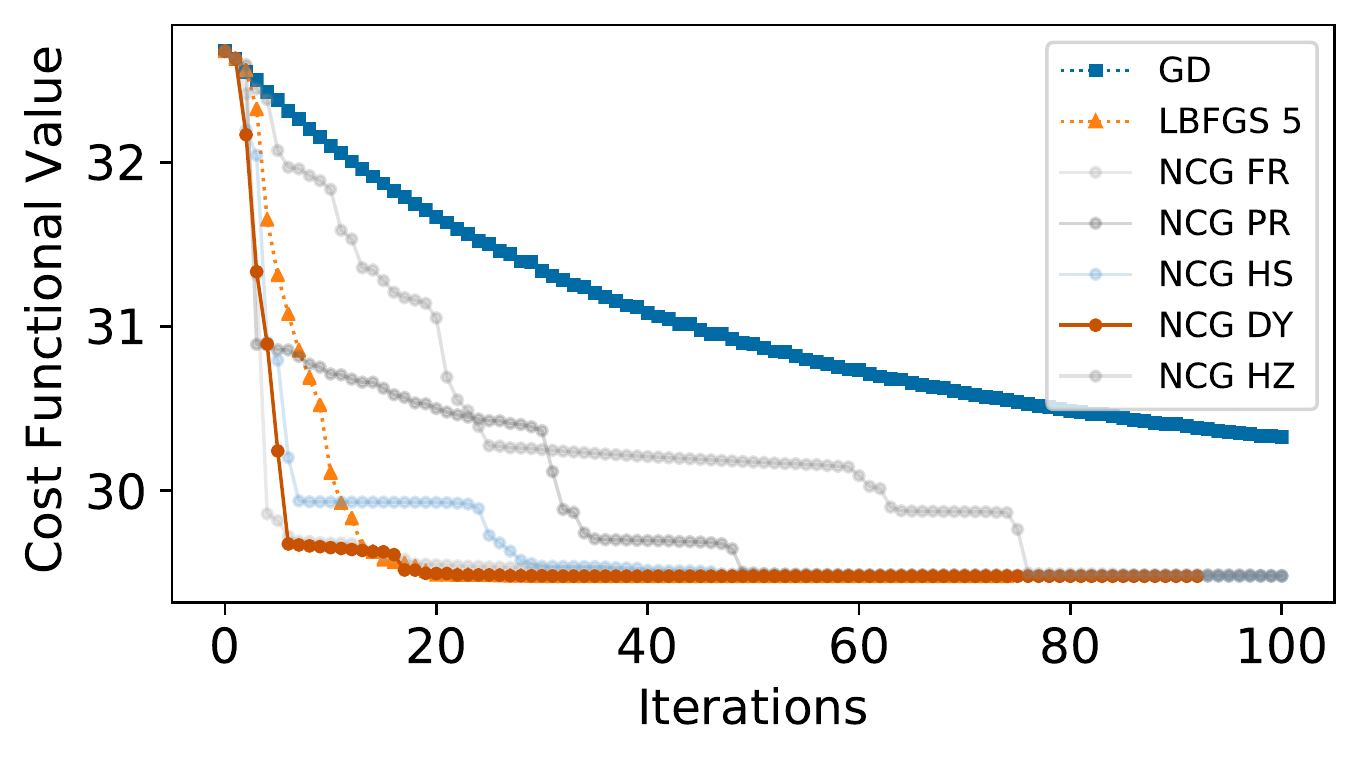}
		\caption{History of the cost functional for the first 100 iterations.}
	\end{subfigure}
	\hfil
	\begin{subfigure}[t]{0.49\textwidth}
		\includegraphics[width=\textwidth]{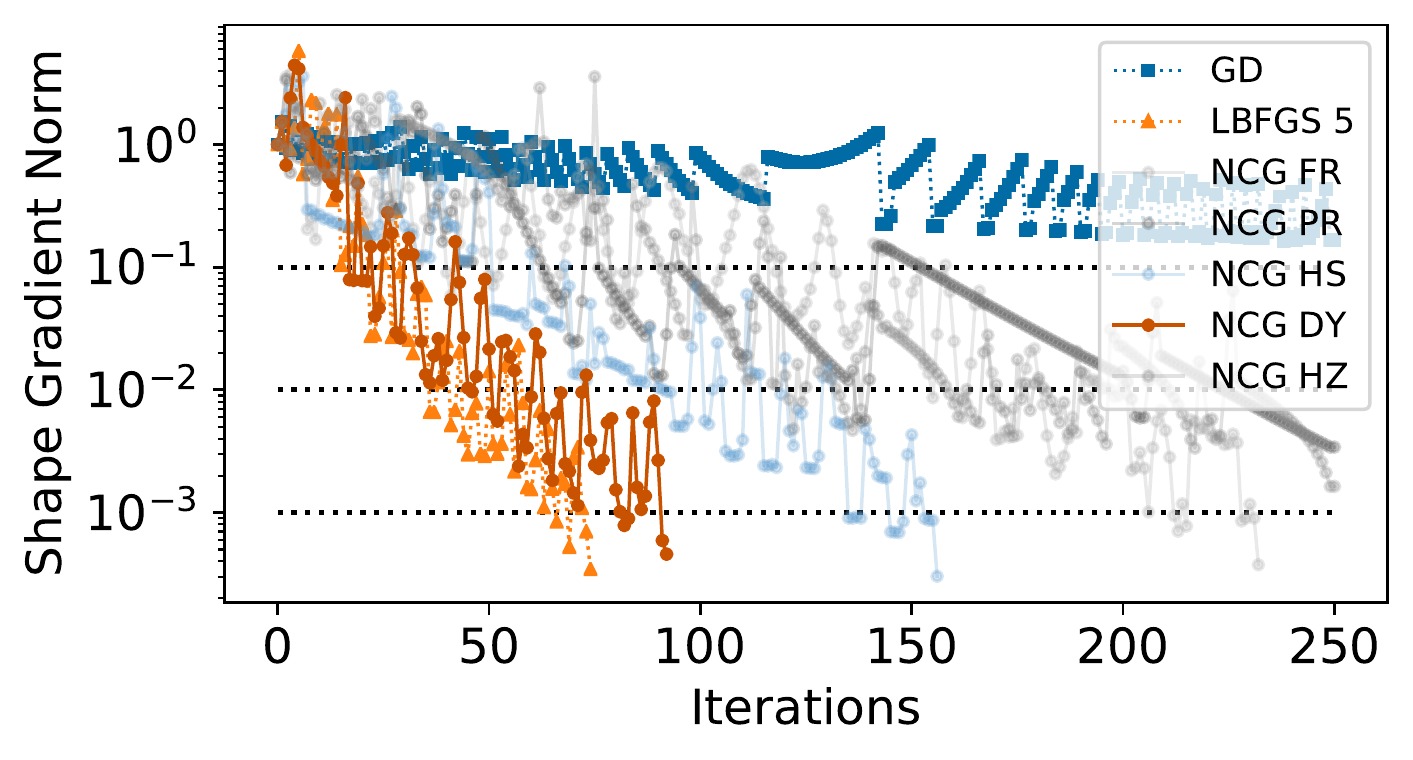}
		\caption{History of the relative shape gradient norm.}
	\end{subfigure}
	\caption{History of the optimization algorithms for the Stokes problem \eqref{eq:problem_stokes}.}
	\label{fig:history_stokes}
\end{figure}

\begin{table}[b]
	\centering
	{\footnotesize
		\rowcolors{2}{\tablegray}{white}
		\setlength{\tabcolsep}{1em}
		\begin{tabular}{c r r r r r r c c}
			\toprule
			tol & \num{1e-01} & \num{5e-02} & \num{1e-02} & \num{5e-03} & \num{1e-03} & \num{5e-04}  &  & state / adjoint solves \\
			\midrule
			GD & - & - & - & - & - & - &  & 504 / 250 \\
			\midrule
			L-BFGS 1 & 26 & 32 & 87 & 88 & 108 & 125 &  & 186 / 126 \\
			L-BFGS 3 & 28 & 30 & 70 & 76 & 112 & 112 &  & 147 / 113 \\
			L-BFGS 5 & 22 & 22 & 36 & 44 & 66 & 74 &  & 95 / 75 \\
			\midrule
			CG FR & 40 & 81 & 155 & 170 & 212 & 232 &  & 467 / 233 \\
			CG PR & 63 & 69 & 137 & 240 & - & - &  & 501 / 250 \\
			CG HS & 51 & 51 & 92 & 106 & 135 & 156 &  & 314 / 157 \\
			CG DY & 17 & 23 & 46 & 57 & 82 & 92 &  & 185 / 93 \\
			CG HZ & 79 & 80 & 121 & 122 & - & - &  & 502 / 250 \\
			\bottomrule
		\end{tabular}
	}
	\caption{Amount of iterations required to reach a prescribed tolerance for the Stokes problem \eqref{eq:problem_stokes}.}
	\label{tab:stokes}
\end{table}

As initial geometry, we choose $\Omega\iidx{0} = (-3, 6) \times (-2, 2) \setminus \Omega\subobs\iidx{0}$, with $\Omega\subobs\iidx{0}$ being a circle with center $(0, 0)$ and radius \num{0.5}, in analogy to \cite{Schulz2016Computational}. This is discretized by a non-uniform mesh consisting of 6525 vertices and 12326 triangles, of which 620 form the deformable boundary $\Gamma\subobs$. For the numerical solution of the state and adjoint systems we use a mixed finite element method with piecewise quadratic Lagrange elements for the velocity component and piecewise linear Lagrange elements for the pressure component, which is LBB-stable for the saddle point structure of \eqref{eq:pde_stokes}. 
For $\velocity\subin$ we choose the parabolic profile 
\begin{equation*}
	\velocity\subin(x) = \nicefrac{1}{4} \left( 2 - x_2 \right) \left( 2 + x_2 \right),
\end{equation*}
so that the maximum inlet velocity is \num{1}. The solution of the Stokes system \eqref{eq:pde_stokes} on the initial and optimized geometries, obtained by the L-BFGS~5 method, can be seen in Figures~\ref{fig:state_stokes_velocity} and~\ref{fig:state_stokes_pressure}, where the magnitude of the velocity and the pressure are depicted, respectively. The weights $\nu_i$ for the regularization are chosen as $\nu_1 = \num{1e4}$ and $\nu_2 = \num{1e2}$. For the parameter $\lamesecond$ needed for the bilinear form $a_\Omega$ we follow the approach described in \cite{Schulz2016Computational} and compute it as the solution of the following Laplace problem
\begin{equation}
	\label{eq:lame_siebenborn}
	\begin{alignedat}{2}
		- \laplace \lamesecond &= 0 \quad &&\text{ in } \Omega, \\
		\lamesecond &= \mu_\textrm{max} \quad &&\text{ on } \Gamma\subobs, \\
		\lamesecond &= \mu_\textrm{min} \quad &&\text{ on } \Gamma\subin \cup \Gamma\subwall \cup \Gamma\subout,
	\end{alignedat}
\end{equation}
where we choose $\mu_\textrm{max} = \num{500}$ and $\mu_\textrm{min} = \num{1}$ as in \cite{Schulz2016Computational}. The remaining parameters for Algorithm~\ref{algo:transformed} are shown in Table~\ref{tab:parameters_stokes}.

\begin{figure}[t]
	\centering
	\begin{subfigure}[t]{0.325\textwidth}
		\includegraphics[width=\textwidth]{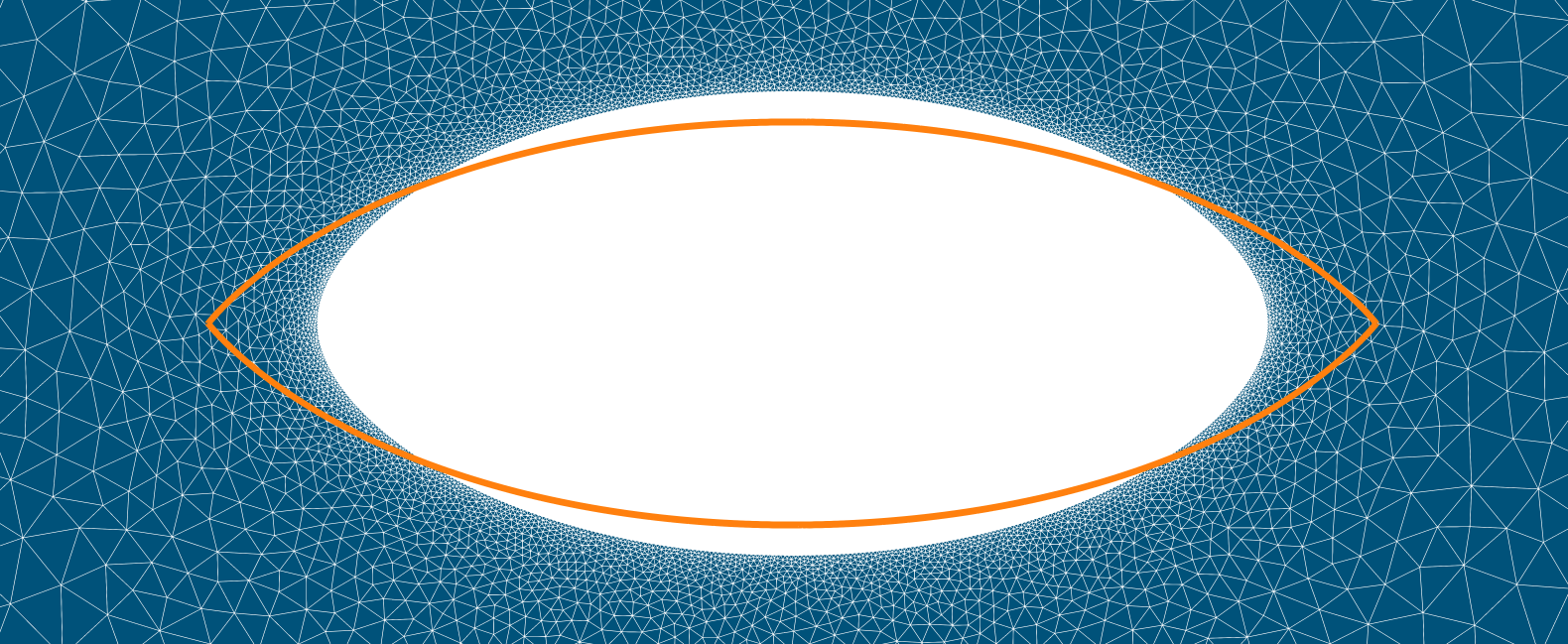}
		\caption{Gradient descent method.}
	\end{subfigure}
	\hfil
	\begin{subfigure}[t]{0.325\textwidth}
		\includegraphics[width=\textwidth]{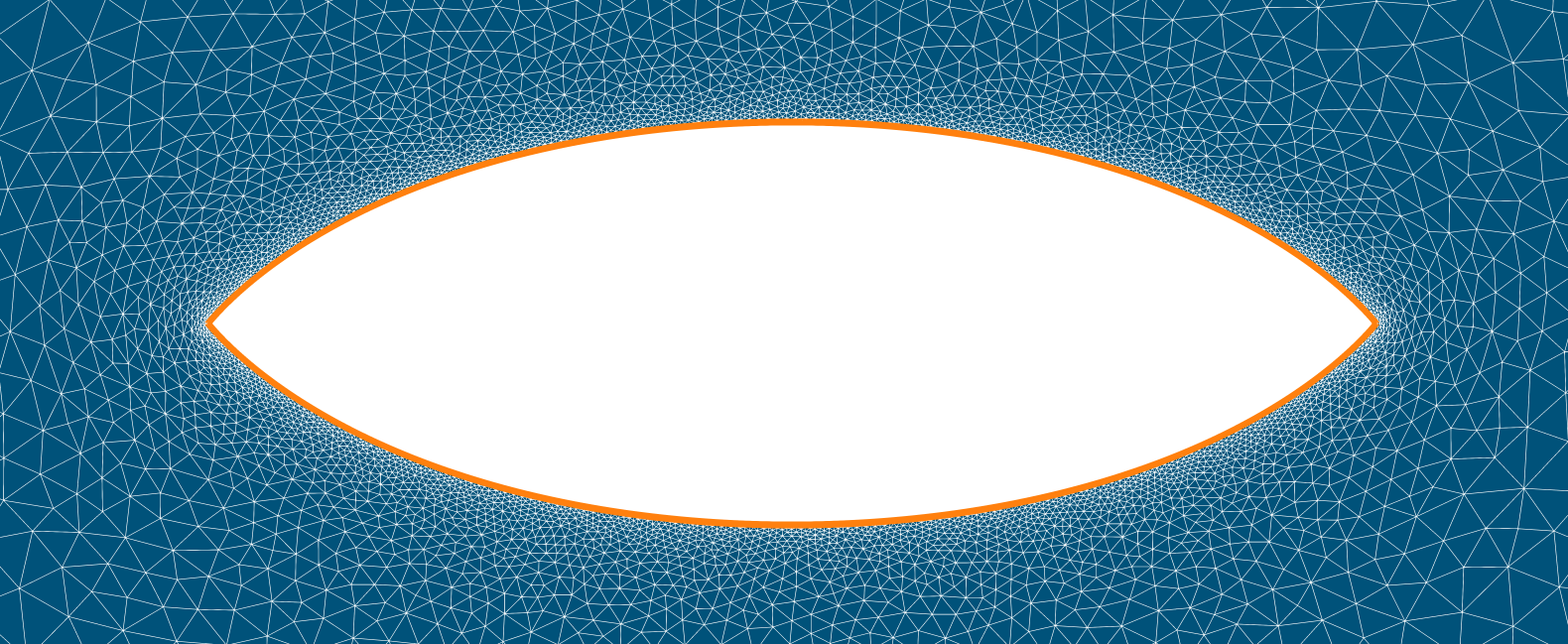}
		\caption{Polak-Ribi\`ere NCG method.}
	\end{subfigure}
	\hfil
	\begin{subfigure}[t]{0.325\textwidth}
		\includegraphics[width=\textwidth]{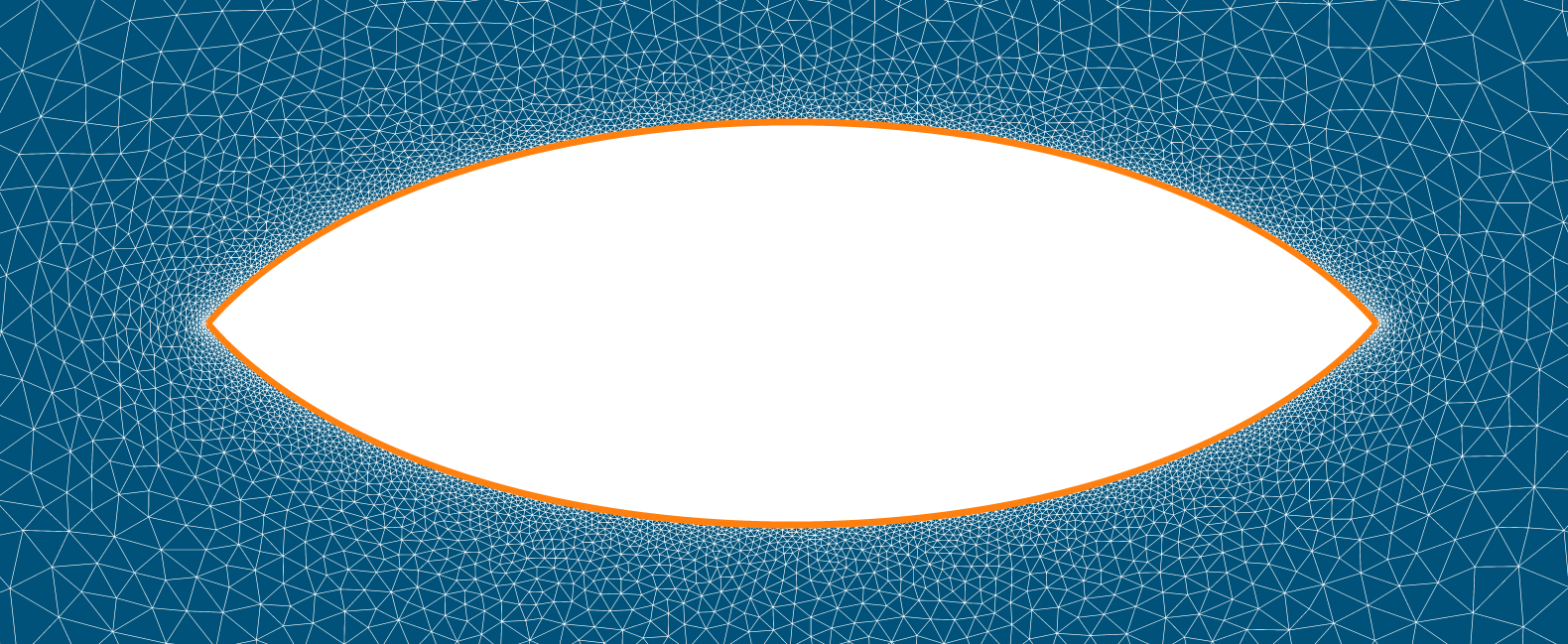}
		\caption{Dai-Yuan NCG method.}
	\end{subfigure}
	\vskip\baselineskip
	\begin{subfigure}[t]{0.49\textwidth}
		\frame{\includegraphics[width=\textwidth]{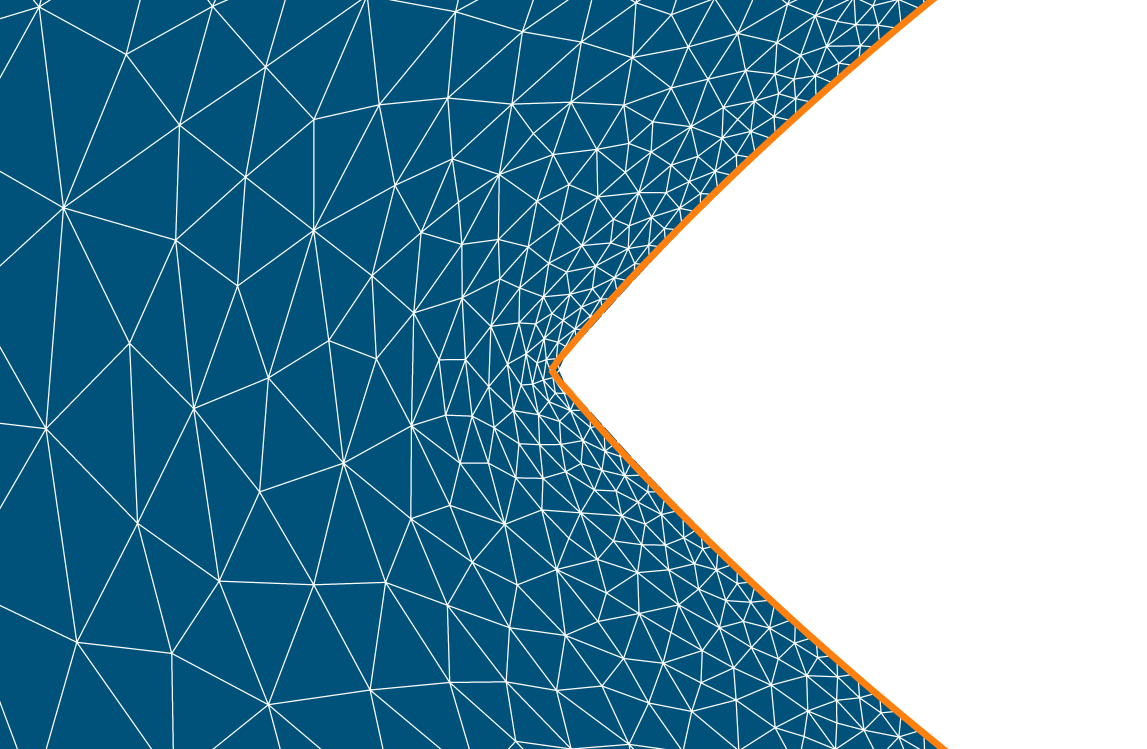}}
		\caption{Polak-Ribi\`ere NCG method.}
		\label{sfig:PR}
	\end{subfigure}
	\hfil
	\begin{subfigure}[t]{0.49\textwidth}
		\frame{\includegraphics[width=\textwidth]{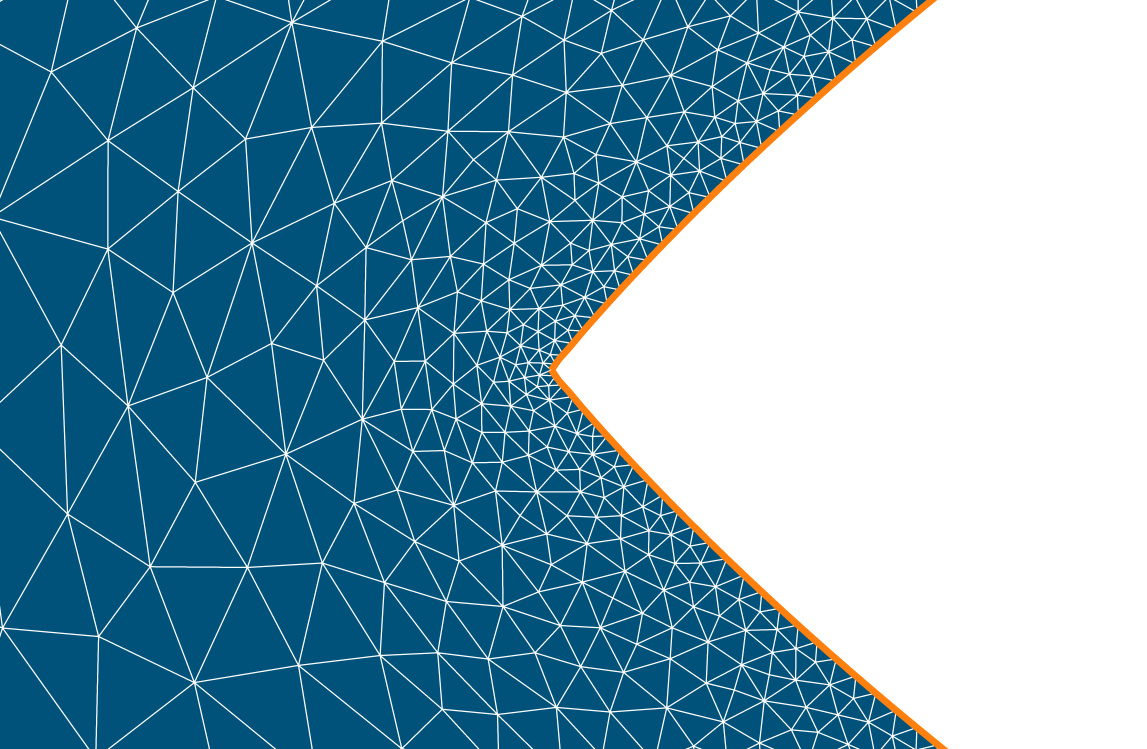}}
		\caption{Dai-Yuan NCG method.}
	\end{subfigure}
	\caption{Optimized Shapes $\Omega\subobs$ (white) and $\Omega$ (blue) compared to the solution of the L-BFGS~5 method (orange) for the Stokes problem \eqref{eq:problem_stokes}.}
	\label{fig:optimized_stokes}
\end{figure}

As before, the history and performance of the methods are shown in Figure~\ref{fig:history_stokes} and Table~\ref{tab:stokes} (cf. Section~\ref{ssec:implementation}). Again, we observe that the NCG methods perform well for this problem. Comparing them to the gradient descent method, we see that all of them give significantly better results for this problem as the former fails to reach even the tolerance of \num{1e-1}, whereas the NCG methods reach much lower tolerances. In particular, all NCG methods except for the Polak-Ribi\`ere and Hager-Zhang variants reach the desired tolerance of \num{5e-4}. Moreover, we see that the NCG methods perform similarly to the L-BFGS methods. In particular, the Dai-Yuan NCG method, which performs nearly as good as the L-BFGS~5 method, works very well. Moreover, the Hestenes-Stiefel NCG method yields results that are very similar to the L-BFGS~1 method, but performs a bit weaker than the L-BFGS~3 and~5 ones. As before, we note that the L-BFGS methods need less solves for the state system due to the built-in scaling of the search directions.

The optimized geometries, depicted in Figure~\ref{fig:optimized_stokes} for the gradient descent, Polak-Ribi\`ere, and Dai-Yuan methods, confirm our previous findings. Whereas the geometry is still far from the optimal ogive for the gradient descent method, we observe basically no difference anymore between the solution of the L-BFGS~5 and the Polak-Ribi\`ere method, even though the latter performed worst of all NCG methods. We see that the overall shape of the ogive is very well approximated by the Polak-Ribi\`ere method, and that there are only very subtle differences occurring at the front and back, which are barely even visible
(cf. Figure~\ref{sfig:PR}). Further, we observe that the Dai-Yuan method approximates the ogive perfectly, including the front and back wedges, and yields basically identical results to the L-BFGS~5 method.

\subsection{Shape Optimization of a Pipe}
\label{ssec:pipe}

Finally, we investigate the shape optimization of a pipe, based on the problems considered in \cite{Schmidt2010Efficient, Ham2019Automated}. Similarly to the previous problem, our objective is to minimize the dissipated energy of the flow, which is now governed by the incompressible Navier-Stokes equations. We also have a geometric constraint for this problem, namely we fix the volume of the pipe so that the geometry cannot degenerate. We denote the domain of the pipe by $\Omega$ and its boundary by $\Gamma$. The latter is divided into three parts, namely the inlet $\Gamma\subin$, the wall boundary $\Gamma\subwall$, and the outlet $\Gamma\subout$. We assume that the inlet and outlet are fixed, and that only a subset $\Gamma\subdef \subset \Gamma\subwall$ is deformable, whereas the part $\Gamma\subfix = \Gamma\subwall \setminus \Gamma\subdef$, located near the in- and outlet, is fixed, as in \cite{Schmidt2010Efficient}. This is done to avoid a degeneration of the geometry near the in- and outlet sections of the pipe.

\begin{figure}[b]
	\centering
	\includegraphics[width=0.75\textwidth]{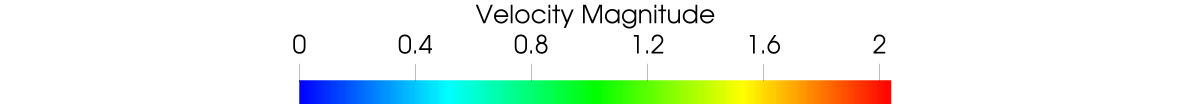}
	\begin{subfigure}{0.49\textwidth}
		\centering
		\includegraphics[width=0.95\textwidth]{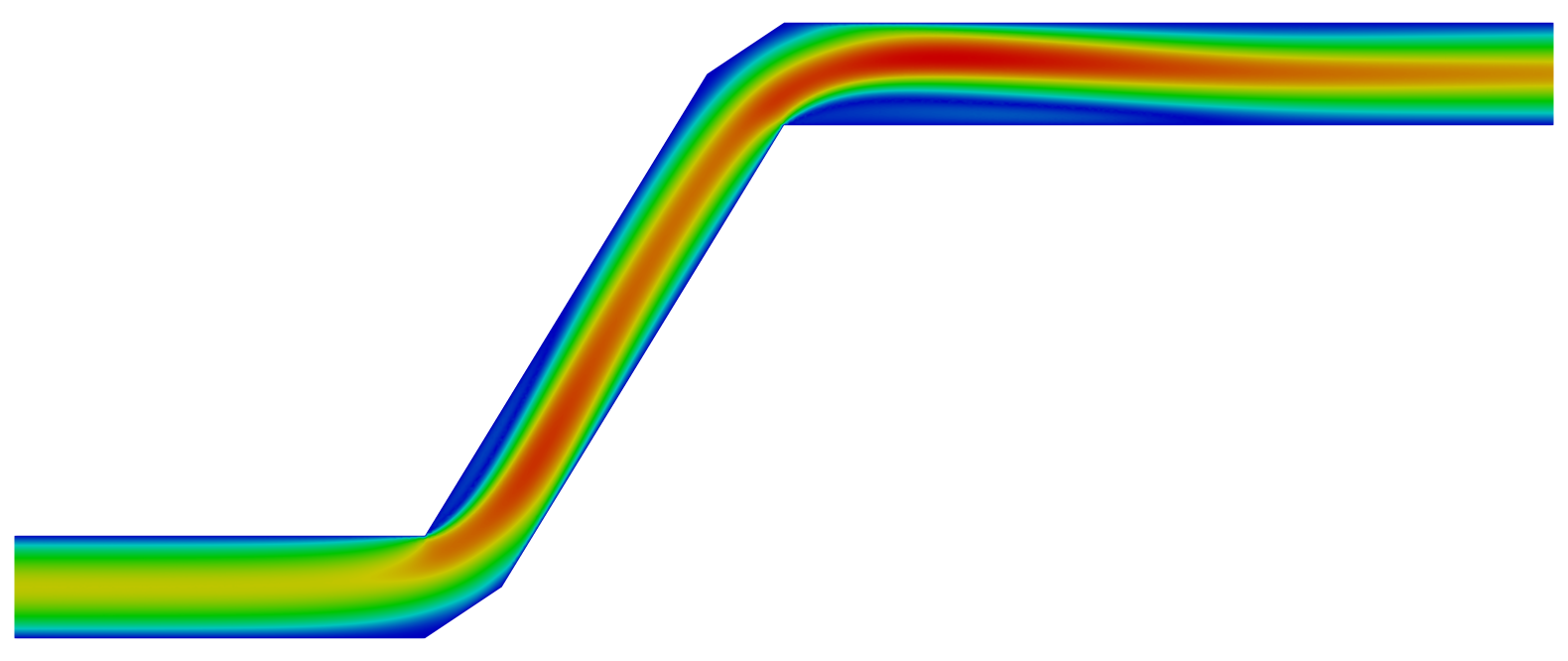}
		\caption{Initial geometry.}
	\end{subfigure}
	\hfil
	\begin{subfigure}{0.49\textwidth}
		\centering
		\includegraphics[width=0.95\textwidth]{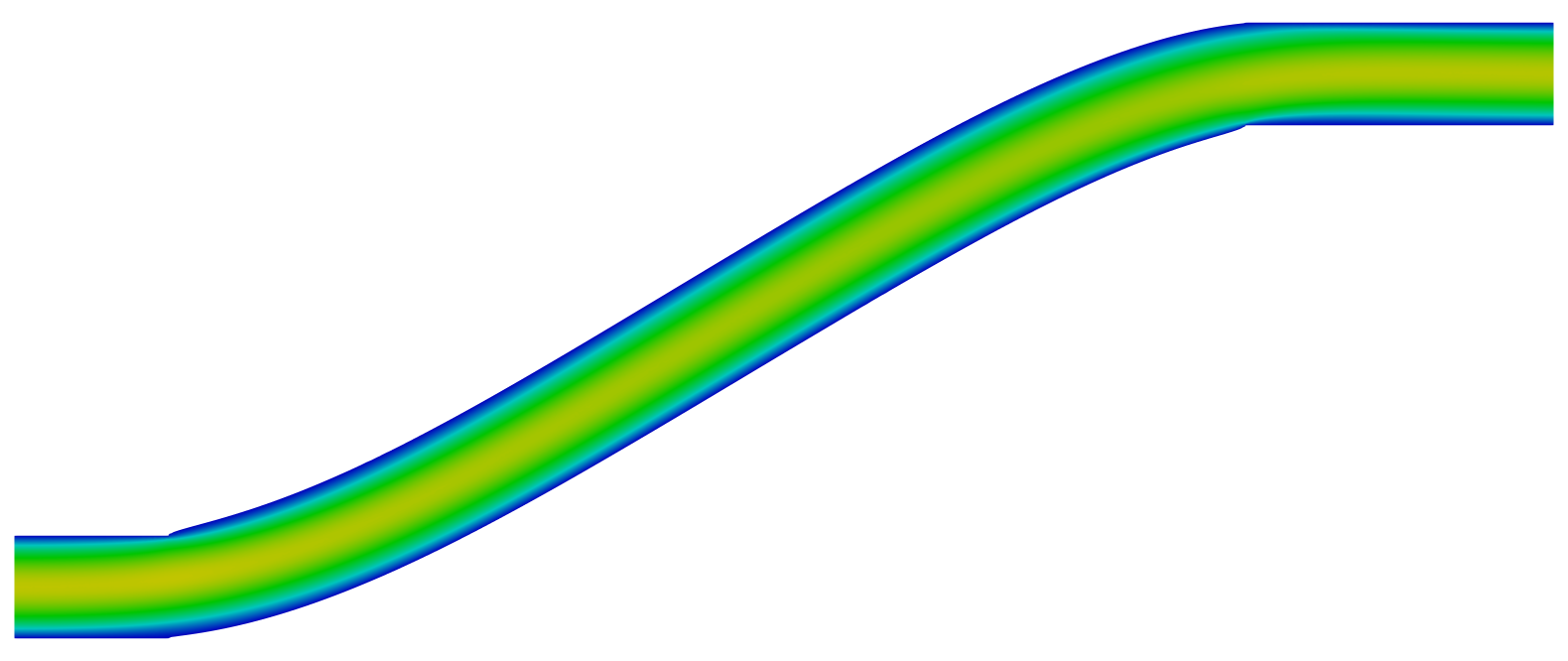}
		\caption{Optimized Geometry.}
	\end{subfigure}
	\caption{Magnitude of the velocity $\norm{\velocity}{}$ for the Navier-Stokes problem \eqref{eq:pde_pipe} on the initial and optimized geometries, obtained by the L-BFGS~5 method.}
	\label{fig:state_navier_stokes_velocity}
\end{figure}

\begin{figure}[b]
	\centering
	\includegraphics[width=0.75\textwidth]{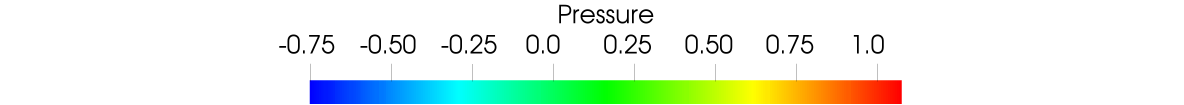}
	\begin{subfigure}{0.49\textwidth}
		\centering
		\includegraphics[width=0.95\textwidth]{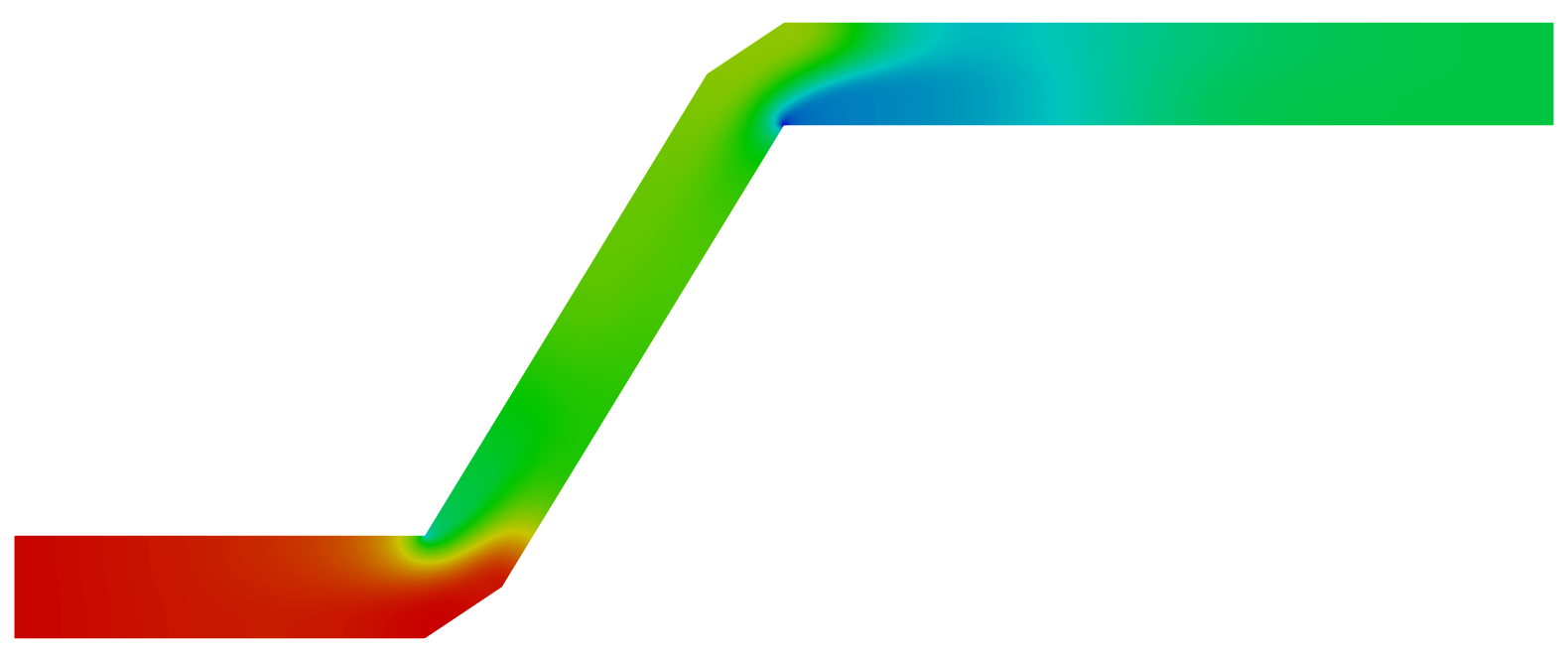}
		\caption{Initial geometry.}
	\end{subfigure}
	\hfil
	\begin{subfigure}{0.49\textwidth}
		\centering
		\includegraphics[width=0.95\textwidth]{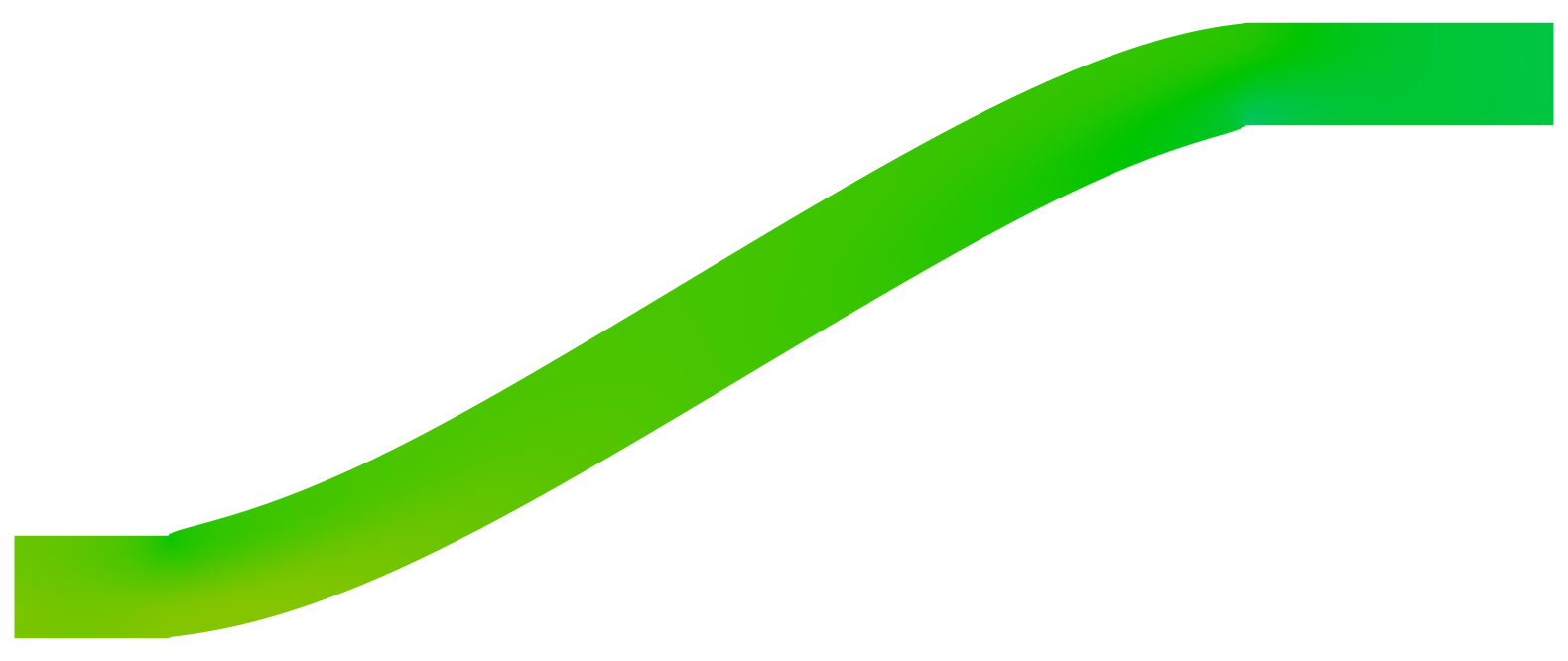}
		\caption{Optimized Geometry.}
	\end{subfigure}
	\caption{Pressure $\pressure$ for the Navier-Stokes problem \eqref{eq:pde_pipe} on the initial and optimized geometries, obtained by the L-BFGS~5 method.}
	\label{fig:state_navier_stokes_pressure}
\end{figure}

The non-dimensionalized incompressible Navier-Stokes system describing the flow for this problem reads
\begin{equation}
	\label{eq:pde_pipe}
	\begin{alignedat}{2}
		- \frac{1}{\reynolds} \laplace \velocity + \left( \velocity \cdot \grad \right) \velocity + \grad \pressure &= 0 \quad &&\text{ in } \Omega,\\
		\divergence{\velocity} &= 0 \quad &&\text{ in } \Omega, \\
		\velocity &= \velocity\subin \quad &&\text{ on } \Gamma\subin, \\
		\velocity &= 0 \quad &&\text{ on } \Gamma\subwall, \\
		\frac{1}{\reynolds} \partial_\normal \velocity - \pressure\normal &= 0 \quad &&\text{ on } \Gamma\subout,
	\end{alignedat}	
\end{equation}
where $\reynolds > 0$ denotes the Reynolds number. The corresponding shape optimization problem reads
\begin{equation}
	\label{eq:problem_pipe}
	\min_{\Omega \in \admissiblegeom}\ \costfunction(\Omega, \velocity) = \integral{\Omega} \nicefrac{1}{\reynolds} \norm{\grad \velocity}{F}^2 \dmeas{x} + \nicefrac{\nu}{2} \left( \text{vol}(\Omega) - \text{vol}(\Omega\iidx{0}) \right)^2 \qquad \text{ s.t. } \eqref{eq:pde_pipe},
\end{equation}
where we again have regularized the geometric constraints in the cost functional, and $\Omega\iidx{0}$ denotes the initial geometry of the pipe with corresponding boundary $\Gamma\iidx{0}$. The set of admissible geometries for this problem is given by
\begin{equation*}
	\admissiblegeom = \Set{\Omega \subset \R^d | \Omega \subset \holdall,\ \Gamma\subin = \Gamma\subin\iidx{0},\ \Gamma\subout = \Gamma\subout\iidx{0}, \text{ and } \Gamma\subfix = \Gamma\subfix\iidx{0}}.
\end{equation*}
For the corresponding shape derivative and adjoint system we refer the reader to \cite{Schmidt2010Efficient}.

\begin{table}[t]
	\centering
	{\footnotesize
		\rowcolors{2}{\tablegray}{white}
		\setlength{\tabcolsep}{1em}
		\begin{tabular}{l r}
			\toprule
			parameter & value \\
			\midrule
			initial step size $t\iidx{0}$ & 5e-3 \\
			%
			%
			maximum number of iterations $k_{\text{max}}$ & 50 \\
			iterations for NCG restart $\cgiter$ & $\infty$ \\
			tolerance for NCG restart $\cgtol$ & 0.25 \\
			\midrule
			first Lam\'e parameter $\lamefirst$ & 0.0 \\
			second Lam\'e parameter  $\lamesecond$ & 1.0 \\
			damping parameter $\damping$ & 0.0\\
			\bottomrule
		\end{tabular}
		\caption{Parameters for Algorithm~\ref{algo:transformed} for the Navier-Stokes problem \eqref{eq:problem_pipe}.}
		\label{tab:parameters_pipe}
	}
\end{table}

\begin{figure}[b]
	\centering
	\begin{subfigure}{0.49\textwidth}
		\includegraphics[width=\textwidth]{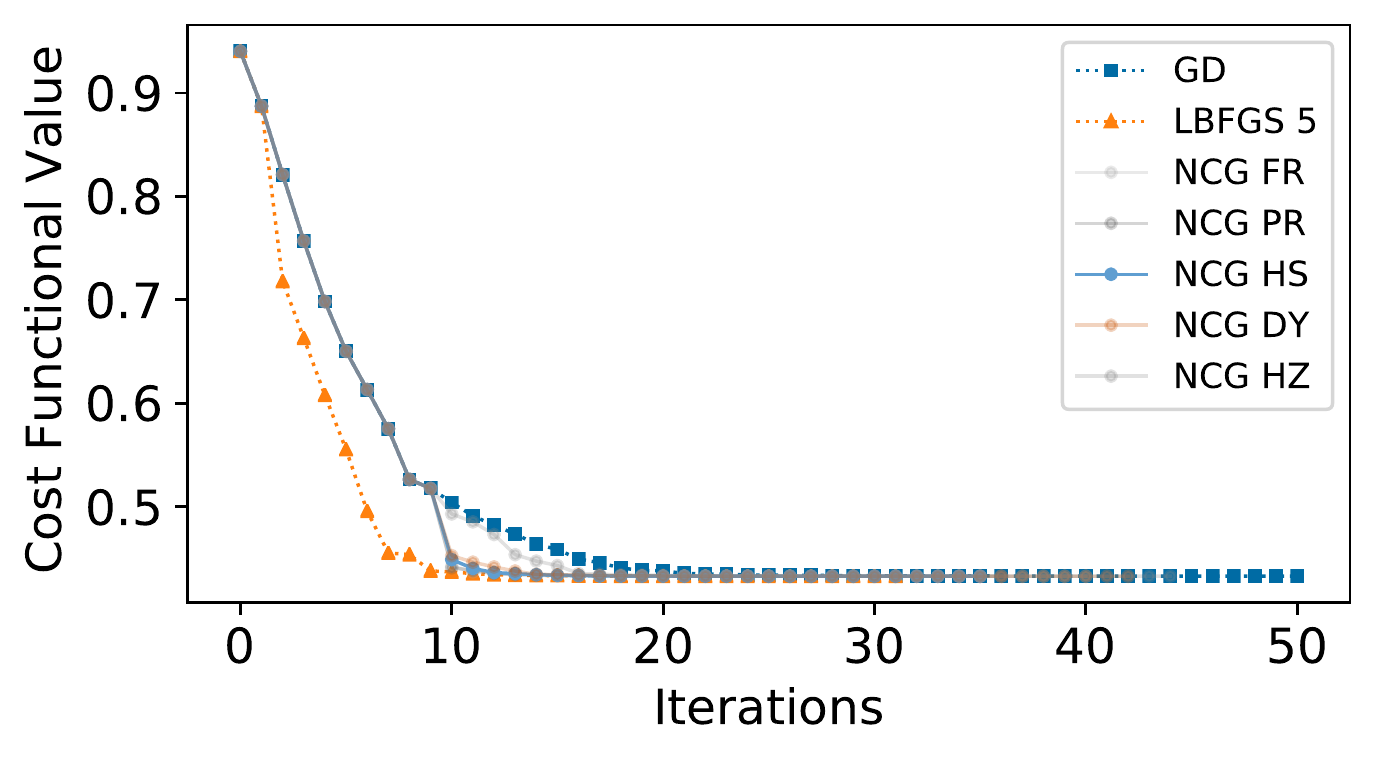}
		\caption{History of the cost functional.}
	\end{subfigure}
	\hfil
	\begin{subfigure}{0.49\textwidth}
		\includegraphics[width=\textwidth]{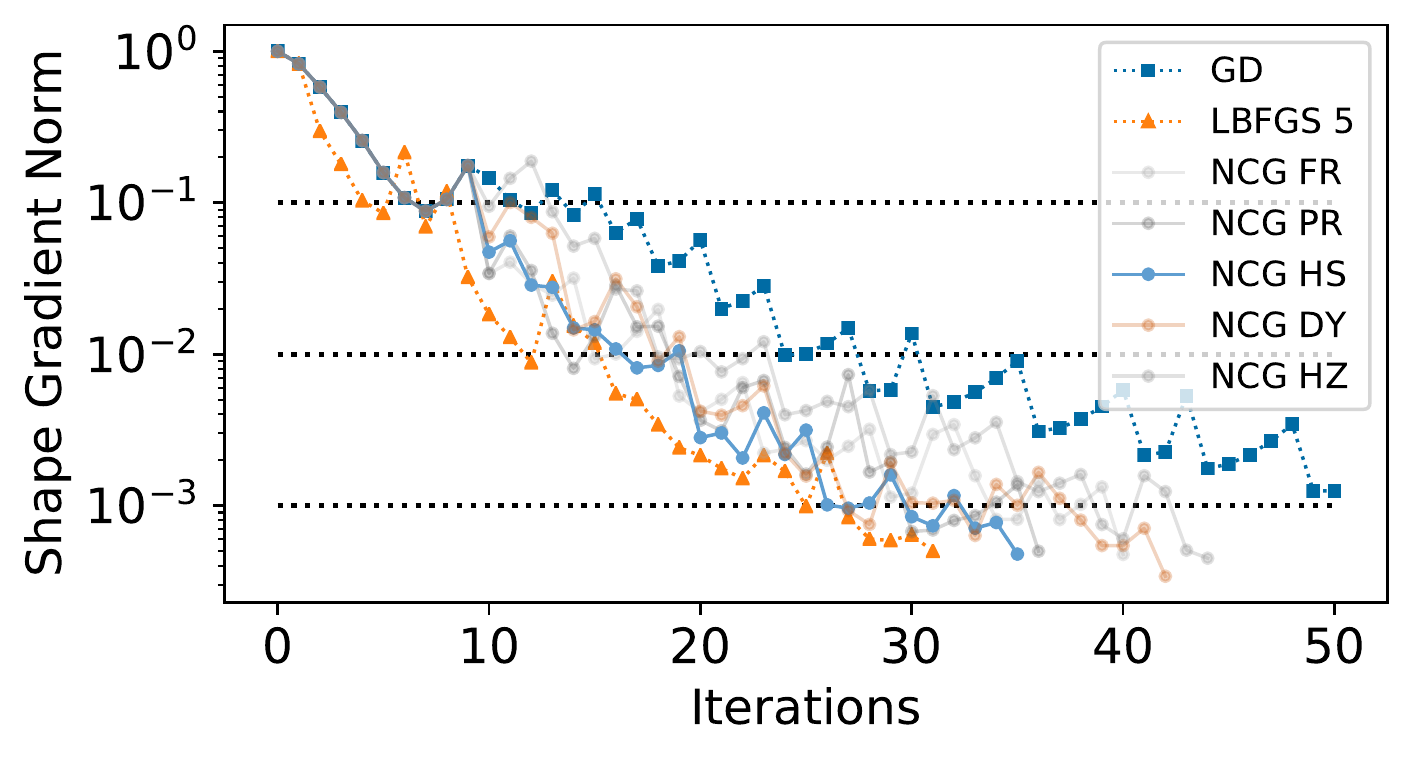}
		\caption{History of the relative shape gradient norm.}
	\end{subfigure}
	\caption{History of the optimization algorithms for the Navier-Stokes problem \eqref{eq:problem_pipe}.}
	\label{fig:history_pipe}
\end{figure}

\begin{table}[b]
	\centering
	{\footnotesize
		\rowcolors{2}{\tablegray}{white}
		\setlength{\tabcolsep}{1em}
		\begin{tabular}{c r r r r r r c c}
			\toprule
			tol & \num{1e-01} & \num{5e-02} & \num{1e-02} & \num{5e-03} & \num{1e-03} & \num{5e-04}  &  & state / adjoint solves \\
			\midrule
			GD & 7 & 18 & 24 & 31 & - & - &  & 94 / 50 \\
			\midrule
			L-BFGS 1 & 5 & 8 & 11 & 13 & 27 & 35 &  & 38 / 36 \\
			L-BFGS 3 & 5 & 10 & 17 & 17 & 23 & 32 &  & 35 / 33 \\
			L-BFGS 5 & 5 & 9 & 12 & 18 & 25 & 31 &  & 37 / 32 \\
			\midrule
			CG FR & 7 & 10 & 15 & 20 & 34 & 40 &  & 69 / 41 \\
			CG PR & 7 & 10 & 14 & 20 & 30 & 36 &  & 57 / 37 \\
			CG HS & 7 & 10 & 17 & 20 & 27 & 35 &  & 58 / 36 \\
			CG DY & 7 & 14 & 18 & 20 & 27 & 42 &  & 69 / 43 \\
			CG HZ & 7 & 16 & 18 & 24 & 39 & 44 &  & 78 / 45 \\
			\bottomrule
		\end{tabular}
	}
	\caption{Amount of iterations required to reach a prescribed tolerance for the Navier-Stokes problem \eqref{eq:problem_pipe}.}
	\label{tab:pipe}
\end{table}

The magnitude of the velocity $\norm{\velocity}{}$ and the pressure $\pressure$ on the initial and optimized geometries are shown in Figures~\ref{fig:state_navier_stokes_velocity} and~\ref{fig:state_navier_stokes_pressure}, respectively. Note, that the initial geometry is based on the one used in \cite{Schmidt2010Efficient}, in particular, it has a boundary with several kinks, which induce strong forces on the flow. This geometry is discretized with a mesh consisting of 16652 nodes and 32300 triangles. As in \cite{Schmidt2010Efficient}, the Reynolds number is chosen as $\reynolds = \num{400}$ and the inlet velocity is given by a parabolic profile with mean inlet velocity of \num{1}. For the discretization of the Navier-Stokes problem \eqref{eq:pde_pipe} we proceed analogously to the previous section and use a mixed finite element method with piecewise quadratic Lagrange elements for the velocity component and piecewise linear Lagrange elements for the pressure component due to their LBB-stability for the underlying saddle point structure. The weight $\nu$ for the regularization of the volume constraint is chosen as $\nu = 1$. The remaining parameters for Algorithm~\ref{algo:transformed} can be found in Table~\ref{tab:parameters_pipe}.

The history of the optimization and the performance of the methods are shown in Figure~\ref{fig:history_pipe} and Table~\ref{tab:pipe}, as remarked in Section~\ref{ssec:implementation}. As before, we observe that the NCG methods are very efficient at solving the shape optimization problem. They perform better than the gradient descent method, having lower function values and gradient norms throughout the optimization as well as needing less iterations to reach the investigated tolerances. In particular, the gradient descent method only reaches a tolerance of \num{5e-3}, whereas all NCG methods reach the desired tolerance of \num{5e-4}. In contrast to the previous problems, we now also restart the NCG methods via the criterion \eqref{eq:restart_ncg}. In particular, we observe that due to the restarting procedure, the first nine iterates of the gradient descent and the NCG methods coincide. This is due to the fact that the geometry changes a lot in the first couple of iterations which makes it hard for the NCG methods to generate (nearly) orthogonal gradient deformations. After these initial iterations, changes of the geometry become smaller and the NCG methods do not have to be restarted as often, which leads to a superior convergence speed. Comparing the NCG methods with the L-BFGS ones we observe that they again yield quite comparable results, with the L-BFGS methods needing just slightly fewer iterations to reach a specified tolerance. In particular, the Polak-Ribi\`ere and Hestenes-Stiefel methods perform very similarly to the L-BFGS~1 method, and the L-BFGS~3 and~5 methods only need 3 to 4 iterations less to reach a tolerance of \num{5e-4}. So, overall, there are not many differences between the performance of the NCG and the L-BFGS methods for this problem, except for the fact that the L-BFGS methods need less state equation solves due to the built-in scaling of the search direction.

\begin{figure}[t]
	\centering
	\begin{subfigure}[t]{0.49\textwidth}
		\includegraphics[width=\textwidth]{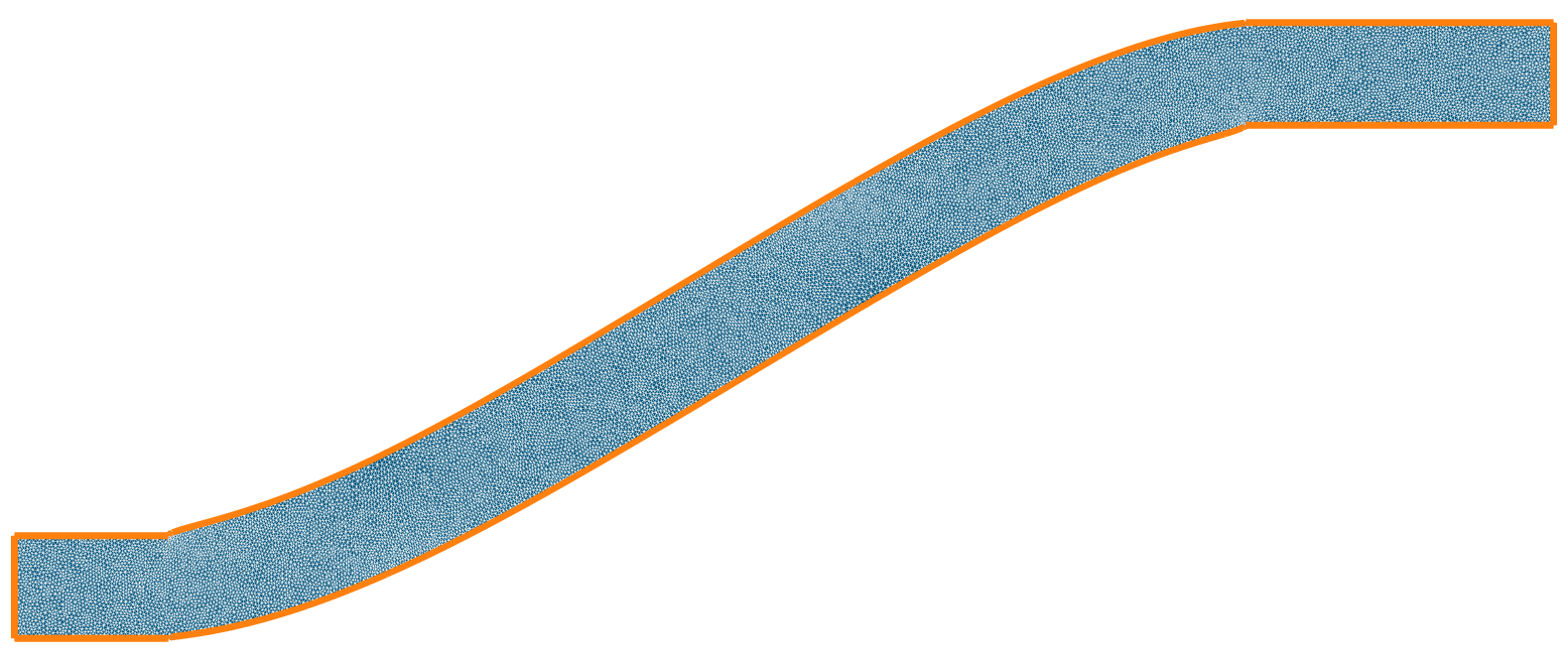}
		\caption{Gradient descent method.}
	\end{subfigure}
	\hfil
	\begin{subfigure}[t]{0.49\textwidth}
		\includegraphics[width=\textwidth]{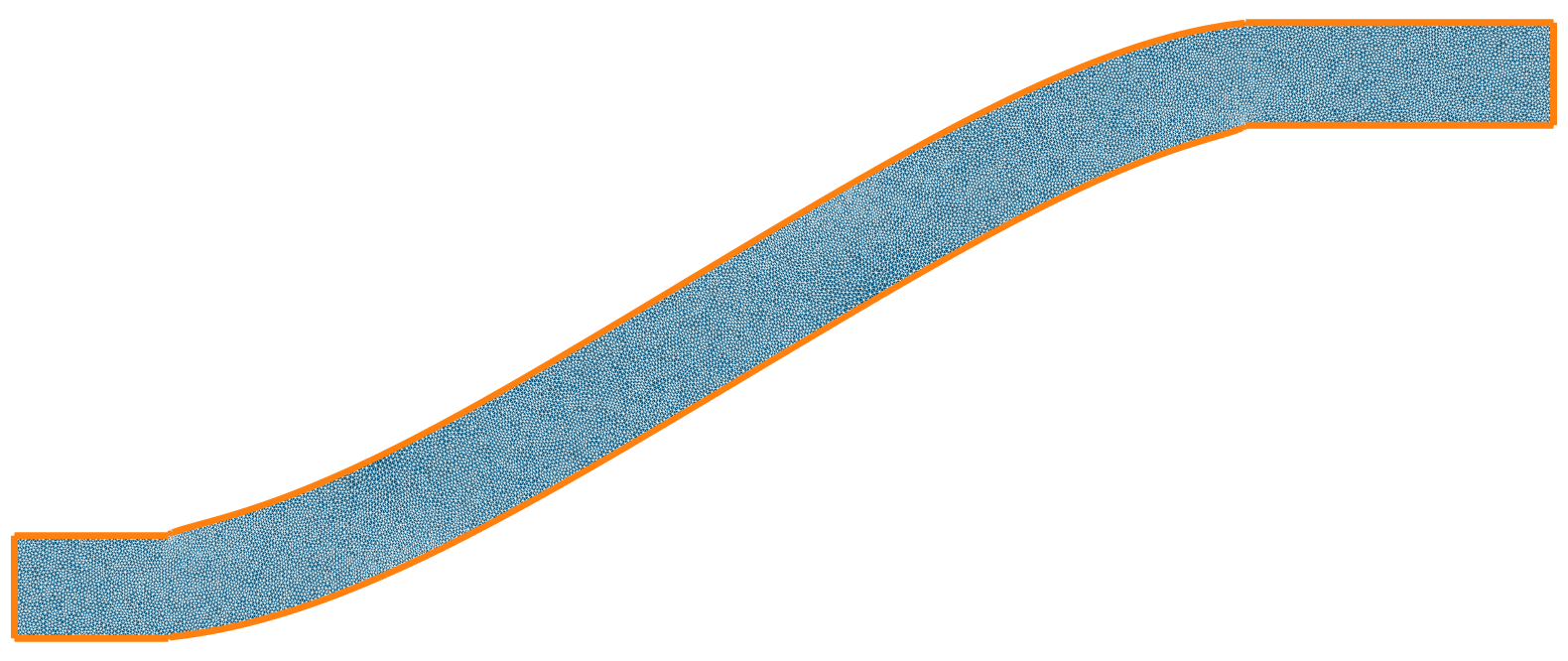}
		\caption{Hager-Zhang NCG method.}
	\end{subfigure}
	\caption{Optimized Shapes (blue) compared to the solution of the L-BFGS~5 method (orange) for the Navier-Stokes problem \eqref{eq:problem_pipe}.}
	\label{fig:optimized_pipe}
\end{figure}

Finally, in Figure~\ref{fig:optimized_pipe} the optimized geometries obtained by the gradient descent and Hager-Zhang NCG methods are compared to the one from the L-BFGS~5 method. There, we do not observe any visual difference between the geometries even for the gradient descent method. Moreover, we note that all methods were able to handle the kinks in the initial geometry well since the middle sections of the optimized pipes are smooth. Note, that the new kinks in the geometry arise at the points where the deformable and fixed boundary are joint together and could potentially be avoided by using a curvature or surface regularization.

\subsection{Summary of the Numerical Comparison}

Our numerical results for the proposed NCG methods suggest that they perform very well for the numerical solution of shape optimization problems. In particular, each of the NCG methods yields significantly better results than the gradient descent method for all considered test cases. In fact, the NCG methods perform similarly to the L-BFGS methods and are about as efficient as the L-BFGS~1 method. Note, that even though the L-BFGS~3 and~5 methods yield the best results for almost all test cases, this increased efficiency comes at the price of a considerably increased memory usage, which can become prohibitive for very large-scale problems (cf. \cite{Kelley1999Iterative}). Hence, in such settings the L-BFGS~1 and NCG methods are very attractive as they have a significantly better performance than the gradient descent method, while using only slightly more memory. In conclusion, the NCG methods proposed in this paper are an efficient and attractive addition to already established gradient-based shape optimization algorithms, and they are particularly interesting for very large-scale problems such as the ones arising from industrial applications.

\section{Conclusion and Outlook}
\label{sec:conclusion}

In this paper, we have proposed and investigated nonlinear conjugate gradient (NCG) methods for shape optimization. After recalling shape calculus and the Steklov-Poincar\'e metrics from \cite{Schulz2016Efficient} for Riemannian shape optimization, we presented a general algorithmic framework for the solution of shape optimization problems. We formulated novel NCG methods for shape optimization in the context of our algorithmic framework and detailed its numerical discretization. Finally, we investigated the proposed NCG methods numerically on four benchmark problems and compared their performance to the gradient descent and L-BFGS methods. The results of this investigation show that the NCG methods significantly outperform the gradient descent method, while needing only slightly more memory, and that they are comparable to the L-BFGS methods. This makes them an efficient and attractive addition to gradient-based optimization methods for the numerical solution of shape optimization problems, particularly for large industrial problems where memory is an issue.

For future research, a theoretical analysis of NCG methods on infinite-dimensional manifolds, particularly for shape optimization, is of great interest, as this is not well-developed yet (cf. \cite{Ring2012Optimization}). Moreover, the application and investigation of the proposed NCG methods to other shape optimization problems, in particular large-scale industrial problems, is also of interest.

\section*{Acknowledgments}


The author thanks Ren\'e Pinnau and Christian Leith\"auser for helpful discussions and comments, and thankfully acknowledges financial support from the Fraunhofer Institute for Industrial Mathematics ITWM.

\bibliographystyle{siamplain}
\bibliography{lit}

\end{document}